%% file: gln.tex
\documentclass[a4paper,11pt,french]{article} 
        \usepackage[french,english]{babel}
\usepackage[T1]{fontenc}
\usepackage[latin1]{inputenc}
\usepackage[margin=1in]{geometry}
        \usepackage{amsfonts}
        \usepackage{amsmath}
        \usepackage{amssymb}
        \usepackage{amsthm}
        \usepackage{amstext}
        \usepackage{mathtools}
        \usepackage{enumerate}
        \usepackage{exscale}
        \usepackage{relsize}
        \usepackage{hyperref}
        \usepackage[all]{xy}
\usepackage[all]{xy}

\include{abbrvs}

\usepackage{latexsym}
\usepackage{amsfonts}
\usepackage{amssymb}
\usepackage{amstext}
\usepackage{amsmath}

\usepackage[normalem]{ulem}
\usepackage{zref-perpage}

\zmakeperpage{footnote}

  \usepackage{hyperref}
\usepackage[dvipsone]{bookmark}
        \usepackage[all]{xy}
        \usepackage{paralist}
                \usepackage{dsfont}
                \usepackage{bbm}
                \usepackage{upgreek}
                \usepackage{mathrsfs}

\numberwithin{equation}{section}
\begin{document}
\selectlanguage{french}
\title{La variante infinitésimale 
de la formule des traces de Jacquet-Rallis pour
les groupes linéaires}
\date{}
\author{Micha\l \ Zydor
\footnote{Université Paris Diderot
Institut de Mathématiques de Jussieu-Paris Rive Gauche
UMR7586
Bâtiment Sophie Germain
Case 7012
75205 PARIS Cedex 13
France 
mail: \texttt{michal.zydor@imj-prg.fr}
}}
\maketitle 
\selectlanguage{english}
\begin{abstract} 
We establish an infinitesimal version of the 
Jacquet-Rallis trace formula for general linear groups. 
Our formula is obtained by integrating a kernel
truncated à la Arthur multiplied by the 
absolute value of the determinant to the power $s \in \C$.
It has a geometric side which is a 
sum of distributions $I_{\ol}(s, \cdot)$ indexed by 
the invariants of the adjoint action of 
$\Gl_{n}(\rmF)$ on $\gll_{n+1}(\rmF)$
as well as a "spectral side" 
consisting of the Fourier transforms 
of the aforementioned distributions. 
 We prove that the distributions $I_{\ol}(s, \cdot)$ 
are invariant and depend only on the choice of 
the Haar measure on $\Gl_n(\A)$. 
For regular semi-simple classes $\ol$, $I_{\ol}(s, \cdot)$ is 
a relative orbital integral of Jacquet-Rallis. 
For classes $\ol$ called relatively regular semi-simple, we express $I_{\ol}(s, \cdot)$
in terms of relative orbital integrals regularised by means of zeta functions.
\end{abstract}
\selectlanguage{french}
\begin{abstract}
Nous établissons une variante infinitésimale 
de la formule des traces de Jacquet-Rallis 
pour les groupes linéaires. 
Notre formule s'obtient par 
intégration d'un noyau tronqué 
à la Arthur multiplié par la valeur absolue du déterminant à la puissance $s \in \C$. 
Elle possède un côté géométrique 
qui est une somme de distributions $I_{\ol}(s, \cdot)$
indexée par les invariants de l'action de 
$\Gl_{n}(\rmF)$ sur $\gll_{n+1}(\rmF)$
ainsi qu'un "côté spectral" 
formé des transformées de Fourier 
des distributions précédentes. 
On démontre que les distributions $I_{\ol}(s, \cdot)$ 
sont invariantes et ne dépendent que du choix 
de la mesure de Haar sur $\Gl_{n}(\A)$. 
Pour des classes $\ol$ semi-simples régulières $I_{\ol}(s, \cdot)$ 
est une intégrale
orbitale relative de Jacquet-Rallis.
Pour les classes $\ol$ dites relativement semi-simples régulières, 
on exprime $I_{\ol}(s, \cdot)$
en terme des intégrales orbitales relatives régularisées à l'aide des fonctions zêta. 
\end{abstract}

\tableofcontents
\setcounter{section}{-1}

\input{intro}
\input{prelegomenes}
\input{general_gln}
\input{quantitative_gln}
\input{fourier_gln}
\input{semisimple_orbs}
\input{koSansT_gln}

\bibliographystyle{alpha}
\bibliography{bibliography}
\end{document}

%% file: abbrvs.tex
\theoremstyle{plain} 
\newtheorem{theo}{Th\'eor\'eme}[section]
\newtheorem{cor}[theo]{Corollaire}
\newtheorem{prop}[theo]{Proposition}
\newtheorem{lem}[theo]{Lemme}

\theoremstyle{definition} 
\newtheorem{defi}[theo]{Définition}

\newtheorem{rem}[theo]{Remarque}

\newcommand{\beqyn}{\begin{eqnarray}}
\newcommand{\enqyn}{\end{eqnarray}}
\newcommand{\blem}{\begin{lem}}
\newcommand{\elem}{\end{lem}}
\newcommand{\brop}{\begin{prop}}
\newcommand{\erop}{\end{prop}}
\newcommand{\bcor}{\begin{cor}}
\newcommand{\ecor}{\end{cor}}
\newcommand{\brem}{\begin{rem}}
\newcommand{\erem}{\end{rem}}
\newcommand{\brems}{\begin{rems}}
\newcommand{\erems}{\end{rems}}
\newcommand{\benum}{\begin{enumerate}}
\newcommand{\enum}{\end{enumerate}}
\newcommand{\nident}{\noindent}

\newcommand{\btheo}{\begin{theo}}
\newcommand{\etheo}{\end{theo}}
\newcommand{\bdem}{\begin{proof}}
\newcommand{\edem}{\end{proof}}
\newcommand{\bdefi}{\begin{defi}}
\newcommand{\edefi}{\end{defi}}

\newcommand{\rmE}{\mathrm{E}}
\newcommand{\rmF}{\mathrm{F}}

\newcommand{\rmI}{\mathrm{I}}

\newcommand{\rmK}{\mathrm{K}}

\newcommand{\rmO}{\mathrm{O}}

\newcommand{\rmV}{\mathrm{V}}

\newcommand{\rmw}{\mathrm{w}}

\newcommand{\Hom}{\mathrm{Hom}}

\newcommand{\Lie}{\mathrm{Lie}}

\newcommand{\Ad}{\mathrm{Ad}}

\newcommand{\Tr}{\mathrm{Tr}}

\newcommand{\Gl}{\mathrm{GL}}

\newcommand{\Res}{\mathrm{Res}}

\newcommand{\vol}{\mathrm{vol}}
\newcommand{\Rel}{\mathrm{Re}}

\newcommand{\ind}{\mathds{1}}

\newcommand{\upla}{\uprho}
\newcommand{\uphi}{\upphi}
\newcommand{\al}{\alpha}

\newcommand{\la}{\lambda}
\newcommand{\La}{\Lambda}

\newcommand{\brtau}{\bar \tau}
\newcommand{\brnl}{\bar \nl}

\newcommand{\inv}{^{-1}}

\newcommand{\bilif}{\langle \cdot,\cdot \rangle}

\newcommand{\matx}[4]{\begin{pmatrix} #1 & #2 \\ #3 & #4 \\ \end{pmatrix}}
\newcommand{\dsl}{\displaystyle \left(}
\newcommand{\rb}{\right)}

\newcommand{\bsl}{\backslash}

\newcommand{\bs}{$\Box$}

\newcommand{\rar}{\rightarrow}

\newcommand{\hrar}{\hookrightarrow}
\newcommand{\smin}{\smallsetminus}
\newcommand{\sbs}{\subseteq}
\newcommand{\sbn}{\subsetneq}
\newcommand{\nsbs}{\nsubseteq}
\newcommand{\sps}{\supseteq}

\newcommand{\gll}{\mathfrak{gl}}

\newcommand{\all}{\mathfrak{a}}

\newcommand{\fl}{\mathfrak{l}}
\newcommand{\gl}{\mathfrak{g}}
\newcommand{\hl}{\mathfrak{h}}
\newcommand{\ml}{\mathfrak{m}}
\newcommand{\nl}{\mathfrak{n}}

\newcommand{\sn}{\mathfrak{s}}
\newcommand{\tl}{\mathfrak{t}}
\newcommand{\ul}{\mathfrak{u}}

\newcommand{\ol}{\mathfrak{o}}

\newcommand{\N}{\mathbb{N}}
\newcommand{\Z}{\mathbb{Z}}
\newcommand{\A}{\mathbb{A}}
\newcommand{\C}{\mathbb{C}}

\newcommand{\Q}{\mathbb{Q}}
\newcommand{\R}{\mathbb{R}}

\newcommand{\Gm}{\mathbb{G}_m}
\newcommand{\Ga}{\mathbb{G}_a}

\newcommand{\calD}{\mathcal{D}}

\newcommand{\calF}{\mathcal{F}}

\newcommand{\calI}{\mathcal{I}}
\newcommand{\calJ}{\mathcal{J}}
\newcommand{\calK}{\mathcal{K}}

\newcommand{\calO}{\mathcal{O}}
\newcommand{\calP}{\mathcal{P}}

\newcommand{\calS}{\mathcal{S}}

\newcommand{\calV}{\mathcal{V}}

\newcommand{\calZ}{\mathcal{Z}}

\newcommand{\cad}{c'est-\`a-dire }

\newcommand{\dem}{$\text{\it D\'emonstration.}$ }

\newcommand{\nin}{\notin}

\newcommand{\tlvpi}{\widetilde{\varpi}}

\newcommand{\tlB}{\widetilde{B}}

\newcommand{\tlG}{\widetilde{G}}

\newcommand{\tlI}{\widetilde{I}}

\newcommand{\tlK}{\widetilde{K}}

\newcommand{\tlP}{\widetilde{P}}
\newcommand{\tlR}{\widetilde{R}}
\newcommand{\tlQ}{\widetilde{Q}}
\newcommand{\tlS}{\widetilde{S}}

\newcommand{\tla}{\widetilde{a}}

\newcommand{\tlk}{\widetilde{k}}

\newcommand{\tls}{\widetilde{s}}

\newcommand{\tlul}{\widetilde{\ul}}

\newcommand{\tlgl}{\widetilde{\gl}}

\newcommand{\tlzero}{\widetilde{0}}
\newcommand{\tlone}{\widetilde{1}}
\newcommand{\tltwo}{\widetilde{2}}

\newcommand{\htau}{\hat \tau}
\newcommand{\hf}{\hat f}

\newcommand{\hDelta}{\widehat \Delta}

\newcommand{\brg}{\bar g}

\newcommand{\brP}{\overline{P}}

\newcommand{\brX}{\overline{X}}

\newcommand{\brLa}{\overline{\La}}

\newcommand{\allts}{\all_{\tlzero}^{*}}
\newcommand{\tlBmin}{\tlB_{0}}

\newcommand{\relPb}{\calF(M_{\tlzero}, B)}

\newcommand{\FFi}{\rmF_{I} \times \rmF_{-I}}

\newcommand{\IoB}{I_{0}}

\newcommand{\sbse}{\subseteq_{\epsilon}}

\newcommand{\acalI}{|\mathcal{I}|}
\newcommand{\acalK}{|\mathcal{K}|}

\newcommand{\acalJ}{|\mathcal{J}|}

%% file: intro.tex
\section{Introduction}

\subsection{Contexte}

Dans \cite{jacqrall}, Jacquet et Rallis proposent une approche 
via une formule des traces relative 
à la conjecture globale de Gan-Gross-Prasad (GGP) pour les groupes unitaires \cite{ggp}. 
Ils proposent une existence de deux formules des traces relatives, une pour les groupes unitaires 
et une pour les groupes linéaires. Les deux formules ont des côtés géométriques et spectraux. 
Le comparaison de côtés géométriques de ces deux formules doit mener aux résultats spectraux 
qui entament la conjecture de Gan-Gross-Prasad. Une version simple de cette formule, \cad 
une version valable pour une certaine classe de fonctions test, a été utilisé par Zhang 
pour démontrer une partie substantielle de la conjecture GGP \cite{zhang2, zhang1} ainsi que 
certains cas du raffinement de la 
 conjecture GGP du à Ichino et Ikeda \cite{ichinoIkeda} et N. Harris \cite{harris}.
 
 \subsection{Nos résultats}
 
Afin qu'on puisse étendre les résultat de Zhang il faut des formules des traces valables 
 pour toutes les fonctions lisses à support compact. 
Cet article étudie la variante infinitésimale du côté géométrique de la formule des traces relative de Jacquet-Rallis pour les groupes linéaires
et il est précédé par \cite{leMoi} où on étudie la variante infinitésimale du côté géométrique pour les groupes unitaires.
Les côtés géométriques pour les groupes ainsi que les côtés spectraux seront présentés dans \cite{leMoi3}.
  
Soient $\rmE/\rmF$ une extension quadratique 
de corps de nombres, $\A$ l'anneau des adèles de $\rmF$ et 
$\eta : \rmF^{*} \bsl \A^{*} \rar \C$ le caractère associé à l'extension $\rmE/\rmF$ 
par la théorie de corps de classes. Soit $n \in \N$ et notons $G = \Gl_{n}$
que l'on voit comme un sous-groupe de $\tlG = \Gl_{n+1}$ par le plongement 
diagonal $g \mapsto \matx{g}{}{}{1}$.
Soit $S_{n+1}$ la sous-$\rmF$-variété de $\Res_{\rmE/\rmF} \tlG$ composée 
de $g \in \Res_{\rmE/\rmF} \tlG$ tels que $g \brg = 1$ où $\brg$ c'est le conjugué de $g$
par l'élément non-trivial du groupe de Galois de $\rmE/\rmF$. 
On voit $G$ et $\tlG$ comme des sous-$\rmF$-groupes de $\Res_{\rmE/\rmF} \tlG$. 
L'action de $G$ sur $\Res_{\rmE/\rmF} \tlG$ stabilise alors $S_{n+1}$.

Du côté géométrique pour le groupe linéaire de la formule des traces relative de Jacquet-Rallis
on étudie l'intégrale formelle 
\begin{equation}\label{eq:intGpsG}
\int_{G(\rmF) \bsl G(\A)} k_{f}(x)\eta (\det x)dx
\end{equation}
où $f \in C_{c}^{\infty}(S_{n+1}(\A))$ est une fonction lisse à support compact sur $S_{n+1}(\A)$ 
et $k_{f}(x)$ égale $\sum_{\gamma \in S_{n+1}(\rmF} f(x^{-1} \gamma x)$ pour $x \in G(\rmF) \bsl G(\A)$.
L'intégrale n'est pas bien définie pour toutes les fonctions $f$. 
Zhang utilise cet intégrale pour des fonctions vérifiant certaines 
conditions du support locale pour lesquelles elle converge et admet une décomposition en une somme des intégrales 
orbitales relatives. 

Dans cet article on s'intéresse à la version infinitésimale de l'intégrale \eqref{eq:intGpsG}. 
Soient $\tlgl = \Lie(\tlG)$ et $\sn_{n+1}$ l'espace tangent à l'identité de $S_{n+1}$. 
Alors, $\sn_{n+1}$, vu comme un sous-espace de $\Res_{\rmE/\rmF}\tlgl$, c'est:
\[
\sn_{n+1} = \{X \in \Res_{\rmE/\rmF}\tlgl_{n+1}| X + \brX = 0\}.
\]
Si on choisit un $\tau \in \rmE$ tel que $\overline{\tau} = -\tau$, la multiplication 
par $\tau$ induit un isomorphisme entre $\sn_{n+1}$ et $\tlgl$.
Le groupe $G$ agit sur $\sn_{n+1}$ et sur $\tlgl$ par adjonction. Cette identification 
est alors $G$-équivariante. 
Soit $s \in \C$. L'analogue infinitésimal de l'intégrale \eqref{eq:intGpsG} que l'on 
considère c'est
\[
I(s,f) = \int_{G(\rmF) \bsl G(\A)}k_{f}(x) |\det x|_{\A}^{s}\eta (\det x)dx, \quad k_{f}(x) = \sum_{\xi \in \tlgl(\rmF)}f(x^{-1} \xi x), \ 
f \in \calS(\tlgl(\A)),
\]
où $|\cdot|_{\A}$ c'est la valeur absolue standard sur le groupe des idèles de $\A$ 
et $\calS(\tlgl(\A))$ c'est l'espace des fonctions de type Bruhat-Schwartz sur $\tlgl(\A)$.  
Dans cet article on définit une régularisation de cette intégrale par un processus de troncature à la Arthur 
et on étudie les propriétés de la distribution ainsi obtenue.
Décrivons brièvement cette troncature maintenant. 

Fixons $B$ un sous-groupe de Borel 
de $G$ ainsi que sa décomposition de Levi $B= M_{0}N_{0}$ avec $M_{0}$ une partie 
de Levi de $B$ et $N_{0}$ son radical unipotent. 
On note $\calF(B)$ l'ensemble des sous-groupes paraboliques de $G$ contenant $B$. 
Tout $P \in \calF(B)$
admet alors 
une unique décomposition de Levi $P = M_{P}N_{P}$. 
Soit $M_{\tlzero}$ l'unique sous-groupe de Levi minimal de $\tlG$ contenant $M_{0}$. 
On note $\relPb$ l'ensemble des sous-groupes paraboliques de $\tlG$ contenant $B$. 
Tout $\tlP \in \relPb$ admet alors une unique décomposition de Levi $\tlP = M_{\tlP}N_{\tlP}$ 
où $M_{0} \sbs M_{\tlP}$. 
Pour tout $\tlP \in \relPb$ on a $\tlP \cap G \in \calF(B)$. On note alors $P := \tlP \cap G$. 


Au début de la section \ref{sec:RTFtlgl} 
on introduit une 
décomposition de $\tlgl(\rmF)$ en classes géométriques, 
dont on note l'ensemble $\calO$, 
qui sont stables par action adjointe de $G(\rmF)$.
Pour tous $\tlP \in \relPb$, $\ol \in \calO$
et $f \in \calS(\tlgl(\A))$
on pose
\[
k_{f,\tlP,\mathfrak{o}}(x) =  \sum_{\xi \in 
\Lie(M_{\tlP})(\rmF) \cap \mathfrak{o}}
\int_{\Lie(N_{\tlP})(\A)}f(x\inv (\xi + U)x)dU, \ 
x \in M_{\tlP}(\rmF)N_{\tlP}(\A)\backslash \tlG(\A).
\]
Le noyau tronqué est défini alors comme
\[
k^{T}_{f,\mathfrak{o}}(x) = 
\sum_{\tlP \in \relPb} 
(-1)^{d_{\tlP}- d_{\tlG}}\sum_{\delta \in P(\rmF)\backslash G(\rmF)}
\htau_{\tlP}(H_{\tlP}(\delta x) - T)k_{f,\tlP,\mathfrak{o}}(\delta x), \ 
x \in \tlG(\rmF)\backslash \tlG(\A),
\]
où, si l'on pose $\all_{\tlP} := 
\Hom_{\Z}(\Hom_{\rmF}(M_{\tlP}, \Gm), \R)$, alors 
$H_{\tlP} : \tlG(\A) \rar \all_{\tlP}$ 
c'est l'application de Harish-Chandra,
$\htau_{\tlP}$ 
est la fonction caractéristique
d'un cône obtus dans $\all_{\tlP}$, 
$d_{\tlP} = \dim_{\R} \all_{\tlP}$ et 
$T$ est un paramètre dans $\all_{\tlzero} := \Hom_{\Z}(\Hom_{\rmF}(M_{\tlzero}, \Gm), \R)$
(voir le paragraphe \ref{par:prelimstraceG}). 
Notre premier résultat, démontré dans le paragraphe \ref{par:convergenceG}, 
c'est alors:
\begin{theo}[cf. théorème \ref{thm:MainConvG}]\label{thm:cvgGIntro}
Pour tout
 $T \in T_{+} + \mathfrak{a}_{\tlzero}^{+}$ et tout $\sigma \in \R$ on a:
\[
\sum_{\mathfrak{o} \in \mathcal{O}}
\int_{G(\rmF) \bsl G(\A)}|k_{f,\mathfrak{o}}^{T}(x)| |\det x|_{\A}^{\sigma}dx < \infty.
\]
\end{theo} 

Ensuite, on s'intéresse au comportement de l'application 
$T \mapsto \int_{G(\rmF) \bsl G(\A)}k_{f,\mathfrak{o}}^{T}(x) |\det x|_{\A}^{s} \eta (\det x)dx$. 
Dans le paragraphe \ref{par:asymptChaptG} on obtient:
\begin{theo}[cf. théorème \ref{thm:mainQualitThm}]\label{thm:polExpIntro}
Soit $\ol \in \calO$ et $s \in \C$. La fonction
\[
T \mapsto\  I^{T}_{\ol}(s, f) := 
\int_{G(\rmF) \bsl G(\A)} k_{f,\mathfrak{o}}^{T}(x) |\det x|_{\A}^{s} \eta(\det x) dx,
\] 
où $T$ parcourt $T_{+} + \mathfrak{a}_{\tlzero}^{+}$,
est un polynôme-exponentielle
dont la partie purement polynomiale est constante si $s \neq -1,1$.
\end{theo}

On note $I_{\ol}(s, f)$ la partie constante du 
polynôme-exponentielle $I_{\ol}^{T}(s,f)$ pour $s \in \C \smin \{-1,1\}$. 
Il s'avère que la distribution
$I_{\ol}(s, \cdot)$ a des propriétés remarquables. 
On obtient:
\btheo\label{thm:invIntro} Soient $f \in \calS(\tlgl(\A))$, $\ol \in \calO$ et $s \in \C \smin \{-1,1\}$.
\begin{enumerate}
\item[(cf. théorème \ref{invarianceTheoG}).] Pour $y \in G(\A)$, 
notons $\upphi^{y} \in \calS(\tlgl(\A))$ définie 
par $\upphi^{y}(X) = \upphi(\Ad(y)X)$. On a alors
$$I_{\ol}(s,f^{y}) = |\det y|_{\A}^{s} \eta(\det y) I_{\ol}(s,f).$$ 
\item[(cf. paragraphe \ref{par:noChoixMadeG}).]
La distribution $I_{\ol}(s, \cdot)$ ne dépendent 
que des choix des mesures de Haar.  
\end{enumerate}
\etheo

Dans le paragraphe \ref{par:regOrbsChapG}, on constate que 
pour les classes $\ol$ dites semi-simples régulières, 
sur lesquelles $G(\rmF)$ opère fidèlement et transitivement, 
la distribution $I_{\ol}^{T}(s, \cdot )$ ne dépend pas de $T$ et $I_{\ol}^{T}(s,f)$ égale
une intégrale orbitale relative qui 
apparaît déjà dans \cite{jacqrall}.

Dans la section \ref{sec:FourierTrans} on démontre la 
version infinitésimale de la formule des traces relative pour les groupes linéaires
\btheo[cf. théorème \ref{thm:RTFJRIG}]  
Pour tout $f \in \calS(\tlgl(\A))$ et tout $s \in \C \smin \{-1,1\}$ on a
\[
\sum_{\mathfrak{o} \in \mathcal{O}}I_{\mathfrak{o}}(s,f) = 
\sum_{\mathfrak{o} \in \mathcal{O}}I_{\mathfrak{o}}(s,\hat f).
\]
\etheo
\noindent
Ici $\hat f$ 
c'est une transformée de Fourier (on en considère plusieurs)
de $f$. 

Notre dernier résultat, démontré dans la section 
\ref{sec:orbRrssSG}, ce sont des formules explicites pour certaines 
nouvelles distributions $I_{\ol}(s, \cdot)$. 
Tout $X \in \tlgl(\rmF)$ s'écrit comme $\matx{B}{u}{v}{d}$ où $B \in \Lie(G)(\rmF)$, $u,v \in \rmF^{n}$ et $d \in \rmF$. 
On écrit alors $B_{X} := B$.
Soient $P \in \calF(B)$ et $B \in \Lie(M_{P})(\rmF)$ un élément elliptique 
(qui n'est contenu dans aucun sous-algèbre parabolique propre de $\Lie(M_{P})$).
Soit $\ol \in \calO$ une classe 
contenant un élément $X$ tel que $B_{X} = B$.
Alors, $\ol$ est une réunion finie d'orbites pour l'action de $G(\rmF)$. 
En plus, les orbites de dimension maximale dans $\ol$ ont des centralisateurs triviaux. 
Fixons un représentant $X_{0}$ de l'orbite fermée dans $\ol$ tel que 
$B_{X_{0}} = B$ et choisissons un ensemble $\rmO_{\ol} \sbs \ol$
de représentants des orbites de dimension maximale dans $\ol$ 
de sorte que tout $X \in \rmO_{\ol}$ vérifie $B_{X} = B$.
Soit $T_{0}$ le centralisateur de $X_0$ dans $G$.
Pour $f \in \calS(\tlgl(\A))$ et $X \in \rmO_{\ol}$ 
on pose:
\[
\zeta_{X}(f)(\la) = 
\int_{G(\A)}
f(\Ad(x^{-1})X)e^{\la(H_{P}(x))} \eta (\det x)dx,
\]
où $\la \in \Hom_{\rmF}(T_{0}, \Gm) \otimes_{\Z} \C =: \all_{\C}^{*}$ que l'on peut voir naturellement
comme un sous-espace de $\Hom_{\R}(\all_{P}, \C)$.
L'intégrale définissant $\zeta_{X}(f)$
converge sur un ouvert
de $\all_{\C}^{*}$ qui dépend de $X$
et admet un prolongement méromorphe 
à $\all_{\C}^{*}$ noté aussi $\zeta_{X}(f)$. 
Notons que
le déterminant est naturellement un élément de $\Hom_{\rmF}(T_{0}, \Gm)$, on a donc 
un élément associé $\det \in \all_{\C}^{*}$. On obtient alors:
\btheo[cf. théorème \ref{thm:theThmOrbsG}] 
Soit $\ol$ comme ci-dessus. 
Alors, la droite $s \det \in \all_{\C}^{*}$, où $s \in \C$, 
n'est contenu dans aucun hyperplan singulier de la 
fonction méromorphe 
\[
\all_{\C}^{*} \ni \la \mapsto \dsl \sum_{X \in \rmO_{\ol}} \zeta_{X}(f) \rb (\la).
\]
En plus, la fonction méromorphe
$
\C \ni s \mapsto \dsl \sum_{X \in \rmO_{\ol}} \zeta_{X}(f) \rb (s \det)
$
est holomorphe pour $s \neq -1,1$ et l'on a:
\[
I_{\ol}(s, f) = \dsl \sum_{X \in \rmO_{\ol}} \zeta_{X}(f) \rb (s \det), \quad 
\forall \ s \in \C.
\]
\etheo
\noindent 
Dans la section 5 de \cite{leMoi} on obtient un résultat complètement analogue pour les groupes unitaires. 

Commentons finalement l'apparition du déterminant à la puissance complexe dans la définition de la distribution 
$I_{\ol}(s, \cdot)$. On ajoute ce terme en vu de possibles applications 
aux dérivées des intégrales orbitales relatives de Jacquet-Rallis. 
Pour les classes semi-simples régulières, les dérivées des analogues locaux des distributions 
$I_{\ol}(s, \cdot)$ définies sur la variété $S_{n+1}$ ont d'abord été introduites dans 
 \cite{zhang3} et ensuite étudiées dans \cite{zhangRapTers, zhangRapSmith} 
 pour des fonctions test particulières. 
Ces travaux lient ces dérivées à des nombres d'intersection sur certains espaces de Rapoport-Zink 
et s'inscrivent dans le cadre de la variante arithmétique de la conjecture de Gan-Gross-Prasad.
 Récemment, dans \cite{mihatsch}, 
 les dérivées des analogues locaux des distributions infinitésimales 
$I_{\ol}(s, \cdot)$ ont aussi été étudiées. On espère que notre formule 
va avoir un intérêt dans ces questions.

\textbf{Remerciements}. 
Je voudrais remercier mon directeur de thèse, Pierre-Henri Chaudouard, pour toute son aide. 

%% file: prelegomenes.tex
\section{Prolégomènes}\label{sec:prolego}

\subsection{Préliminaires pour la formule des traces}\label{par:prelimstraceG}

Soient $\rmF$ un corps de nombres et 
$G$ un $\rmF$-groupe algébrique 
réductif que l'on suppose déployé sur $\rmF$. 
Pour tout sous-$\rmF$-groupe de Levi $M$ 
de $G$ (\cad un facteur de Levi d'un $\rmF$-sous-groupe parabolique 
de $G$) soit $\calF(M)$ l'ensemble de $\rmF$-sous-groupes 
paraboliques de $G$ contenant $M$ et $\calP(M)$ le sous-ensemble 
de $\calF(M)$ composé de sous-groupes paraboliques 
admettant $M$ comme facteur de Levi.
On fixe 
un sous-groupe de Levi minimal $M_{0}$ 
de $G$. 
On appelle les éléments de $\calF(M_{0})$ 
les sous-groupes paraboliques semi-standards 
et les éléments de $\calP(M_{0})$ les sous-groupes de Borel.
On utilisera toujours le symbole $P$, avec des indices éventuellement, 
pour noter un sous-groupe parabolique semi-standard.
Pour tout $P \in \calF(M_{0})$ soit $N_{P}$ le radical unipotent 
de $P$ et $M_{P}$ le facteur de Levi de $P$ contenant $M_{0}$. 
On a alors $P = M_{P}N_{P}$. On note $A_{P}$ 
le tore central de $M_{P}$ déployé sur $\rmF$ 
et maximal pour cette propriété. Si $P$ est un sous-groupe de Borel, 
on a alors $A_{P} = M_{0}$ et on pose dans ce cas $A_{0} := A_{P} = M_{P}$.
Pour $P_{1} \in \calF(M_{0})$, quand il n'y aura pas d'ambiguïté, on écrit 
$N_{1}$ au lieu de $N_{P_{1}}$, $M_{1}$ au lieu de $M_{P_{1}}$ etc.

Soit  $P \in \calF(M_{0})$.
On définit le $\R$-espace vectoriel 
$\mathfrak{a}_{P} := \Hom_{\Z}(\Hom_{\rmF}(M_{P}, \Gm), \R)$, 
isomorphe à $\Hom_{\Z}(\Hom_{\rmF}(A_{P}, \Gm), \R)$ grâce à l'inclusion $A_{P} \hrar M_{P}$, 
ainsi que son espace dual $\all_{P}^{*} = \Hom_{\rmF}(M_{P}, \Gm) \otimes_{\Z} \R$ et on pose
 \begin{equation}\label{eq:dPDefG}
 d_{P} = \dim_{\R} \all_{P}, \quad d_{Q}^{P} = d_{Q} - d_{P}, \ 
 Q \sbs P.
\end{equation}

Si $P_{1}\subseteq P_{2}$, on a un
homomorphisme injectif canonique
$\mathfrak{a}_2^{*} \hookrightarrow \mathfrak{a}_1^{*}$ 
qui donne la projection 
$\mathfrak{a}_1  \twoheadrightarrow \mathfrak{a}_2$, 
dont on note 
$\mathfrak{a}_{1}^{2} =\mathfrak{a}_{P_{2}}^{P_{2}}$ le noyau.
On a aussi l'inclusion 
$\mathfrak{a}_{{2}}\hookrightarrow \mathfrak{a}_{{1}}$, 
qui est une section de
$\mathfrak{a}_1 \twoheadrightarrow \mathfrak{a}_2$, 
grâce à la restriction des caractères de $A_{{1}}$ à
$A_{{2}}$. 
Il s'ensuit  que si 
$P_{1} \subseteq P_{2}$ alors
\begin{equation}\label{eq:decomp}
\mathfrak{a}_1 = \mathfrak{a}_1^2\oplus \mathfrak{a}_2.
\end{equation}
Conformément à cette décomposition, on pose aussi 
$(\all_{1}^{2})^{*} = 
\{\la \in \all_{1}^{*}| \la(H) = 0 \ \forall H \in \all_{2}\}$. 
On aura besoin aussi de 
$(\all_{1,\C}^{2})^{*} := (\all_{\tlQ}^{\tlG})^{*} \otimes_{\R} \C$ 
et de $\all_{1,\C}^{*} = \all_{1}^{*} \otimes_{\R} \C$.

Si $P \sbs Q$ sont des sous-groupes paraboliques semi-standards 
où $P$ est un sous-groupe de Borel on note simplement 
$\all_{0} = \all_{P}$,
$\all_{0}^{Q} = \all_{P}^{Q}$, $\all_{0}^{*} = \all_{P}^{*}$ etc. Cela ne dépend pas du choix 
de $P$. 
En général donc, si $P_{1} \sbs P_{2}$, 
grâce à la décomposition (\ref{eq:decomp}) ci-dessus, on considère les espaces
$\all_{1}$ et $\all_{1}^{2}$
(resp. $\all_{1}^{*}$ et $(\all_{1}^{2})^{*}$) comme des sous-espaces de $\all_{0}$ (resp. $\all_{0}^{*}$).

Notons $\Delta_{P}^{G} = \Delta_{P}$ 
l'ensemble de racines 
simples pour l'action de $A_{P}$ sur $N_{P}$. 
Il y a une correspondance bijective entre les sous-groupes paraboliques $P_{2}$
contenant $P_{1}$ et les sous-ensembles 
$\Delta_{1}^{2} = \Delta_{P_{1}}^{P_{2}}$ de 
$\Delta_1 = \Delta_{P_{1}}$. En fait, 
$\Delta_1^2$
est l'ensemble de racines simples 
pour l'action de $A_1$ sur $N_1 \cap M_2$ et 
l'on a
\begin{displaymath}
\mathfrak{a}_2 = \{H \in \mathfrak{a}_1| \al(H) = 0 
\ \forall \al \in \Delta_1^2\}. 
\end{displaymath}
En plus $\Delta_{1}^{2}$ (les restrictions de ses éléments 
à $\all_{1}^{2}$) est une base de $(\all_{1}^{2})^{*}$.

Fixons $P_{1} \subseteq P_{2}$ et soit $B \in \calP(M_{0})$ 
contenu dans $P_{1}$. On a alors l'ensemble 
$\Delta_{B}^{\vee} = \{\al^{\vee} \in \all_{0}| \al  \in \Delta_{B}\}$ 
de coracines simples associées aux racines simples $\Delta_{0}$. 
Soit $(\Delta_{{1}}^{{2}})^{\vee}$ l'ensemble de projections d'éléments de 
$\Delta_{B}^{\vee}$ à $\all_{1}^{2}$ privé de zéro. Cela ne dépend pas du choix de $B$. 
L'ensemble 
$\Delta_{{1}}^{{2}}$ est en bijection canonique avec $(\Delta_{{1}}^{{2}})^{\vee}$, 
la bijection étant: à $\al \in \Delta_1^2$ on associe 
l'unique $\al^{\vee} \in (\Delta_{{1}}^{{2}})^{\vee}$ tel 
que $\al(\al^{\vee}) >0$. 
Notons également $\hDelta_1^2$ et 
$(\hDelta_1^2)^{\vee}$ 
les bases de $(\mathfrak{a}_{{1}}^{{2}})^{*}$ et 
$\mathfrak{a}_{{1}}^{{2}}$ duales \`a 
$(\Delta_{{1}}^{{2}})^{\vee}$ et
$\Delta_{{1}}^{{2}}$ 
 respectivement. 
Si $P_{2} = G$ on note simplement $\Delta_{{1}}, 
\Delta_{{1}}^{\vee}$ etc. 

Soient $P, P_{1}, P_{2} \in \calF(M_{0})$, on note
\begin{displaymath}
\mathfrak{a}_{P}^{+} = 
\{H \in \mathfrak{a}_{P}| \alpha(H) > 0 \
\forall \alpha \in \Delta_{P}\}
\end{displaymath}
et si $P_{1} \subseteq P_{2} $ 
notons $\tau_{1}^{2}$, $\htau_1^2$
les fonction caractéristiques de
\begin{equation*}
\{H \in \mathfrak{a}_0| \alpha(H) > 0 \
\forall \alpha \in \Delta_1^2\}, \quad  
\{H \in \mathfrak{a}_{0}| \varpi(H) > 0 \
\forall \varpi\ \in \hDelta_{{1}}^{{2}}\}
\end{equation*}
respectivement. 
On note $\tau_{P}$ pour $\tau_{P}^{G}$ et 
$\htau_{P}$ pour $\htau_{P}^{G}$. 

Soit $\A = \A_{\rmF}$ l'anneau des adèles de $\rmF$ 
et soit $|\cdot |_{\A}$ la valeur absolue standard sur le groupe 
des idèles $\A^{*}$.
 Pour tout $P \in \calF(M_{0})$, posons 
 $H_{P} : M_{P}(\A) \rightarrow \mathfrak{a}_{P}$
défini comme
\begin{equation*}
\langle H_{P}(m),\chi \rangle = \log (|\chi (m)|_{\A}), \quad \chi \in 
\Hom_{\rmF}(M_{P}, \Gm), \ m \in M_{P}(\A).
\end{equation*}
C'est un homomorphisme continu et surjectif, 
donc si l'on note $M_{P}(\A)^{1}$ 
son noyau, on obtient la suite exacte suivante
\begin{displaymath}
1 \rightarrow M_{P}(\A)^{1} \rightarrow M_{P}(\A) 
\rightarrow \mathfrak{a}_{P} \rightarrow 0.
\end{displaymath}
Soit $A_{P}^{\infty}$ la composante neutre du groupe des $\R$-points du 
tore déployé et défini sur $\Q$ maximal pour cette propriété dans 
le $\Q$-tore 
$\Res_{\rmF/\Q}A_{P}$. Alors, comme $\rmF \otimes_{\Q} \R$ 
s'injecte dans $\A$, on a 
naturellement $A_{P}^{\infty}\hookrightarrow A_{P}(\A) \hookrightarrow M_{P}(\A)$. En plus, 
la restriction de $H_{P}$ \`a $A_{P}^{\infty}$ est 
un isomorphisme donc 
$M_{P}(\A)$ est un produit direct de 
$M_{P}(\A)^{1}$ et $A_{P}^{\infty}$. 
Pour $Q \in \calF(M_{0})$ contenant $P$ on pose $A_{P}^{Q,\infty} = A_{P}^{\infty} \cap M_{Q}(\A)^{1}$. 
L'application $H_{P}$ induit alors un isomorphisme entre $A_{P}^{Q,\infty}$ et $\all_{P}^{Q}$.

Fixons $K$ un sous-groupe compact maximal admissible de 
$G(\A)$ par rapport à $M_{0}$. La notion d'admissibilité 
par rapport à un sous-groupe de Levi minimal est définie 
dans le paragraphe 1 de \cite{arthur2}. 
On a donc, que pour tout sous-groupe parabolique semi-standard $P$, 
$K \cap M_{P}(\A)$ est admissible dans $M_{P}(\A)$ 
et on obtient aussi la décomposition d'Iwasawa
$G(\A) = P(\A)K = N_{P}(\A)M_{P}(\A) K$ ce 
qui nous permet d'étendre $H_{P}$ à $G(\A)$ en posant 
$H_{P}(x) = H_{P}(m)$ où $x = nmk$ avec 
$m \in M_{P}(\A), n \in N_{P}(\A), k \in K$. Dans ce cas $H_{P}(x)$ 
ne dépend pas du choix de $m$.

On note $\Omega^{G}$ 
le groupe de Weyl de $(G,M_{{0}})$. 
Pour tout $s \in \Omega^{G}$, on choisit un représentant $w_{s}$
dans l'intersection de $U(\rmF) \cap K$ 
avec le normalisateur de $M_{{0}}$. 
Ceci n'est pas toujours possible donc on impose la condition sur 
$K$ que cela est possible. On peut toujours trouver un tel $K$ dans le 
cas où $G = \Gl_{n}$ et c'est le cas qui va nous intéresser. 
Pour un $\rmF$-sous-groupe $H$ de $G$ et $s \in \Omega^{G}$ 
on note $sH$ le $\rmF$-sous-groupe $w_{s}Hw_{s}^{-1}$. 
Le groupe $\Omega^{G}$ agit donc ainsi sur $\calF(M_{0})$.
Pour tout $B \in \calP(M_{0})$ l'application 
$\Omega^{G} \ni s \mapsto sB \in \calP(M_{0})$ est une bijection.
Pour tout $P \in \calF(M_{0})$ 
soit $\Omega^{P}$ le sous-groupe de 
$\Omega^{G}$ stabilisant $P$. On a donc 
$\Omega^{P} = \{s \in \Omega^{G}| w_{s} \in M_{P}(\rmF)\}$. 

Notons finalement, que parfois, pour économiser l'espace,
 on utilisera la notation $[H]$ pour noter $H(\rmF) \bsl H(\A)$.
 
 \subsection{Le domaine de Siegel}\label{par:Siegel}

Soit $B$ un sous-groupe de Borel de $G$ semi-standard. 
Pour un compact $\omega \subseteq N_{B}(\A)M_{0}(\A)^{1}$ et 
un réel négatif $c$, on définit le domaine de Siegel 
$\mathfrak{S}_{B}(\omega,c)$ dans $G(\A)$ comme:
\begin{displaymath}
\mathfrak{S}_{B}(\omega,c)= 
 \{mak \in G(\A)| m \in \omega, k \in K, a \in A_{B}^{\infty}(c)\}
\end{displaymath}
où $A_{B}^{\infty}(c) = A_{B}^{\infty}(G,c) 
= \{a \in A_{B}^{\infty}|
\al(H_{B}(a)) > c, \forall \al \in \Delta_{B}\}$.
En général, pour un sous-groupe parabolique 
semi-standard
$P$ de $G$ contenant $B$ on définit:
\begin{displaymath}
\mathfrak{S}_{B}^{P}(\omega,c) = 
\{mak \in G(\A)|  
m \in \omega, k \in K, a \in A_{B}^{\infty}(P,c)\}.
\end{displaymath}
où $A_{B}^{\infty}(P,c) = \{a \in A_{B}^{\infty}|
\al(H_{B}(a)) > c, \forall \al \in \Delta_{B}^{P}\}$.

On utilisera le résultat suivant 
de la théorie de réduction, qu'on peut trouver, par exemple 
dans \cite{godement}:

\brop\label{prop:siegelDomain}
Il existe un réel négatif $c_{0}$ 
et pour tout sous-groupe de Borel semi-standard $B$ de $G$
un compact $\omega_{B} 
 \subseteq N_{B}(\A)M_{0}(\A)^{1}$ tels que pour 
tout sous-groupe parabolique semi-standard
$P$ contenant  $B$ l'on a:
\begin{displaymath}
G(\A) = P(\rmF) \mathfrak{S}_{B}^{P}(\omega_{B},c_{0}).
\end{displaymath}
\erop

Fixons la constante $c_{0}$ comme ci-dessus. 
Pour tout sous-groupe de Borel semi-standard $B$ 
on fixe aussi un $\omega_{B}$ comme dans la proposition 
\ref{prop:siegelDomain} de façon que si 
$B' \in \calP(M_{0})$ est tel que $sB = B'$ 
on a $\omega_{B'} = w_{s}\omega_{B}w_{s}^{-1}$.
Les définitions de ce paragraphe sont valables en particulier 
pour les sous-groupes de Levi de $G$. On voit donc qu'on peut 
fixer un $\omega_{B}$ de façon que 
pour tout sous-groupe parabolique semi-standard $P$ et tout $B \in \calP(M_{0})$ 
le contenant
le compact $M_{P}(\A) \cap \omega_{B}$ 
ainsi que le sous-groupe $K_{P} := M_{P}(\A) \cap K$ 
jouissent des rôles de $\omega_{B}$ et $K$ ci-dessus 
par rapport au groupe réductif $M_{P}$ et son 
sous-groupe de Borel $B \cap M_{P}$, la constante $c_{0}$ 
restant la même.  

Soient $B \in \calP(M_{0})$, $P \sps B$ et $T \in \all_{0}$. On définit 
$F_{B}^{P}(x,T)$ comme la fonction 
caractéristique de l'ensemble:
\begin{displaymath}
\{x \in G(\A)| \ \exists \delta \in P(\rmF) \
\delta x \in \mathfrak{S}_{B}^{P}(\omega_{B},c_{0})
, \ \varpi(H_{B}(\delta  x)-T) < 0 \
\forall \varpi \in \hDelta_{B}^{P}\}.
\end{displaymath}
Autrement dit, si l'on pose
\begin{equation*}
A_{B}^{\infty}(P,c_{0},T) := 
\{a \in A_{B}^{\infty}(P,c_{0})| 
\varpi(H_{B}(a)-T) < 0, \ \forall \varpi \in \hDelta_{B}^{P}\},
\end{equation*}
on a alors que $F_{B}^{P}(\cdot, T)$ 
est la fonction caractéristique de 
la projection de $\omega_{B} A_{B}^{\infty}(P,c_{0},T)K$ 
sur 
$A_{P}^{\infty}N_{P}(\A)M_{P}(\rmF)\backslash G(\A)$.

\subsection{Algèbres de Lie}\label{par:algLie}

Soit $\gl = \Lie(G)$, pour tout $P, P_{1}, P_{2} \in \calF(M_{0})$ 
tels que $P_{1} \sbs P_{2}$
on note $N_{1}^{2} = N_{1} \cap M_{2}$, 
$\ml_{P} = \Lie(M_{P})$, $\nl_{P} = \Lie(N_{P})$ et 
$\nl_{1}^{2} = \Lie(N_{1}^{2})$. Soit $\brP \in \calF(M_{P})$ 
le sous-groupe parabolique opposé à $P$ 
(i.e. tel que $\brP \cap P = M_{P}$). 
On note alors $\brnl_{P} = \nl_{\brP}$ 
et $\brnl_{1}^{2} = \ml_{2} \cap \brnl_{1}$. 

On fixe une forme bilinéaire  non-dégénérée
$\langle \cdot \ , \cdot\rangle$ 
sur $\gl$ invariante par adjonction. 
Pour tout $P \in \calF(M_{0})$
la restriction de la forme $\bilif$ à 
$\brnl_{P} \times \nl_{P}$ 
est alors non-dégénérée, donc 
l'espace $\brnl_{P}$
s'identifie à l'espace dual à $\nl_{P}$ 
grâce à cette forme.

\subsection{Fonctions de Bruhat-Schwartz}\label{bschwartzG}

On note  
$\A_{f}$ l'anneau des adèles finis de $\rmF$ et 
$\rmF_{\infty} := \rmF \otimes_{\Q} \R$, de façon qu'on a 
$\A = \rmF_{\infty}\times \A_{f}$. 
Notons $\mathcal{S}(\gl(\A))$ l'ensemble des 
fonctions des Bruhat-Schwartz sur $\gl(\A)$ 
i.e. l'espace de fonctions sur 
$\gl(\A)$ engendré par des fonctions du type 
$f_{\infty}\otimes \chi^{\infty}$ où 
$f_{\infty}$ est une fonction de la classe de 
Schwartz sur $\gl(\rmF_{\infty})$ et 
$\chi^{\infty}$ est une fonction caractéristique d'un 
compact ouvert de $\gl(\A_{f})$. 

\subsection{Les mesures de Haar}\label{haarSect}

Soit $P$ un sous-groupe parabolique semi-standard de $G$.
On fixe $dx$ une mesure de Haar sur $G(\A)$, 
ainsi que 
pour tout sous-groupe connexe $V$ de $N_{P}$ 
(resp. toute sous-alg\`ebre $\mathfrak{h}$ de $\nl_{P}$)
l'unique mesure de Haar sur $V(\A)$ (resp. $\mathfrak{h}(\A)$)
pour laquelle le volume de $V(\rmF)\backslash V(\A)$ 
(resp. $\mathfrak{h}(\rmF)\backslash \mathfrak{h}(\A)$)
soit $1$. Choisissons la mesure de Haar $dk$ sur $\rmK$ 
normalisé de m\^eme façon. 

On fixe aussi une norme 
euclidienne $\|\cdot \|$ sur $\mathfrak{a}_{0}$ 
invariante par le groupe de Weyl
$\Omega^{G}$
et sur tout sous-espace de $\all_{0}$ la mesure de Haar compatible 
avec cette norme. 
Pour tout $Q \in \calF(M_{0})$ tel que $Q \sps P$, 
on en déduit
les mesures de Haar sur 
$A_{P}^{Q,\infty}$ et $A_{P}^{\infty}$
via l'isomorphisme 
$H_{P}$. 

Soit $dp$ la mesure de Haar sur $P(\A)$ invariante 
\`a gauche normalisé de façon que $dx = dpdk$ 
(grâce à la décomposition d'Iwasawa).
Notons $\rho_{P}^{G} = \rho_{P}$ 
l'élément de $(\all_{P}^{G})^{*}$ 
tel que $d(\Ad(m^{-1})n) = e^{2\rho_{P}(H_{P}(m))}dn$ 
pour $m \in M_{P}(\A)$ et $n \in N_{P}(\A)$.
Il s'ensuit 
qu'il existe une unique 
mesure de Haar $dm$ sur 
$M_{P}(\A)$ telle 
que si l'on écrit 
$p = nm$ 
où $p \in P(\A)$, $n \in N_{P}(A)$ et
$m \in M_{P}(\A)$ 
alors $dp = e^{-2\rho_{P}(H_{P}(m))}dndm$. Les mesures de Haar sur 
$M_{P}(\A)$ et $A_{P}^{\infty}$ induisent alors une unique mesure de Haar sur 
$M_{P}(\A)^{1}$, que l'on fixe, telle que la mesure de Haar sur $M_{P}(\A)$ 
soit le produit de mesures sur $A_{P}^{\infty}$ et sur $M_{P}(\A)^{1}$.

\subsection{\texorpdfstring{$\Gl_{n} \hrar \Gl_{n+1}$}{Gl(n) -> Gl(n+1)}}\label{par:glnGlnplus1}

Soit $W$ un $\rmF$-espace vectoriel de dimension finie $n+1$ 
et soit $V \sbs W$ un sous-espace de dimension $n$, où $n \in \N$.
Notons $\tlG = \Gl(W)$.
Fixons un vecteur 
$e_{0} \in W \smallsetminus V$ et notons $D_{0}$ la droite 
qu'il engendre. 
On a alors $W = V \oplus D_{0}$ ce qui permet 
d'identifier $G = \Gl(V)$ comme un sous-groupe de 
$\tlG$ stabilisant $V$ et fixant $e_{0}$. 
Choisissons $n$-droites $D_{1}, \ldots, D_{n}$ dans $V$ qui engendrent $V$.
Soit $M_{0}$ le stabilisateur 
dans $G$ des droites 
$(D_{i})_{i=1, \ldots, n}$. 
C'est un sous-$\rmF$-groupe de Levi minimal de $G$. 
Soit alors $M_{\tlzero}$ l'unique sous-$\rmF$-groupe de Levi minimal de $\tlG$ 
contenant $M_{0}$. Alors $M_{\tlzero}$ est le stabilisateur des droites
$(D_{i})_{i=0,\ldots, n}$ dans $\tlG$.

Les résultats des paragraphes précédentes s'appliquent 
aux groupes $G$ et $\tlG$ et leurs sous-groupes de Levi minimales $M_{0}$
et $M_{\tlzero}$. Les objets associés à $\tlG$ seront notés toujours avec un tilde. 
Pour le choix du sous-groupe compact maximal, on fixe des vecteurs non-nuls 
$e_{i} \in D_{i}$ pour tout $i =1, \ldots, n$ ce qui avec le choix du vecteur $e_{0}$ 
défini les isomorphismes $\tlG \cong \Gl_{n+1}$ et $G \cong \Gl_{n}$. 
On pose alors $\tlK = \prod_{v}\tlK_{v}$ où, pour une place fini $v$ de $\rmF$ on note 
$\tlK_{v} = \Gl_{n+1}(\calO_{v})$, où $\calO_{v}$ c'est l'anneau des entiers de la complétion de $\rmF$ en $v$, 
pour une place réelle $v$ on pose $\tlK_{v} = O(n+1)$-le groupe orthogonal anisotrope et 
pour une place complexe on met $\tlK_{v} = U(n+1)$- le groupe unitaire anisotrope. 
On pose aussi $K = \tlK \cap G(\A)$. 
Dans ce cas $\tlK$ et $K$ vérifient les conditions 
du paragraphe (\ref{par:prelimstraceG}). 
Les inclusions $G \hrar \tlG$  et $M_{0} \hrar M_{\tlzero}$ induisent l'inclusion 
$\Omega^{G} \hrar \Omega^{\tlG}$.
On choisit aussi des représentants du groupe de Weyl 
$\Omega^{\tlG}$ de $\tlG$ comme les éléments permutants les vecteurs $e_{i}$. 
On a alors pour tout $\tls \in \Omega^{\tlG}$ que $w_{\tls} \in \tlG(\rmF) \cap \tlK$ 
et si $s \in \Omega^{G}$ alors $w_{s} \in G(\rmF) \cap K$.

On identifie
$\all_{0}$ et $\all_{0}^{*}$ avec des sous-espaces de $\all_{\tlzero}$ 
et $\all_{\tlzero}^{*}$ respectivement. En particulier la mesure de Haar et la norme 
euclidienne sur $\all_{0}$ sont celles d'un sous-espace de $\all_{\tlzero}$.

Posons $V_{i} = \oplus_{j = 1}^{i} D_{i} \sbs V$. On fixe 
$B \in \calP(M_{0})$ le stabilisateur du drapeau 
\begin{equation}\label{eq:drapSt}
0 \sbn V_{1} \sbn \cdots \sbn V_{n} = V.
\end{equation}
On appelle $P \in \calF(M_{0})$ tels que $P \sps B$ les sous-groupes paraboliques 
standards de $G$. Tout $P \in \calF(M_{0})$ standard est alors défini comme le stabilisateur 
d'un sous-drapeau du drapeau (\ref{eq:drapSt}) ci-dessus.

Notons $V^{*} = \Hom_{\rmF}(V, \rmF)$ et $W^{*} = \Hom_{\rmF}(W, \rmF)$.
Pour un $\rmF$-sous-espace $\calV \sbs V$ de 
type $\calV = \bigoplus_{i = k}^{l}D_{i}$ où $1 \le k \le l \le n$ 
on pose $\calV^{*} = \{\la \in V^{*}| \, \la|_{D_{i}} = 0, \, 1 \le i < k, \ l < i \le n \}$ 
et $\calV^{\bot} = \{ \la \in V^{*}| \la|_{\calV} = 0\}$.

Pour tout $\tlP \in \calF(M_{\tlzero})$ on admet la notation:
\begin{equation}
P := \tlP \cap G. 
\end{equation}
Alors $P$ est un sous-groupe parabolique de $G$ semi-standard.
Le groupe $G$ (resp. $\tlG$) agit aussi naturellement sur $V^{*}$ (resp. $W^{*}$) donc 
aussi sur $V \times V^{*}$ (resp. $W \times W^{*}$). 
Pour $\tlP \in \calF(M_{\tlzero})$
on note
$\calV_{\tlP} \sbs V \times V^{*}$ 
le plus grand sous-espace de $V \times V^{*}$ stabilisé par $\tlP$ vu comme un sous-espace 
de $W \times W^{*}$.
On note aussi $Z_{\tlP} \sbs V$ le plus petit sous-espace de $V$ 
tel que $M_{\tlP}$ stabilise $Z_{P} \oplus D_{0} \sbs W$. 
On pose $\calZ_{\tlP} = Z_{\tlP} \times Z_{\tlP}^{*} \sbs V \times V^{*}$.

On note $\relPb = \{P \in \calF(M_{\tlzero}) | \, B \sbs \tlP\}$. On 
appelle les éléments de $\relPb$ les sous-groupes 
paraboliques \textit{relativement standards} de 
$\tlG$. 

Soit $\tlP \in \relPb$. Il est défini alors par 
une suite des entiers $0=i_{0} \le \cdots \le i_{l}=n$ 
et un entier $k$ tels que 
$1 \le k \le l \le n+1$ et 
$k$ vérifie la propriété que si
$i_{j-1} = i_{j}$ pour un $1 \le j \le l$ 
alors $j=k$. Le sous-groupe $\tlP$ est alors le stabilisateur 
du drapeau:
\begin{equation}\label{eq:drapeauRel}
0 = V_{i_{0}} \subsetneq \cdots \subsetneq  V_{i_{k-1}}
\subsetneq 
V_{i_{k}} \oplus D_{0} \subsetneq \cdots \subsetneq V_{i_{l}} 
\oplus D_{0} =W.
\end{equation}
Dans ce cas $P$ est le stabilisateur de 
\begin{equation*}
0 = V_{i_{0}} \subsetneq \cdots \subsetneq  V_{i_{k-1}}
\subseteq 
V_{i_{k}} \subsetneq \cdots \subsetneq V_{i_{l}} =V.
\end{equation*}
On a en plus
\begin{equation}
\calV_{\tlP} = V_{i_{k-1}} \times V_{i_{k}}^{\bot} \sbs V \times V^{*}, \quad 
Z_{\tlP} = \bigoplus_{i_{k-1} < j \le i_{k}}D_{i}.
\end{equation}
On voit qu'on a des isomorphismes:
\begin{gather}
M_{\tlP} =
\prod_{1 \le j < k} \Gl_{i_{j}-i_{j-1}}
\times \Gl(Z_{\tlP} \oplus D_{0})
\times \prod_{k < j \le l} \Gl_{i_{j}-i_{j-1}} 
\cong  \tlG_{\tlP} \times H_{\tlP} \label{eq:MPtildedecomp}
\\
M_{P} =
\prod_{1 \le j < k} \Gl_{i_{j}-i_{j-1}}
\times \Gl(Z_{\tlP})
\times \prod_{k < j \le l} \Gl_{i_{j}-i_{j-1}} 
\cong  G_{\tlP} \times H_{\tlP}\label{eq:MPdecomp}
\end{gather}
où $\tlG_{\tlP} = \Gl(Z_{\tlP} \oplus D_{0})$, 
$G_{\tlP} = \Gl(Z_{\tlP})$ et 
l'on voit $H_{\tlP}$ comme un 
sous-groupe de $M_{P}$, donc aussi de $M_{\tlP}$, 
qui agit trivialement sur 
$Z_{\tlP}$. Dans ce contexte on note
$A_{\tlP}^{st,\infty} := A_{H_{\tlP}}^{\infty}$ - c'est la partie standard du tore $A_{\tlP}^{\infty}$ 
puisque $A_{\tlP}^{\infty} \cap G(\A) = A_{\tlP}^{st,\infty}$, d'où la notation.
On pose aussi $\all_{\tlP}^{st} := \all_{H_{\tlP}}$, donc 
$\all_{\tlP}^{st} \sbs \all_{\tlzero}$ 
ce qui détermine les mesures de Haar sur $\all_{\tlP}^{st}$ et sur $A_{\tlP}^{st, \infty}$. 
On a alors:
\[
M_{P}(\A) = A_{\tlP}^{st,\infty} (H_{\tlP}(\A)^{1} \times G_{\tlP}(\A)), \quad 
M_{\tlP}(\A) = A_{\tlP}^{st,\infty} (H_{\tlP}(\A)^{1} \times \tlG_{\tlP}(\A)).
\]
On fixe donc l'unique mesure de Haar sur $H_{\tlP}(\A)^{1} \times G_{\tlP}(\A)$ 
de façon que la mesure de Haar sur $M_{P}(\A)$ choisie soit produit de cette mesure et celle sur $A_{\tlP}^{st,\infty}$.

On fixe un $c_{0} < 0$  tel que la proposition \ref{prop:siegelDomain} ci-dessus soit vraie
pour tous les groupes $\Gl_{k}$ définis
sur $\rmF$ o\`u $k \le n+1$. 
On suppose alors que pour tout sous-groupe de Borel relativement standard $\tlB$ de $\tlG$ 
on a $\omega_{B} \sbs \omega_{\tlB}$. 

On remarque finalement une conséquence de notre choix des représentants 
du groupe de Weyl. On suppose donc que pour tout 
$\tls \in \Omega^{\tlG}$ on a $w_{\tls} \in \tlG(\rmF) \cap \tlK$. 
Soient $\tls \in \Omega^{\tlG}$, 
$\tlB \in \calP(M_{\tlzero})$ et $\tlP, \tlQ \in \calF(M_{\tlzero})$ tels que $\tlP \sbs \tlQ$.
Alors pour tout $x \in \tlG(\A)$ et tout $T \in \all_{\tlzero}$ 
les formules suivantes sont vraies:
\begin{gather}
\htau _{\tls \tlP}^{\tls \tlQ}(H_{\tls \tlP}(x)-T) =
\htau_{\tlP}^{\tlQ}(H_{\tlP}(w_{\tls}^{-1}x)-\tls^{-1}T),\label{eq:weyl2} \\
F^{\tls P}_{\tls B}(x,T) = F^{\tlP}_{\tlB}(w_{\tls}^{-1}x,\tls^{-1}T)\label{eq:weyl3}.
\end{gather}

\subsection{Relation entre \texorpdfstring{$\all_{\tlP}^{st}$ et $\all_{\tlP}$}{certains espaces}}\label{par:HQtlQ}

L'application naturelle $\all_{\tlP}^{st} \rar \all_{\tlP}^{\tlG}$, induit par la projection $\all_{\tlP} \rar \all_{\tlP}^{\tlG}$, 
est un isomorphisme. 
On note $j_{\tlP}$ son Jacobien
et on note $\iota_{\tlP}^{st} : (\all_{\tlP,\C}^{\tlG})^{*} \rar (\all_{\tlP,\C}^{st})^{*} := \Hom_{\R}(\all_{\tlP}^{st}, \C) = \all_{H_{\tlP},\C}^{*}$ 
l'isomorphisme induit.
Si $\tlP = \tlG$ on met $j_{\tlP} = 1$.

Ainsi, pour toute fonction $\uphi$ sur $\all_{\tlP}$ qui est $\all_{\tlG}$-invariante, 
et tout $\la \in (\all_{\tlP,\C}^{st})^{*}$
on a:
\begin{equation}\label{eq:HQtlQ}
\int_{\all_{\tlP}^{st}} e^{\la(H)} \uphi(H)dH = j_{\tlP}^{-1}
\int_{\all_{\tlP}^{\tlG}}
e^{\iota_{\tlP}^{st}(\la)(H)}
\uphi(H)dH.
\end{equation}



\subsection{Quelques résultats de la théorie de la réduction}\label{par:reductionTheory}

On fixe un sous-groupe de Borel semi-standard $\tlBmin \in \calP(M_{\tlzero})$ de $\tlG$. 
Pour tout $H \in \all_{\tlzero}$ et tout $\tlP \in \calF(M_{\tlzero})$ on note 
$H_{\tlP} \in \all_{\tlP}$ la projection de $\tls H$ à $\all_{\tlP}$ 
où $\tls \in \Omega^{\tlG}$ est tel que $\tls \tlBmin \sbs \tlP$. Cette définition ne dépend pas 
du choix de $\tls$. 

La relation (\ref{eq:weyl3}) nous permet de poser la définition 
suivante: soit $\tlP \in \calF(M_{\tlzero})$ 
et $\tls \in \Omega^{\tlG}$ tel que $\tls \tlBmin \sbs \tlP$, on pose donc:
\begin{equation}
F^{\tlP}(x,T) := F_{\tls \tlBmin}^{P}(x,T_{\tls \tlB_{0}}), \quad 
x\in \tlG(\A), \ T \in \all_{\tlzero}.
\end{equation}
Cette définition ne dépend pas du choix de 
$\tls \in \Omega^{\tlG}$ car $F^{\tlP}_{\tlB}$ est 
$M_{\tlP}(\rmF)$-invariant pour tout sous-groupe de 
Borel semi-standard $\tlB$ contenu dans $\tlP$. 
Remarquons 
aussi que si le paramètre $T$ appartient 
à la chambre $\all_{\tlBmin}^{+}$ alors $T_{\tls \tlBmin}$ est dans la chambre 
$\all_{\tls \tlBmin}^{+}$. On note désormais:
\[
\all_{\tlzero}^{+} := \all_{\tlBmin}^{+}.
\]

Le lemme suivant est 
l'application du lemme 6.4 de \cite{arthur3} à notre 
situation:

\blem\label{lem:lemme64}
Il existe un $T_{+} \in \all_{\tlzero}^{+}$ 
tel que pour tout sous-groupe parabolique 
semi-standard $\tlQ$ de $\tlG$, 
tout sous-groupe de Borel semi-standard 
$\tlB$ contenu dans $\tlQ$, 
tout $T \in T_{+}+\all_{\tlzero}^{+}$ 
et tout $x \in \tlG(\A)$ l'on a:
\begin{displaymath}
\sum_{\tlB \subseteq  \tlP \subseteq \tlQ}
\sum_{\delta \in \tlP(\rmF)\backslash \tlQ(\rmF)}
F^{\tlP}(\delta x,T)
\tau_{\tlP}^{\tlQ}(H_{\tlP}(\delta x)-T_{\tlP}) =1.
\end{displaymath}
\elem
\noindent
Fixons alors un $T_{+} \in \all_{\tlzero}^{+}$ 
comme dans le lemme \ref{lem:lemme64} ci-dessus.

On a alors l'analogue relatif du lemme \ref{lem:lemme64}.

\brop\label{harderNarashiman}
 Soit $\tlQ$ un sous-groupe parabolique 
 relativement standard
 de $\tlG$.  
Alors, pour tout $T \in T_{+} + \all_{\tlzero}^{+}$ 
et tout $x \in G(\A)$, l'on a:
\begin{displaymath}
\sum_{B \sbs \tlP \subseteq \tlQ}
\sum_{\eta \in P(\rmF)\backslash Q(\rmF)}
F^{\tlP}(\eta x,T)\tau_{\tlP}^{\tlQ}
(H_{\tlP}(\eta x)-T_{\tlP})=1.
\end{displaymath}
où la somme porte sur tous les sous-groupes paraboliques 
relativement standards de 
$\tlG$ contenus dans $\tlQ$ et 
on note $P = \tlP \cap G$ etc. 
\bdem 
Soient $\tlQ$ et $T$ comme dans l'énoncé de 
la proposition. Fixons un sous-groupe 
de Borel relativement standard $\tlB$ contenu dans $\tlQ$.
 Pour un $\tlP \in  \relPb$ tel que 
$\tlB \subseteq \tlP \subseteq \tlQ$. On fixe
$\Omega^{\tlG}_{\tlP, \tlQ}$  un ensemble
de représentants dans $\{\tls \in \Omega^{\tlG}| \tls^{-1}\tlP \sbs  \tlQ \text{ et }\tls^{-1} \tlP \in  \relPb\}$ 
pour la relation $\tls_{1} \sim \tls_{2}$ 
si est seulement si $\tls_{1}\tls_{2}^{-1}\tlP = \tlP$. 
Pour $\tls \in \Omega^{\tlG}_{\tlP, \tlQ}$ on note 
$P_{\tls} = (\tls^{-1}\tlP) \cap G$.

On voit donc que la somme dans le lemme égale
\[
\sum_{\tlB \subseteq  \tlP \subseteq \tlQ} 
\sum_{\tls \in \Omega^{\tlG}_{\tlP, \tlQ}}
\sum_{\eta  \in P_{\tls}(\rmF)\backslash Q(\rmF)}
F^{\tlP}_{\tlB}(w_{\tls}\eta x, T_{\tlB})
\tau_{\tlP}^{\tlQ}
(H_{\tlP}(w_{\tls}\eta x)-T_{\tlP}).
\]
Sous cette forme le résultat est démontré dans \cite{ichYamGl}, lemme 2.3. 
\edem
\erop

Soit $\tlP$ un sous-groupe parabolique 
relativement standard de $\tlG$. On note 
\begin{equation}\label{eq:relBrelToP}
\calP(M_{\tlzero}, B, \tlP) = \{\tlB \in \calP(M_{\tlzero}) \cap \relPb| \tlB \sbs \tlP\}.
\end{equation}
Posons 
$\omega_{P} = \omega_{B} \cap M_{P}(\A)$ 
et $K_{P} = K \cap M_{P}(\A)$. 
On a alors

\blem\label{lem:thatsAllThereIs} Soit $\tlP \in  \relPb$. 
Le sous-ensemble suivant de $M_{P}(\A)$
\begin{equation*}
\bigcup_{ \tlB \in \calP(M_{\tlzero}, B, \tlP)}
\omega_{P}
(A_{B}^{\infty}(P,c_{0}) \cap A_{\tlB}^{\infty}(\tlP,c_{0}))K_{P}
\end{equation*}
se surjecte sur $M_{P}(\rmF)\backslash M_{P}(\A)$. 
\bdem
En vertu de la proposition \ref{prop:siegelDomain} il suffit de montrer que 
pour tout 
$a \in A_{B}^{\infty}(P, c_{0})$ il existe un $\tlB \in \calP(M_{\tlzero}, B, \tlP)$ tel que 
$a \in A_{\tlB}^{\infty}(\tlP, c_{0})$. On se ramène facilement au cas $\tlP= \tlG$ et dans ce 
cas c'est l'assertion (2.5) de \cite{ichYamGl}.
\edem 
\elem

On note aussitôt:
\bcor\label{cor:thatsAllThereIs}
Soit $\tlP \in  \relPb$. L'ensemble
\begin{equation*}
\bigcup_{ \tlB \in \calP(M_{\tlzero}, B, \tlP)}
\omega_{P}
(A_{B}^{\infty}(P,c_{0}) \cap A_{\tlB}^{\infty}(\tlP,c_{0},T_{\tlB}))
K_{P},
\end{equation*}
où $\calP(M_{\tlzero}, B, \tlP)$ est défini par (\ref{eq:relBrelToP}) ci-dessus,
se surjecte sur l'ensemble des 
$m \in M_{P}(\rmF)\backslash M_{P}(\A)$ 
tels que $F^{\tlP}(m,T)=1$. 
\bdem 
Soit
$m \in M_{P}(\rmF)\backslash M_{P}(\A)$
et $\tlB \in \calP(M_{\tlzero}, B, \tlP)$ tel que 
 $k_{1}ak_{2} \in \omega_{P} (A_{B}^{\infty}(P,c_{0}) \cap A_{\tlB}^{\infty}(\tlP,c_{0}))K_{P}$ est
un représentant de $m$
dans $M_{P}(\A)$, 
comme dans le lemme \ref{lem:thatsAllThereIs} ci-dessus. 
Il résulte du lemme \ref{lem:lemme64} et du fait qu'il existe une constante $c$ telle que
$\varpi(H_{\tlB}(w_{\tls} n)) \le c$ pour tout $\varpi \in \hDelta_{\tlB}$ 
et tout $\tls \in \Omega^{\tlG}$ que si 
$\tlk_{1} \in \omega_{\tlB} \cap M_{\tlP}(\A)$, $\tla \in A_{\tlB}^{\infty}(\tlP, c_{0})$ 
et $\tlk \in \tlK \cap M_{\tlP}(\A)$ sont tels 
que $F^{\tlP}_{\tlB}(\tlk_{1} \tla \tlk, T) = 1$ alors $F^{\tlP}_{\tlB}(\tla, T) = 1$. 
En particulier, puisque $\omega_{B} \sbs \omega_{\tlB}$ pour tout $\tlB \in \relPb \cap \calP(M_{\tlzero})$ 
et $K \sbs \tlK$, 
si l'on suppose que 
$F^{\tlP}(m,T)= F^{\tlP}_{\tlB}(k_{1}ak_{2}, T_{\tlB}) = 1$ 
on a $F^{\tlP}_{\tlB}(a,T_{\tlB})=1$ d'où $a \in A_{B}^{\infty}(P,c_{0}) \cap A_{\tlB}^{\infty}(\tlP,c_{0})$ 
et le corollaire suit. 
\edem
\ecor

\subsection{Majorations cruciales}\label{majorCrucSec}

Dans ce paragraphe on se propose de donner une 
majoration de  l'intégrale du type:
\begin{equation}\label{eq:majCruc}
\int\limits_{M_{P}(\rmF)\backslash M_{P}(\A)}
e^{\upla(H_{P}(m))} 
F^{\tlP}(m,T)|\Upphi(m)| dm
\end{equation}
où $\tlP \in \relPb$, $\upla \in \all_{P}^{*}$, 
et 
$\Upphi$ est une fonction sur 
$M_{\tlP}(\A)$ invariante à gauche par 
$M_{\tlP}(\rmF)A_{\tlG}^{\infty}$.
Pour $\tlB \in \calP(M_{\tlzero}, B, \tlP)$ notons:
\[
A_{\tlB}^{\tlG, \infty}(\tlP, c_{0}, T_{\tlB}) := A_{\tlB}^{ \infty}(\tlP, c_{0}, T_{\tlB})  \cap A_{\tlB}^{\tlG,\infty}.
\]
On a alors:

\blem\label{lem:majCruc}
 Il existe un $c >0$ tel que 
 pour toute fonction $\Upphi$ comme ci-dessus, tout $T \in T_{+} + \all_{\tlzero}^{+}$, 
 tout $\upla \in \all_{P}^{*}$ 
et tout $\tlB \in \calP(M_{\tlzero}, B, \tlP)$ il existe
des $\upla_{\tlB, \tlP} \in (\all_{\tlB}^{\tlG})^{*}$ tels que l'intégrale (\ref{eq:majCruc})
est majorée par $c$ fois
\[
\sum_{
\tlB \in \calP(M_{\tlzero}, B, \tlP)
}
\int\limits_{
A_{\tlB}^{\tlG,\infty}(\tlP,c_{0},T_{\tlB})}
e^{\upla_{\tlB, \tlP}(H_{\tlB}(a))}
\sup_{\begin{subarray}{c}
k_{1} \in \omega_{\tlB} \cap M_{\tlP}(\A) \\
k_{2} \in \tlK \cap M_{\tlP}(\A)
\end{subarray}}
|\Upphi(k_{1}ak_{2})|da.
\]
\bdem
Soient $B_{P} = B\cap M_{P}$ et $B_{P}(\A)^{1} = N_{B}^{P}(\A)M_{B}(\A)^{1}$.
Alors $\omega_{P} \subseteq B_{P}(\A)^{1}$.
Le compact $K_{P}$ est un sous-groupe compact maximal 
de $M_{P}(\A)$. La 
décomposition d'Iwasawa 
implique qu'on a $M_{P}(\A) = B_{P}(\A)^{1}A_{B}^{\infty}K_{P}$.
 Il existe alors une unique mesure $db'$ sur $B_{P}(\A)^{1}$ 
invariante à gauche telle que $dm = e^{-2 \rho_{B}^{P}(H_{P}(a))} db' da dk$. 

En utilisant le corollaire \ref{cor:thatsAllThereIs}
on voit alors que l'intégrale 
(\ref{eq:majCruc}) est majorée par:
\begin{equation}\label{eq:majCruc2}
\sum_{
\tlB \in \calP(M_{\tlzero}, B, \tlP)
}
\vol(\omega_{P}) \qquad 
\int\limits_{
\mathclap{A_{B}^{\infty}(P,c_{0}) \cap A_{\tlB}^{\infty}
(\tlP,c_{0},T_{\tlB})}} \ 
e^{(\upla-2 \rho_{B}^{P})(H_{B}(a))}
\sup_{\begin{subarray}{c}
k_{1} \in \omega_{P} \\
k_{2} \in K_{P}
\end{subarray}}
|\Upphi(k_{1}ak_{2})| da.
\end{equation}
Soit $\tlB \in  \calP(M_{\tlzero}, B, \tlP)$.
On peut appliquer les résultats du paragraphe \ref{par:HQtlQ} à $\tlP = \tlB$ 
car on a $\all_{B} = \all_{\tlB}^{st}$. 
Si l'on note 
$\upla_{\tlP, \tlB} = \iota_{\tlB}^{st}(\upla-2 \rho_{B}^{P}) \in (\all_{\tlB}^{\tlG})^{*}$, 
on obtient, en vertu de l'égalité (\ref{eq:HQtlQ}), que 
(\ref{eq:majCruc2}) égale
\begin{equation*}
\sum_{
\tlB \in \calP(M_{\tlzero}, B, \tlP)
}
j_{\tlB}^{-1}
\vol(\omega_{P}) \qquad 
\int\limits_{
A_{\tlB}^{\tlG,\infty}(\tlP,c_{0},T_{\tlB})}
e^{\upla_{\tlB, \tlP}(H_{\tlB}(a))}
\sup_{\begin{subarray}{c}
k_{1} \in \omega_{P} \\
k_{2} \in K_{P}
\end{subarray}}
|\Upphi(k_{1}ak_{2})| da.
\end{equation*}
On conclut en remarquant que $\omega_{P} \sbs \omega_{\tlB} \cap M_{\tlP}(\A)$ 
et $K_{P} \sbs \tlK \cap M_{\tlP}(\A)$.
\edem
\elem


%% file: general_gln.tex
\section{Formule des traces relative pour \texorpdfstring{$\gll(n+1) // \Gl(n)$}{gl(n+1) // Gl(n)}}\label{sec:RTFtlgl}

On note $\gl = \Lie(G)$ et $\tlgl = \Lie(\tlG)$. 
Conformément au paragraphe \ref{par:algLie}, pour tout 
$\tlP \in \calF(M_{\tlzero})$, on note $\ml_{\tlP} = \Lie(M_{\tlP})$, 
$\nl_{\tlP} = \Lie(N_{\tlP})$, etc. Pour la forme bilinéaire 
$\bilif$ on choisit la forme trace. 

\subsection{Les invariants}\label{par:lesInvs}

Soit $X \in \tlgl$, 
suivant la décomposition 
$W =  V \oplus D_{0}$ 
on écrit:
\begin{equation}\label{eq:xisamatrix}
X = 
\begin{pmatrix}
B & u \\
v & d \\
\end{pmatrix} 
\end{equation}
où $B \in \gl$, 
$u \in \Hom_{\rmF}(D_{0} , V)$, 
$v \in \Hom_{\rmF}(V,  D_{0})$ et $d \in \Hom_{\rmF}(D_{0}, D_{0})$.
On rappelle qu'on a fixé un vecteur non-nul $e_{0} \in D_{0}$.
On note alors $e_{0}^{*} \in W^{*}$ défini par $e_{0}^{*}(e_{0})=1$  
et $e_{0}^{*}|_{V} \equiv 0$ et 
on identifie donc $d$ avec $e^{*}_{0}(d(e_{0})) \in \Ga$, 
$u$ avec $u(e_{0}) \in V$ et  
$v$ avec l'élément de $V^{*} = \Hom_{\rmF}(V,\rmF)$ 
défini 
par $x \mapsto e_{0}^{*}(v(x))$.

On dit que $X \in \tlgl$ 
est \textit{semi-simple régulier} s'il vérifie les conditions 
de la proposition suivante, due à \cite{rallSchiff}, théorème 6.1 et 
proposition 6.3. 

\brop\label{prop:releltG}
Soit $X = \matx{B}{u}{v}{d} \in \tlgl$,
alors les conditions 
suivantes sont équivalentes:
\begin{enumerate}[1)]
\item $\det(a_{ij}) \neq 0$ o\`u $a_{ij} = vB^{i+j}u$, 
$0 \le i,j \le n-1$.
\item Le stabilisateur de $X$ dans 
$G$ 
est trivial et 
l'orbite de $X$ dans $\tlgl$ 
pour l'action de $G$ 
est fermée pour 
la topologie de Zariski.
\end{enumerate}
\erop

On introduit alors les invariants suivants pour l'action de 
$G$ sur $\tlgl$. 
Soit $X = \matx{B}{u}{v}{d} \in \tlgl$ 
comme ci-dessus. On pose $A_{0}(X) = d$ et 
$A_{i}(X) = vB^{i-1}u$ 
pour $i = 1,2,\ldots ,n$ ainsi 
que $B_{j}(X) = \Tr \bigwedge^{j} B$ 
pour $j = 1,\ldots, n$. 
Alors, le lemme 3.1 de \cite{zhang2} dit 
que les fonctions $A_{i}$, $B_{j}$ engendrent 
l'anneau des polynômes sur $\tlgl$ invariants sous 
l'action de $G$. Ils définissent alors une 
relation d'équivalence sur $\tlgl(\rmF)$, 
moins fine que la 
relation de conjugaison par $G(\rmF)$, 
où $X,Y \in \tlgl(\rmF)$ sont dans la même classe si et 
seulement si $A_{i}(X) = A_{i}(Y)$ et 
$B_{j}(X) = B_{j}(Y)$ 
pour $i=0,\ldots, n$  et $j = 1, \ldots, n$. 
 Notons 
 $\mathcal{O}$ l'ensemble des classes d'équivalence pour 
 cette relation. Cette 
 relation pour des éléments semi-simples réguliers coïncide 
 avec la relation de conjugaison sous $G(\rmF)$ comme il est démontrée dans 
 \cite{rallSchiff}, proposition 
 6.2.

\blem\label{orbRegConj}
 Soient $X $ et $ Y $ deux éléments semi-simples réguliers de 
$\tlgl(\rmF)$. Ils appartiennent à la même classe dans 
$\mathcal{O}$ si et seulement s'ils 
sont conjugués par $G(\rmF)$.
\elem



Dans le paragraphe \ref{par:glnGlnplus1}, 
on a introduit pour tout $\tlP \in \calF(M_{\tlzero})$ 
des sous-espaces $\calV_{\tlP}$ et $\calZ_{\tlP}$ 
de $V \times V^{*}$. Voici leur rapport 
avec la décomposition (\ref{eq:xisamatrix}). 

\blem\label{shortLemmeOrbG}
 Soient $X = \matx{B}{u}{v}{d} \in 
 \tlgl$ et  
 $\tlP \in \calF(M_{\tlzero})$
Alors
\begin{enumerate}[i)]
\item $X \in \ml_{\tlP}$ 
si et seulement si 
$B \in \ml_{P}$ et
 $(u,v) \in \calZ_{\tlP}$.
\item $X \in \nl_{\tlP}$ 
si et seulement si 
$B \in \nl_{P}$ et
$(u,v) \in \calV_{\tlP}$.
\end{enumerate} 
 \elem

\dem Il s'agit d'un calcul matriciel direct. \bs

On étudiera maintenant 
les intersections des classes $\ol \in \calO$ 
avec les algèbres de 
Lie
de sous-groupes paraboliques relativement standards.

\brop\label{classIntUnipG}
 Soit $\tlP$ un sous-groupe parabolique relativement standard
 de $\tlG$. 
Alors, pour tous
$X \in \ml_{\tlP}(\rmF)$,
$N \in \nl_{\tlP}(\rmF)$ et $\ol \in \calO$ on a:
\begin{displaymath}
X \in \ol \iff X+N \in \ol.
\end{displaymath}
\erop

\dem La preuve s'appuie sur
le lemme \ref{shortLemmeOrbG} ci-dessus et 
est essentiellement identique à celle de la proposition 1.5 de \cite{leMoi}. \bs

\bcor\label{invsIntParG}
Soient $\tlP \in \relPb$ et 
$\ol \in \calO$. Alors
pour tout $A \subseteq \ml_{\tlP}(\rmF)$ 
et tout $B \subseteq \nl_{\tlP}(\rmF)$ on a:
\begin{equation*}
\ol \cap (A \oplus B) = (\ol \cap A) \oplus B.
\end{equation*}
\ecor


\subsection{Convergence du noyau modifié}\label{par:convergenceG}

On définit $\det \in \allts$ comme le déterminant du tore $A_{\tlzero}$ pour son action sur $W$. 
Notons que pour tout $x \in G(\A)$ et tout $\sigma \in \R$ 
on a alors $|\det x|_{\A}^{\sigma} = e^{\sigma \det (H_{G}(x)) }$.

Soit 
$f \in \calS(\tlgl(\A))$. 
Pour tout $\tlP \in \relPb$ et toute classe $\mathfrak{o} \in \mathcal{O}$ 
posons
\begin{equation}\label{noyauTroncPar}
k_{\tlP,\mathfrak{o}}(x) = k_{f, \tlP,\mathfrak{o}}(x) =
\sum_{\xi \in \ml_{\tlP}(\rmF) \cap \mathfrak{o}}
\int\limits_{
\mathrlap{\nl_{\tlP}(\A)}}
\ 
f(x^{-1}(\xi + U_{\tlP})x)dU_{\tlP}, \ 
x \in M_{P}(\rmF)N_{P}(\A) \bsl G(\A)
\end{equation} 
Pour $T \in \all_{\tlzero}$ on pose alors
\begin{equation}\label{noyauTronc}
k_{\mathfrak{o}}^{T}(x) = k_{f, \mathfrak{o}}^{T}(x) = 
\sum_{\tlP \in \relPb}
(-1)^{d_{\tlP}^{\tlG}}
\sum_{
\mathclap{\delta \in P(\rmF)\backslash G(\rmF)}}
\htau_{\tlP}(H_{\tlP}(\delta x)-T_{\tlP})
k_{\tlP,\mathfrak{o}}(\delta x), \ x \in G(\rmF)\backslash G(\A)
\end{equation}
où $d_{\tlP}^{\tlG}$ est défini par (\ref{eq:dPDefG}). N.B. la somme sur $\delta$ 
dans $P(\rmF)\backslash G(\rmF)$ est finie en vertu du lemme 5.1 de \cite{arthur3}.

\btheo\label{thm:MainConvG} 
 On a pour tout $T \in T_{+} + 
\all_{\tlzero}^{+}$ et tout
$\sigma \in \R$
\begin{displaymath}
\sum_{\mathfrak{o} \in \mathcal{O}} 
\int_{G(\rmF)\backslash  G(\A)}|k_{\mathfrak{o}}^{T}(x)|
|\det x|_{\A}^{\sigma} dx < \infty.
\end{displaymath}
\bdem 

En utilisant la proposition (\ref{harderNarashiman}) 
on a que $k_{\mathfrak{o}}^{T}(x)$ égale la somme sur $\tlP \in \relPb$ de $(-1)^{d_{\tlP}^{\tlG}}$ fois
la somme sur $\delta$ dans $P(\rmF)\backslash G(\rmF)$ de
\begin{equation*}
\displaystyle \left(
\sum_{\relPb \ni \tlP_{1} \subseteq \tlP}
\sum_{\eta \in P_{1}(\rmF)\backslash P(\rmF)}
F^{\tlP_{1}}(\eta\delta x,T_{\tlP_{1}})\tau_{\tlP_{1}}^{\tlP}
(H_{\tlP_{1}}(\eta \delta x)-T_{\tlP_{1}})
\right)
\htau_{\tlP}(H_{\tlP}(\delta x)-T_{\tlP})
k_{\tlP,\mathfrak{o}}(\delta x).
\end{equation*}
Suivant l'article \cite{arthur3}, les paragraphes 
6 et 7, on a:
\begin{displaymath}
\tau_{\tlP_{1}}^{\tlP}(H)
\htau_{\tlP}(H) = 
\sum_{\tlP_{2} \supseteq \tlP}
\sigma_{\tlone}^{\tltwo}(H) \quad 
H \in \mathfrak{a}_{\tlone}
\end{displaymath}
où
\begin{displaymath}
\sigma_{\tlP_{1}}^{\tlP_{2}}(H) = 
\sigma_{\tlone}^{\tltwo}(H) = 
\sum_{\tlQ \supseteq \tlP_{2}}
(-1)^{d_{\tltwo}^{\tlQ}}\tau_{\tlP_{1}}^{\tlQ}(H)
\htau_{\tlQ}(H),
\ H \in \mathfrak{a}_{\tlone}.
\end{displaymath}
Si l'on pose alors
\begin{align}
&\chi_{\tlone,\tltwo}^{T}(x) = 
F^{\tlone}(x,T)
\sigma_{\tlone}^{\tltwo}(H_{\tlone}(x)-T_{\tlone})\label{chippG}, \ 
x \in \tlP_{1}(\rmF)\backslash \tlG(\A),
\\
&k_{\tlP_{1},{\tlP_2},\mathfrak{o}}(x) =k_{\tlone,\tltwo,
\mathfrak{o}}(x) = \sum_{\tlP_{1} \subseteq \tlP \subseteq \tlP_{2}}
(-1)^{d_{\tlP}^{\tlG}}k_{\tlP,\mathfrak{o}}(x), \ 
x \in P_{1}(\rmF)\backslash G(\A),
\end{align}
on s'aperçoit alors que
\begin{equation}\label{eq:koAlternating}
k_{\mathfrak{o}}^{T}(x) = 
\sum_{\tlP_{1} \subseteq \tlP_{2}}
\sum_{\delta \in P_{1}(\rmF)\backslash G(\rmF)}
\chi_{\tlone,\tltwo}^{T}(\delta x) 
k_{\tlone,\tltwo,\mathfrak{o}}(\delta x)
\end{equation}
où la somme porte sur tous les couples $\tlP_{1}, \tlP_{2} \in \relPb$ tels que 
$\tlP_{1} \sbs \tlP_{2}$.

On fixe alors de tels $\tlP_{1} \subseteq \tlP_{2}$ et 
 l'on s'aperçoit qu'il suffit de démontrer:
\begin{equation}\label{firstSuffit}
\sum_{\mathfrak{o} \in \mathcal{O}}
\int_{P_{1}(\rmF)\backslash G(\A)}
\chi_{\tlone,\tltwo}^{T}(x)
|k_{\tlone,\tltwo,\mathfrak{o}}(x)|
e^{\sigma \det (H_{G}(x))}
dx < \infty.
\end{equation}

Si $\tlP_{1} = \tlP_{2} \neq \tlG$, c'est une conséquence 
du lemme 5.1 de \cite{arthur3} 
que  
$\sigma_{\tlone}^{\tltwo} \equiv 0$
donc l'intégrale (\ref{firstSuffit}) converge.
Si $\tlP_{1} = \tlP_{2} = \tlG$ 
les résultats de la section \ref{majorCrucSec}, 
le lemme \ref{lem:majCruc} en l'occurrence, 
montrent que 
 l'intégrale 
(\ref{firstSuffit}) est majorée par une somme sur $\tlB \in \relPb \cap \calP(M_{\tlzero})$ 
d'une intégrale d'une fonction continue
sur un ensemble $A_{\tlB}^{\tlG,\infty}(\tlG,c_{0},T_{\tlB})$ 
qui est compact donc on a aussi la convergence 
dans ce cas. Si $\tlP_{1} \sbn \tlP_{2}$ on démontrera 
un peu plus:

\btheo\label{thm:mainConv2} 
Soient $f \in \mathcal{S}(\tlgl(\A))$, $\sigma \in \R$
 et 
$\tlP_{1}, \tlP_{2}$ deux sous-groupes 
paraboliques relativement standards de $\tlG$
tels que $\tlP_{1}\subsetneq \tlP_{2}$. Alors
pour tout réel $\varepsilon_{0} > 0$ et tout $N \in \N$ 
il existe une 
constante $C$ qui ne dépend que de $N$, $f$, $\sigma$ et 
$\varepsilon_{0}$ telle que:
\begin{displaymath}
\sum_{\mathfrak{o} \in \mathcal{O}}
\int_{P_{1}(\rmF)\backslash G(\A)}
\chi_{\tlone,\tltwo}^{T}(x)
|k_{\tlone,\tltwo,\mathfrak{o}}(x)|
|\det x|_{\A}^{\sigma} dx
< Ce^{-N\|T\|}
\end{displaymath}
pour tout $T \in T_{+} + \all_{\tlzero}^{+}$ 
tels que $\forall \al \in \Delta_{\tlBmin} \ \al(T) > 
\varepsilon_{0}\|T\|$. 
\etheo

On fixe alors deux sous-groupes relativement standards 
$\tlP_{1} \subsetneq \tlP_{2}$. 
On introduit d'abord quelques notations.
Pour tous sous-groupes paraboliques relativement standards
$\tlQ \subseteq  \tlS$ posons:
\begin{displaymath}
(\bar\nl_{\tlQ}^{\tlS})' = \bar \nl_{\tlQ}^{\tlS}
\smallsetminus \bigcup_{ \tlQ \subseteq \tlR \subsetneq \tlS}
\bar \nl_{\tlQ}^{\tlR}, \quad 
\ml_{\tlS,\tlQ}' =
(\bar \nl_{\tlQ}^{\tlS})' \oplus \ml_{\tlQ} \oplus \nl_{\tlQ}^{\tlS}.
\end{displaymath}
On a alors les décompositions suivantes:
\begin{equation}\label{locFermLinG} 
\bar \nl_{\tlQ}^{\tlS} = \coprod_{Q\subseteq R\subseteq S}
(\bar\nl_{\tlQ}^{\tlR})'
\end{equation}
et 
$\ml_{\tlP} = \coprod_{\tlP_{1} \subseteq \tlS \subseteq \tlP}
(\ml_{\tlS,\tlone}'\oplus \nl_{\tlS}^{\tlP})$
pour tout $\tlP_{1} \subseteq \tlP \subseteq \tlP_{2}$.

Finalement, 
en utilisant le corollaire \ref{invsIntParG} on obtient, 
pour tout $\mathfrak{o} \in \mathcal{O}$ 
et tout $\tlP_{1} \subseteq \tlP \subseteq \tlP_{2}$
\begin{equation*}
\mathfrak{o} \cap \ml_{\tlP}(\rmF) = 
\coprod_{\tlP_{1} \subseteq \tlS \subseteq \tlP} 
(\mathfrak{o} \cap 
(\ml_{\tlS,\tlone}'\oplus \nl_{\tlS}^{\tlP})(\rmF)) = 
\coprod_{\tlP_{1} \subseteq \tlS \subseteq \tlP} 
((\mathfrak{o} \cap \ml_{\tlS,\tlone}'(\rmF)) 
\oplus \nl_{\tlS}^{\tlP}(\rmF)).
\end{equation*}
Grâce \`a cela, 
on peut réécrire $k_{\tlP,\mathfrak{o}}(x)$ comme:
\begin{displaymath}
\sum_{\tlP_{1} \subseteq \tlS \subseteq \tlP}
\sum_{\eta \in \nl_{\tlS}^{\tlP}(\rmF)}
\sum_{\zeta \in \ml_{\tlS,\tlone}'(\rmF) \cap \mathfrak{o}}
\int_{\nl_{\tlP}(\A)}f(x\inv (\eta + \zeta+ U_{\tlP})x)dU_{\tlP}.
\end{displaymath}
Fixons un caractère additif non-trivial $\psi$ sur 
$\rmF \backslash \A$. 
En appliquant la formule sommatoire 
de Poisson pour la somme portant sur
$\eta \in \nl_{\tlS}^{\tlP}(\rmF)$ 
de la fonction:
\begin{displaymath}
\nl_{\tlS}^{\tlP}(\A) \ni Y \longmapsto
\int_{\nl_{\tlP}(\A)}
f(x\inv (Y + \zeta  + U_{\tlP})x)dU_{\tlP}, 
\end{displaymath}
pour tout $\zeta \in \ml_{\tlS,\tlone}'(\rmF) \cap \mathfrak{o}$,
on obtient
\begin{displaymath}
k_{P,\mathfrak{o}}(x) = 
\sum_{\tlP_{1} \subseteq \tlS \subseteq \tlP}
\sum_{\eta \in \bar \nl_{\tlS}^{\tlP}(\rmF)}
\sum_{\zeta \in \ml_{\tlS,\tlone}'(\rmF) \cap \mathfrak{o}}
\Phi_{\tlS}(x,\zeta,\eta),
\end{displaymath}
o\`u:
\begin{displaymath}
\Phi_{\tlS}(x,X,Y) =\mathllap{
\int_{\mathrlap{\nl_{\tlS}(\A)}}}
f(x\inv (X + U_{\tlS})x)\psi 
(\langle Y,U_{\tlS}\rangle)dU_{\tlS}, \  x \in G(\A), 
X \in \ml_{\tlS}(\A), 
Y \in \bar \nl_{\tlS}^{\tltwo}(\A).
\end{displaymath}
En utilisant l'égalité (\ref{locFermLinG}) 
on peut écrire $k_{\tlP,\mathfrak{o}}(x)$ aussi comme:
\begin{displaymath}
\sum_{\tlP_{1} \subseteq \tlS \subseteq \tlR \subseteq  \tlP}
\sum_{\zeta \in \ml_{\tlS,\tlone}'(\rmF) \cap \mathfrak{o}}
\sum_{\eta \in (\bar \nl_{\tlS}^{\tlR})'(\rmF)}
\Phi_{\tlS}(x,\zeta,\eta).
\end{displaymath}
Grâce à cette formule, on a pour tout $\mathfrak{o} \in \mathcal{O}$:
\begin{multline}\label{doubleNoyauLongG}
k_{\tlone,\tltwo,\mathfrak{o}}(x) =  
\sum_{\tlP_{1} \subseteq \tlP \subseteq \tlP_{2}}
(-1)^{d_{\tlP}^{\tlG}}
k_{\tlP,\mathfrak{o}}(x) = \\
\sum_{\tlP_{1} \subseteq \tlS \subseteq \tlR \subseteq \tlP \subseteq \tlP_{2}}
(-1)^{d_{\tlP}^{\tlG}} 
\sum_{\zeta \in \ml_{\tlS,\tlone}'(\rmF) \cap \mathfrak{o}}
\sum_{\eta \in (\bar\nl_{\tlS}^{\tlR})'(\rmF)}
\Phi_{\tlS}(x,\zeta,\eta) = \\
\sum_{\tlP_{1} \subseteq \tlS \subseteq \tlR \subseteq  \tlP_{2}}
\sum_{\zeta \in \ml_{\tlS,\tlone}'(\rmF) \cap \mathfrak{o}}
\sum_{\eta \in( \bar\nl_{\tlS}^{\tlR})'(\rmF)}
\Phi_{\tlS}(x,\zeta,\eta)
\sum_{\tlR \subseteq \tlP \subseteq \tlP_{2}}
(-1)^{d_{\tlP}^{\tlG}}.
\end{multline}
On invoque maintenant l'identité due à Arthur \cite{arthur3}, 
proposition 1.1:
\begin{equation}\label{basicidentityG}
\sum_{\{ \tlP | \tlR \subseteq \tlP \subseteq \tlP_{2}\}}
(-1)^{d_{\tlP}^{\tlP_{2}}} = 
\begin{cases}
0 \text{ si } \tlR \neq \tlP_{2}, \\
1 \text{ sinon}.
\end{cases}
\end{equation}
On en déduit que la somme (\ref{doubleNoyauLongG}) décrivant 
$k_{\tlone,\tltwo,\mathfrak{o}}(x)$ se réduit à:
\begin{displaymath}
(-1)^{d_{\tlP_{2}}^{\tlG}}
\sum_{\tlP_{1} \subseteq \tlS \subseteq \tlP_{2}} 
\sum_{\eta \in (\bar\nl_{\tlS}^{\tltwo})'(\rmF)}
\sum_{\zeta \in \ml_{\tlS,\tlone}'(\rmF) \cap \mathfrak{o}}
\Phi_{\tlS}(x,\zeta,\eta).
\end{displaymath}
Ainsi, pour démontrer le théor\`eme \ref{thm:mainConv2} il suffit 
de majorer:
\begin{equation}\label{mainThmConv3G}
\int_{P_{1}(\rmF)\backslash G(\A)}
\chi_{\tlone,\tltwo}^{T}(x)
\sum_{\eta \in (\bar\nl_{\tlS}^{\tltwo})'(\rmF)}
\sum_{\zeta \in \ml_{\tlS,\tlone}'(\rmF)}
|\Phi_{\tlS}(x,\zeta,\eta)|e^{\sigma \det (H_{G}(x))}dx
\end{equation}
o\`u $\tlP_{1} \subseteq \tlS\subseteq \tlP_{2}$ sont fixés.  
Remarquons que la double somme sur $\mathfrak{o} \in \mathcal{O}$ 
et $\zeta \in 
 \ml_{\tlS,\tlone}'(\rmF) \cap \mathfrak{o}$ 
s'est réduit \`a la somme sur tout 
$\zeta \in  \ml_{\tlS,\tlone}'(\rmF)$.

En utilisant la décomposition d'Iwasawa $G(\A) = P_{1}(\A)K$, 
ainsi que la décomposition $P_{1} = N_{1}M_{1}$ on a 
que (\ref{mainThmConv3G}) égale

\begin{equation*}\label{secondSuffit}
\int\limits_{K}\int\limits_{[M_{1}]}
\int\limits_{[N_{1}]}
F^{\tlone}(m_{1},T)
\sigma_{\tlone}^{\tltwo}(H_{\tlone}(m_{1})-T_{\tlone})
e^{(\sigma \det-2\rho_{1})(H_{1}(m_{1}))}
\sum_{\mathclap{\eta \in (\bar\nl_{\tlS}^{\tltwo})'(\rmF)}} \quad
\sum_{\zeta \in \ml_{\tlS,\tlone}'(\rmF)}
|\Phi_{\tlS}(n_{1}m_{1}k,\zeta,\eta)|
dn_{1}dm_{1}dk.
\end{equation*}

Pour $\tlB \in \calP(M_{\tlzero}, B, \tlP)$ 
(voir définition (\ref{eq:relBrelToP})) soit 
$\la_{\tlB, \tlone, \sigma}$ l'élément de $(\all_{\tlB}^{\tlG})^{*}$ 
associé à $\upla = \sigma\det -2\rho_{1} \in \all_{1}^{*}$ 
par le lemme \ref{lem:majCruc}. En vertu de ce lemme 
il suffit de borner pour un tel $\tlB \sbs \tlP_{1}$  fixé 
l'expression suivante:

\begin{equation*}
\int\limits_{
\mathclap{A_{\tlB}^{\tlG,\infty}(P_{\tlone},c_{0},T_{\tlB})}}
\ \
\int_{[N_{1}]}
\sigma_{\tlone}^{\tltwo}(H_{\tlone}(a)-T_{\tlone})
e^{\la_{\tlB, \tlone, \sigma}(H_{\tlB}(a))}
\sup_{
k_{1} \in \omega_{\tlone}, 
k \in \tlK
}
\sum_{\eta \in (\bar\nl_{\tlS}^{\tltwo})'(\rmF)}
\sum_{\zeta \in \ml_{\tlS,\tlone}'(\rmF)}
|\Phi_{\tlS}(n_{1}k_{1}ak,\zeta,\eta)|dn_{1}da
\end{equation*}
où $\omega_{\tlone} = \omega_{\tlB} \cap M_{\tlone}(\A)$.

L'intégrale 
ci-dessus est essentiellement identique 
à celle qui apparaît dans la preuve 
de la proposition 4.4 dans \cite{chaud} 
qui dit justement qu'une telle expression 
vérifie les conditions du théorème \ref{thm:mainConv2}. 
Plus précisément, on voit qu'elle apparaît quand on passe de l'expression (4.8) à (4.9) 
dans loc. cit. On remarque que dans la preuve dans loc. cit. on a 
$\la_{\tlB, \tlone, \sigma}= 2\rho_{\tlB} - 2\rho_{\tlP_{1}}$, mais en fait 
la preuve marche sans changement pour n'importe quel $\la \in (\all_{\tlB}^{\tlG})^{*}$
ce qui 
démontre les théorèmes \ref{thm:mainConv2} et \ref{thm:MainConvG}. 
\edem
\etheo

%% file: quantitative_gln.tex
\section{Propriétés qualitatives}\label{quantSect}

On fixe une fois pour toute 
$\eta : \rmF \bsl \A^{*} \rar \C^{*}$ un caractère continu, qui est
trivial sur le groupe $\R_{> 0}^{*}$ vu comme un sous-groupe de 
$\A^{*}$ via l'inclusion $\R^{*}_{>0} \hrar (\rmF \otimes_{\Q}  \R)^{*} \hrar \A^{*}$. 
Notons alors pour $s \in \C$ et $x \in \A^{*}$:
\[
\eta_{s}(x) = |x|_{\A}^{s} \eta(x).
\]

Pour une fonction $f \in \mathcal{S}(\tlgl(\A))$ 
et $T \in T_{+} + \all_{\tlzero}^{+}$
on note
\begin{displaymath}
k^{T}(x) = 
k_{f}^{T}(x)  = \sum_{\mathfrak{o} \in \mathcal{O}}
k_{f,\mathfrak{o}}^{T}(x), \quad 
x \in G(\rmF)\backslash G(\A).
\end{displaymath}
Grâce au théorème \ref{thm:MainConvG} les 
distributions suivantes:
\begin{align*}
I_{\mathfrak{o}}^{T}(\eta_{s},f) &= 
\int_{G(\rmF)\backslash G(\A)}k_{f,\mathfrak{o}}^{T}(x)\eta_{s}(\det x)dx, &
\ \mathfrak{o} \in \mathcal{O}, \ T \in T_{+} + \all_{\tlzero}^{+},
\\ 
I^{T}(\eta_{s},f) &= \int_{G(\rmF)\backslash G(\A)}k_{f}^{T}(x) \eta_{s}(\det x)dx&
\end{align*}
sont bien définies pour tout $s \in \C$.

Dans le paragraphe \ref{par:asymptChaptG} 
on démontrera que pour $s \in \C$ la fonction
$T \mapsto I_{\mathfrak{o}}^{T}(\eta_{s},f)$ 
est un polynôme-exponentielle 
et si $s \neq -1,1$ son terme purement polynomial, 
noté $I_{\ol}(\eta_{s},f)$, ne dépend pas de $T$. 
Pour bien énoncer ce résultat on étudie d'abord les fonctions de type 
polynôme-exponentielle 
dans le paragraphe \ref{par:fonsPolExpG}.
et dans le paragraphe \ref{par:JTforLevisG} on introduit
 les distributions $I_{\ol}^{M_{\tlQ},T}(\eta_{s}, \cdot)$ pour tout sous-groupe 
parabolique relativement standard $\tlQ$ de $\tlG$.
La suite de cette section est consacrée 
aux propriétés des distributions $I_{\ol}(\eta_{s}, \cdot)$ pour 
$s \neq -1,1$.

\subsection{Polynômes-exponentielles}\label{par:fonsPolExpG}

Soit $\calV$ un $\R$-espace vectoriel de dimension finie. 
Par un polynôme-exponentielle sur $\calV$ on entend une fonction 
sur $\calV$ de la forme
\[
f(v) = \sum_{\la \in \calV^{*}}e^{\la(v)}P_{\la}(v), \quad v \in \calV
\]
où
$P_{\la}$ est un polynôme sur $\calV$ 
à coefficients complexes, égale $0$ 
pour presque tout $\la \in \calV^{*}$. 
On appelle $\la \in \calV^{*}$ tels que $P_{\la} \neq 0$ les exposants 
de $f$ et le polynôme correspondant à $\la = 0$ le terme 
purement polynomial de $f$.
On a alors le résultat d'unicité suivant: 
si $f$ est comme ci-dessus 
et $g = \sum_{\la \in \calV^{*}}e^{\la}Q_{\la}$ est un 
polynôme-exponentielle sur $\calV$ tel que 
$g(v) = f(v)$ pour tout $v \in \calV$ alors 
pour tout $\la \in \calV^{*}$ on a $P_{\la} = Q_{\la}$.

Pour $i = 0, 1, \ldots,n $ soit $e_{i}^{*} \in \all_{\tlzero}^{*}$ 
le caractère par lequel $A_{\tlzero}$ agit sur $D_{i}$. 
Posons aussi $e_{j}^{\vee} \in \all_{\tlzero}$ les 
éléments tels que $e_{i}^{*}(e_{j}^{\vee}) = \delta_{ij}$ où 
$i,j = 0,1, \ldots, n$.
On pose pour $i=1,\ldots, n$
\[
\tlvpi_{i}^{-} = \dfrac{n+1-i}{n+1}(\sum_{j=1}^{i}e_{j}^{*}) - 
\dfrac{i}{n+1}(e_{0}^{*} + \sum_{j=i+1}^{n} e_{i}^{*}), \quad 
\tlvpi_{i}^{+} = \dfrac{n+1-i}{n+1}(e_{0}^{*} + \sum_{j=1}^{i-1}e_{j}^{*}) - 
\dfrac{i}{n+1}(\sum_{j=i}^{n} e_{i}^{*}).
\]
On définit $\tlvpi_{i}^{-,\vee}$, $\tlvpi_{i}^{+, \vee}$ en remplaçant $^{*}$ par $^{\vee}$.
Alors $\tlvpi_{i}^{-}, \tlvpi_{i}^{+} \in (\all_{\tlzero}^{\tlG})^{*}$. 
On pose aussi $\tlvpi_{l}^{-} = \tlvpi_{l}^{+} = 0$ pour 
$l \nin \{1, \ldots, n\}$.

Fixons un sous-groupe parabolique relativement standard $\tlQ$ 
de $\tlG$ stabilisant le drapeau
\begin{equation*}
0 = V_{i_{0}} \subsetneq \cdots \subsetneq  V_{i_{k-1}}
\subsetneq 
V_{i_{k}} \oplus D_{0} \subsetneq \cdots \subsetneq V_{i_{l}} 
\oplus D_{0} =W.
\end{equation*}
Alors
\begin{gather*}
\hDelta_{\tlQ} = \{\tlvpi_{i_{a}}^{-}, \, \tlvpi_{i_{b}+1}^{+} | \, 
1 \le a \le k-1, \ k \le b \le l-1 \}, \\
\hDelta_{\tlQ}^{\vee} = \{\tlvpi_{i_{a}}^{-,\vee}, \, \tlvpi_{i_{b}+1}^{+,\vee} | \, 
1 \le a \le k-1, \ k \le b \le l-1 \}.
\end{gather*}

Posons 
$\tlvpi_{\tlQ}^{-} := \tlvpi^{-}_{m}$ où 
$m = \max(\{ j | \tlvpi_{j}^{-} \in \hDelta_{\tlQ} \} \cup \{0\})$. 
De même, on pose $\tlvpi_{\tlQ}^{+} := \tlvpi^{+}_{m}$ 
où $m = \min(\{j | \tlvpi_{j}^{+} \in \hDelta_{\tlQ} \} \cup \{0\})$.
Pour $s \in \C$ posons: 
\begin{equation}\label{eq:sQ}
s_{\tlQ} := \dfrac{s(n+1) + i_{k-1} + i_{k} - n}{i_{k}-i_{k-1} +1} 
\end{equation}
et
\begin{equation}\label{eq:uplaQs}
\upla_{\tlQ,s} := (1 + s_{\tlQ})\tlvpi_{\tlQ}^{-} + 
(1 - s_{\tlQ})\tlvpi_{\tlQ}^{+} \in (\all_{\tlQ,\C}^{\tlG})^{*}.
\end{equation}

Avec la notation du paragraphe \ref{par:HQtlQ} on a:
\blem\label{lem:HbarrhosQ} Soient $s \in \C$ et $\tlQ \in \relPb$. Alors
\begin{enumerate}
\item 
\begin{equation*}
\upla_{\tlQ,s} = 
\iota_{\tlQ}^{st}(s\det - 2\rho_{Q}) + 
2\rho_{\tlQ}
\end{equation*}
où l'on voit $s\det - 2\rho_{Q}$ comme l'élément de 
$(\all_{\tlQ,\C}^{st})^{*} $ par restriction. 
\item Pour tout $m \in H_{\tlQ}(\A)^{1} \times G_{\tlQ}(\A)$ on a
\[
e^{- \upla_{\tlQ,s}(H_{\tlQ}(m))} |\det m |^{s}_{\A}= |\det m|_{\A}^{s_{\tlQ}}.
\]
\item
Pour tout $\tlR \in \relPb$ tel que $\tlR \sps \tlQ$, la restriction de $\upla_{\tlQ,s}$ à $\tlR$ 
égale $\upla_{\tlR,s}$.
\end{enumerate}
\bdem
Calcul direct.
\edem
\elem

\blem\label{lem:uplasNonNul}
 Soit $\tlQ \in \relPb$.
 \begin{enumerate}[i)]
 \item Pour tout $s \in \C \smin \{-1,1\}$ et tout 
 $\tlvpi^{\vee} \in \hDelta_{\tlQ}^{\vee}$ 
on a $\upla_{\tlQ,s}(\tlvpi^{\vee}) \neq 0$.
\item Pour tout $s \in \C$ tel que $-1 < \Rel(s) < 1$ 
et tout  $\tlvpi^{\vee} \in \hDelta_{\tlQ}^{\vee}$ 
on a $\Rel(\upla_{\tlQ,s}(\tlvpi^{\vee})) > 0$.
 \end{enumerate}

\bdem 
On a pour $1 \le a \le k-1$ que 
$\upla_{\tlQ,s}(\tlvpi_{i_{a}}^{-, \vee}) = i_{a}(1+s)$ 
et pour $k \le b \le l-1$ que
$\upla_{\tlQ,s}(\tlvpi_{i_{b}+1}^{+, \vee}) = (n - i_{b})(1-s)$, 
d'où les résultats voulus. 
\edem
\elem

On note le corollaire immédiat du lemme \ref{lem:uplasNonNul}. 
\bcor\label{cor:uplasNonNul}
 Soit $\tlQ$ un sous-groupe parabolique relativement standard de $\tlG$.
Alors, pour tout sous-groupe parabolique $\tlR \sps \tlQ$ différent de $\tlG$ 
et tout 
$s \in \C \smin \{-1,1\}$ la restriction de 
$\upla_{\tlQ,s}$ à $\all_{\tlR}^{\tlG}$ est non-nulle. 
\ecor

Soit 
 $v_{\tlQ}$ 
 le volume dans $\all_{\tlQ}^{\tlG}$ du 
 parallélotope engendré par $(\hDelta_{\tlQ})^{\vee}$.
 Suivant le paragraphe 2 de  \cite{arthur2}, posons
\begin{equation}\label{eq:thetaHatDefG}
\hat \theta_{\tlQ}(\mu) = 
v_{\tlQ}^{-1}
\prod_{\tlvpi \in \hDelta_{\tlQ}^{\vee}} \mu(\tlvpi^{\vee}), 
\quad \mu \in \all_{\tlQ,\C}^{*}.
\end{equation}

Supposons $\tlR \sps \tlQ$, pour 
$X \in \all_{\tlQ}$ on note $X_{\tlR}$ sa projection à $\all_{\tlR}$ 
selon la décomposition (\ref{eq:decomp}). 
Suivant loc. cit., posons
\begin{equation}\label{eq:GammaQDefG}
\Gamma_{\tlQ}'(H,X) = 
\sum_{\tlR \sps \tlQ}(-1)^{d_{\tlQ}^{\tlR}}
\htau_{\tlR}(H_{\tlR}-X_{\tlR})\tau_{\tlQ}^{\tlR}(H), \quad H,X \in \all_{\tlQ}.
\end{equation}

\blem\label{lem:pQExplicitG}
 Soit $\tlQ \in \relPb$.
Alors, pour tout $s \in \C$ 
\[
p_{\tlQ,s}(X) := \int_{\all_{\tlQ}^{st}}e^{(s\det + 2\rho_{\tlQ} -2\rho_{Q})(H)}\Gamma_{\tlQ}'(H ,X)dH, \quad 
X \in \all_{\tlQ}
\]
est un polynôme-exponentielle sur $\all_{\tlQ}/\all_{\tlG}$. 
En plus, si $s \neq -1,1$, 
pour tout $\tlR \sps \tlQ$ 
il existe un polynôme $P_{\tlQ,\tlR,s}$ de degré au plus $d_{\tlQ}^{\tlG}$ sur $\all_{\tlR}/\all_{\tlG}$ 
tel que
\[
p_{\tlQ,s}(X) := j_{\tlQ}^{-1}\sum_{\tlR \sps \tlQ}e^{\upla_{\tlR,s}(X_{\tlR})}
p_{\tlQ,\tlR,s}(X_{\tlR})
\]
où $p_{\tlQ,\tlG,s}(X_{\tlG}) = (-1)^{d_{\tlQ}^{\tlG}}\hat \theta_{\tlQ}(\upla_{\tlQ,s})^{-1}$.
En particulier, si $s \neq -1,1$, la fonction $p_{\tlQ,s}$ est un polynôme-exponentielle
dont le terme purement polynomial est constant 
et égale $(-1)^{d_{\tlQ}^{\tlG}}\hat \theta_{Q}(\upla_{\tlQ,s})^{-1}$.

\brem On ne prétend pas que les polynômes $p_{\tlQ,\tlR,s}$ 
sont uniquement déterminés pour tout $\tlR \sps \tlQ$. En effet, il 
arrive que $\upla_{\tlR,s} = \upla_{\tlR',s}$ pour $\tlR \neq \tlR'$. Cependant, $\tlG$ 
est le seul sous-groupe parabolique $\tlR \sps \tlQ$
tel que $\upla_{\tlR,s} = 0$ si $s \neq -1,1$ d'où l'unicité du
terme $p_{\tlQ,\tlG,s}$.
\erem
\bdem

En utilisant le lemme \ref{lem:HbarrhosQ} et l'équation (\ref{eq:HQtlQ}) on a
\[
p_{\tlQ,s}(X) = j_{\tlQ}^{-1}
\int_{\all_{\tlQ}^{\tlG}}e^{\upla_{\tlQ,s}(H)}\Gamma_{\tlQ}'(H ,X)dH.
\]

Il résulte du lemme 2.1 de \cite{arthur2} que 
pour un $X$ fixé, la fonction $\all_{\tlQ}^{\tlG} \ni H \mapsto \Gamma_{\tlQ}'(H,X)$ 
est une fonction caractéristique d'un compact dans $\all_{\tlQ}^{\tlG}$. 
L'intégrale ci-dessus est donc bien définie et 
le lemme 
2.2  dans loc. cit. montre que c'est un polynôme-exponentielle ce qui 
démontre 
la première partie du lemme. 

Quand $s \neq -1,1$, on applique le lemme \ref{lem:uplasNonNul} \textit{i)} et son corollaire
\ref{cor:uplasNonNul} et 
on voit qu'on est dans la même situation que dans le lemme 3.3 de \cite{leMoi}. 
En utilisant ce lemme on conclut. 
\edem
\elem

\subsection{Une généralisation du théorème \ref{thm:MainConvG}}\label{par:JTforLevisG}

Soient $W'$ un $\rmF$-espace vectoriel de dimension $m+1$, 
$V' \sbs W'$ un sous-espace de dimension $m$ et $D_{0}' \sbs (W' \smin V') \cup \{0\}$  
une droite, où $m \in \N$.
Soient $H \cong \prod_{i=1}^{k}\Gl_{n_{i}}$ 
$G' = \Gl(V')$ et $\tlG' = \Gl(W')$ où $k \in \N$, et 
$n_{i} \in \N$ pour $i = 1,\ldots, k$.
On identifie $G'$ avec le sous-groupe de $\tlG'$ qui agit trivialement 
sur $D'_{0}$.
On va généraliser le théorème \ref{thm:MainConvG} au cas de l'inclusion 
$H \times G' \hrar H \times \tlG'$. 

Notons $\hl = \Lie(H)$, $\gl' = \Lie(G')$ 
et $\tlgl' = \Lie(\tlG')$. 
Pour $X \in (\hl \times \tlgl')(\rmF)$ soient $X_{1} \in \hl(\rmF)$ 
et $X_{2} \in \tlgl'(\rmF)$ tels que $X = X_{1}  + X_{2}$. 
Soit $\calO^{H \times \tlG'}$ la relation d'équivalence sur $(\hl \times \tlgl')(\rmF)$ 
définie de la façon suivante. On a $X = X_{1} + X_{2} \sim 
Y = Y_{1} + Y_{2}$ si et seulement si 
les polynômes caractéristiques de $X_{1}$ et $Y_1$ 
coïncident et si
$X_{2}$ et $Y_{2}$ sont dans la même classe 
pour 
la relation d'équivalence dans $\tlgl'(\rmF)$ 
décrite dans le paragraphe 
\ref{par:lesInvs} par rapport à l'inclusion $G' \hrar \tlG'$.

Soit $B$ un sous-$\rmF$-groupe de Borel
de $H \times G'$ et fixons aussi $M_{0}$ une partie de Levi de $B$. 
Soit $M_{\tlzero}$ l'unique sous-groupe de Levi minimal de $H \times \tlG'$ 
tel que $M_{\tlzero} \sps M_{0}$.
On peut alors parler de sous-groupes paraboliques standards de $H \times G'$ 
et semi-standards de $H \times \tlG'$.
Notons $\calF_{H \times G'}(M_{\tlzero}, B)$ le
sous-ensemble de sous-groupes paraboliques semi-standards $\tlP$ de $H \times \tlG'$ 
tels que $\tlP \sps B$.

Fixons un sous-groupe de Borel $\tlBmin \in \calF(M_{\tlzero})$. 
Soit $\tlP \in \calF(M_{\tlzero})$. Pour tout 
$H \in \all_{\tlBmin}$ on note $H_{\tlP}$ la projection de 
$sH$ à $\all_{\tlP}$ où $s$ est un élément du groupe de Weyl de 
$H \times \tlG'$ tel que $s^{-1}\tlP \sps \tlBmin$.

Pour une fonction 
$f \in \mathcal{S}((\hl \times \tlgl')(\A))$,  
un $\tlP \in \calF_{H \times G'}(M_{\tlzero}, B)$
et une classe $\ol \in \calO^{H \times \tlG'}$ on pose
\[
k_{f,\tlP, \ol}(x) = 
\sum_{\xi \in \ml_{\tlP}(\rmF) \cap \ol}
\int_{\nl_{\tlP}(\A)}
f(x^{-1}(\xi + U_{\tlP})x)dU_{\tlP}, \quad 
x \in  (H \times \tlG')(\A).
\]
Pour un $T \in \all_{\tlBmin}^{+}$ on pose donc
\[
k^{T}_{f,\ol}(x) = 
\sum_{\tlP \in \calF_{H \times G'}(M_{\tlzero}, B)} (-1)^{d_{\tlP}^{H \times \tlG'}}
\sum_{\delta \in P(\rmF) \bsl (H \times G')(\rmF)}
\htau_{\tlP}^{H \times \tlG'}(H_{\tlP}(\delta x)-T_{\tlP})
k_{f,\tlP,\ol}(\delta x)
\]
où $P = \tlP \cap (H \times G')$.

\btheo Soit $f \in \mathcal{S}((\hl\times \tlgl')(\A))$, 
alors pour tout $T \in  \all_{\tlBmin}^{+}$ 
suffisamment régulier et tout $\sigma \in \R$ on a
\[
\sum_{\ol \in \calO^{H \times \tlG'}}
\int_{(H \times G')(\rmF) \bsl H(\A)^{1} \times G'(\A)}
|k_{f,\ol}^{T}(x)|
|\det x|_{\A}^{\sigma}
dx < \infty.
\]
\bdem
La preuve est similaire à celle du théorème \ref{thm:MainConvG}.
Les détails sont laissés au lecteur.
Notons que si $f = f_{1} \otimes f_{2}$ 
où $f \in \calS(\hl(\A))$ et $f_{2} \in \calS(\tlgl'(\A))$, 
c'est une conséquence immédiate des théorèmes \ref{thm:MainConvG} ci-dessus
et 3.1 de \cite{chaud}.
\edem 
\etheo

Notons alors pour $s \in \C$, $\ol \in \calO^{H \times \tlG'}$ 
et $f \in \mathcal{S}((\hl\times \tlgl')(\A))$
\[
I_{\ol}^{H \times \tlG',T}(\eta_{s},f) =\int_{(H \times G')(\rmF) \bsl H(\A)^{1} \times G'(\A)}
k_{f,\ol}^{T}(x) \eta_{s}(\det x)dx. 
\]


Revenons dans le contexte de l'inclusion $G \hrar \tlG'$.
Soit $\tlQ$ un sous-groupe parabolique relativement standard de 
$\tlG$. 
Comme il est expliqué dans le paragraphe \ref{par:glnGlnplus1}, 
on a les décompositions 
 $M_{\tlQ} \cong H_{\tlQ} \times \tlG_{\tlQ}$ et 
$M_{Q} \cong H_{\tlQ} \times G_{\tlQ}$ où $H_{\tlQ}$, $G_{\tlQ}$ 
et $\tlG_{\tlQ}$ vérifient les conditions de ce paragraphe. 

Soit $\ol \in \mathcal{O}$, il existe 
$\ol_{\tlQ,1}, \ldots, \ol_{\tlQ,m} \in \calO^{M_{\tlQ}}$,
où $0 \le m < \infty$,
tels que
\begin{equation}\label{eq:mltlQcapolG}
\ml_{\tlQ}(\rmF) \cap \mathfrak{o} = 
\coprod_{i=1}^{m} \ol_{\tlQ,i} \cap \ml_{\tlQ}(\rmF).
\end{equation}
Pour $T \in \all_{\tlzero}^{+}$ et $s \in \C$, on définit alors les distributions $I_{\ol}^{M_{\tlQ},T}(\eta_{s}, \cdot)$ 
et $I^{M_{\tlQ},T}(\eta_{s}, \cdot)$ sur 
$\calS(\ml_{\tlQ}(\A))$ par:
\begin{equation}\label{eq:indDistDefG}
I^{M_{\tlQ},T}_{\mathfrak{o}}(\eta_{s}, f) = 
\sum_{i=1}^{m} 
I^{M_{\tlQ},T}_{\ol_{\tlQ,i}}(\eta_{s_{\tlQ}}, f), \quad
I^{M_{\tlQ},T}(\eta_{s}, f) = \sum_{\ol \in \calO}I^{M_{\tlQ},T}_{\mathfrak{o}}(\eta_{s}, f)
\end{equation}
où $s_{\tlQ}$ est défini par (\ref{eq:sQ}) et
pour $\ol_{\tlQ} \in \calO^{M_{\tlQ}}$, $I^{M_{\tlQ},T}_{\ol_{\tlQ}}(\eta_{s}, \cdot)$ 
c'est la distribution associée à l'inclusion $M_{Q} \hrar M_{\tlQ}$ 
décrite ci-dessus par rapport au sous-groupe de Levi minimal $M_{0}$ de $M_{Q}$
et aux sous-groupes 
de Borel $B \cap M_{Q}$ de $M_{Q}$ et $\tlBmin \cap M_{\tlQ}$ de $M_{\tlQ}$.

Pour $f \in \calS(\tlgl(\A))$ on 
pose
\begin{equation}\label{eq:fQdefG}
f_{\tlQ}(X) = \int_{K}\int_{\nl_{\tlQ}(\A)}f(k\inv(X+U_{\tlQ})k)\eta(\det k)dU_{\tlQ}dk, \quad 
X \in \ml_{\tlQ}(\A);
\end{equation}
alors $f_{\tlQ} \in \calS(\ml_{\tlQ}(\A))$. 
Notons que l'application 
\[
\tlQ \sps \tlP \mapsto M_{\tlQ} \cap \tlP
\]
définit une bijection entre les sous-groupes paraboliques relativement 
standards de $\tlG$
contenus dans $\tlQ$ et les sous-groupes paraboliques 
semi-standards de
$M_{\tlQ}$ contenant $B \cap M_{Q}$. 
En utilisant le lemme \ref{lem:HbarrhosQ}, on s'aperçoit alors que pour tout sous-groupe de Borel 
relativement standard $\tlB \sbs \tlQ$ et tous $T \in \all_{\tlzero}^{+}$
et $s \in \C$ on a
\begin{multline}\label{eq:JMQisThisG}
I^{M_{\tlQ},T}_{\mathfrak{o}}(\eta_{s}, f_{\tlQ}) = 
\int\limits_{\mathclap{M_{Q}(\rmF) \backslash H_{Q}(\A)^{1} \times G_{\tlQ}(\A)}} \qquad \quad
\sum_{i=1}^{m}k_{f_{\tlQ},\ol_{\tlQ,i}}^{T_{\tlB}}(m)\eta_{s_{\tlQ}}(\det m)dm  = 
\int\limits_{\mathclap{M_{Q}(\rmF) \backslash H_{Q}(\A)^{1} \times G_{\tlQ}(\A)}}
e^{-\upla_{\tlQ,s}(H_{\tlQ}(m))}
\sum_{\tlP \sbs \tlQ}(-1)^{d_{\tlP}^{\tlQ}} \\
\sum_{\mathclap{\eta \in (P \cap M_{Q})(\rmF) \bsl M_{Q}(\rmF)}} \ \
\htau_{\tlP}^{\tlQ}(H_{\tlP}(\eta m)-T_{\tlP})
\dsl
\sum_{\xi \in \ml_{\tlP}(\rmF) \cap \ol}
\int_{\nl_{\tlP}^{\tlQ}(\A)}
f_{\tlQ}(\Ad((\eta m)^{-1})(\xi + U_{\tlP}^{\tlQ}))dU_{\tlP}^{\tlQ}
\rb \eta_{s}(\det m)dm.
\end{multline}

\subsection{Le comportement en \texorpdfstring{$T$}{T}}\label{par:asymptChaptG}


On démontre la proposition suivante.

\brop\label{prop:mainQualitPropG}
Soient $f \in \calS(\tlgl(\A))$, $T' \in T_{+} + \all_{\tlzero}^{+}$, 
$s \in \C$, 
$\ol \in \calO$ et $T \in T' + \all_{\tlzero}^{+}$. Alors
\[
I_{\ol}^{T}(\eta_{s},f) = \sum_{\tlQ \in \relPb}
p_{\tlQ,s }(T_{\tlQ} - T'_{\tlQ})e^{\upla_{\tlQ,s}(T'_{\tlQ})}I_{\ol}^{M_{\tlQ}, T'}(\eta_{s},f_{\tlQ})
\]
où pour un sous-groupe parabolique $\tlQ$ 
relativement standard, 
la fonction $p_{\tlQ,s}$ est définie dans le lemme \ref{lem:pQExplicitG}, 
$\upla_{\tlQ,s } \in (\all_{\tlQ}^{\tlG})^{*}$ est défini 
par \ref{eq:uplaQs}, 
la distribution $I_{\ol}^{M_{\tlQ}, T'}$ 
est définie dans le paragraphe \ref{par:JTforLevisG} et 
$f_{\tlQ} \in \calS(\ml_{\tlQ}(\A))$ est définie par (\ref{eq:fQdefG}) 
dans le même paragraphe.

\bdem
Il est démontré dans le 
le paragraphe 2
de \cite{arthur2}, que les fonctions $\Gamma_{\tlQ}'$, 
définies par (\ref{eq:GammaQDefG}), 
vérifient la relation suivante: 
pour tout 
sous-groupe parabolique relativement standard \(\tlP\) de 
\(\tlG\), on a:
\begin{equation}\label{GammaPrimeRecurrenceG}
\htau_{\tlP}(H - X) = 
\sum_{\tlQ \supseteq \tlP}
(-1)^{d_{\tlQ}^{\tlG}}
\htau_{\tlP}^{\tlQ}(H)
\Gamma_{\tlQ}'(H,X), \quad H,X \in \all_{\tlP}.
\end{equation}

Fixons un \(T' \in T_{+}+\all_{\tlzero}^{+} \) 
et soit  $T \in T' + \all_{\tlzero}^{+}$.
En utilisant l'égalité ci-dessus 
dans la définition du noyau $k_{\ol}^{T}$  (\ref{noyauTronc})
avec 
$H = H_{\tlP}(\delta x)-T'_{\tlP}$ et 
$X = T_{\tlP} - T'_{\tlP}$ 
pour tout $\tlP \in \relPb$ et tout $\delta \in P(\rmF) \bsl G(\rmF)$ 
on a
\begin{multline}\label{eq:useGammaJolG}
I_{\ol}^{T}(\eta_{s},f)=
\int_{G(\rmF) \backslash G(\A)}
\sum_{\tlP \in \relPb}
(-1)^{d_{\tlP}^{\tlG}} 
\sum_{\delta \in P(\rmF)\backslash G(\rmF)}
\sum_{\tlQ \supseteq \tlP}(-1)^{d_{\tlQ}^{\tlG}}
\Psi^{T,T'}_{\tlP,\tlQ, \ol}(\delta x) \eta_{s}(\det x)dx = \\
\sum_{\tlQ \in \relPb}
\int\limits_{Q(\rmF)\backslash  G(\A)}
\sum_{\relPb \ni \tlP \subseteq \tlQ }
(-1)^{d_{\tlP}^{\tlQ}}
\sum_{\delta \in (P \cap M_{Q})(\rmF)\backslash M_{Q}(\rmF)}
\Psi^{T,T'}_{\tlP,\tlQ, \ol}(\delta x) \eta_{s}(\det x)dx
\end{multline}
où:
\begin{displaymath}
\Psi^{T,T'}_{\tlP,\tlQ, \ol}(x) = 
k_{\tlP,\ol}(x)
\htau_{\tlP}^{\tlQ}(H_{\tlP}(x)- T'_{\tlP})
\Gamma_{\tlQ}'(H_{\tlP}(x)-T'_{\tlQ},
T_{\tlQ}-T'_{\tlQ}).
\end{displaymath}

Le fait qu'on peut sortir la somme $\sum_{\tlQ}$ 
avant l'intégrale va se déduire du fait qu'on 
va montrer que les intégrales correspondantes 
sont absolument convergentes. 

Fixons $\tlQ \in \relPb$. 
On remplace l'intégrale sur $Q(\rmF)\backslash G(\A)$ 
par l'intégrale sur 
\[
N_{Q}(\rmF) \bsl N_{Q}(\A) \times A_{\tlQ}^{st, \infty} \times  
(M_{Q}(\rmF) \backslash H_{\tlQ}(\A)^{1}\times G_{\tlQ}(\A)) \times K
\]
ce qui donne $dx = e^{-2\rho_{Q}(H_{Q}(am))}dndadmdk$.

Soient $m \in (H_{\tlQ}(\A)^{1}\times G_{\tlQ}(\A))$, 
$a \in A_{\tlQ}^{st,\infty}$, 
$k \in K$, $\delta \in M_{Q}(\rmF)$ et $\tlP \sbs \tlQ$.
On a donc:
\begin{displaymath}
\int_{N_{Q}(\rmF)\backslash N_{Q}(\A)}
\Psi^{T,T'}_{\tlP,\tlQ, \ol}(\delta namk)dn = 
\int_{N_{Q}(\rmF)\backslash N_{Q}(\A)}
\Psi_{P,Q,\mathfrak{o}}^{T,T'}(na \delta mk)dn
\end{displaymath}
car $\delta \in M_{Q}(\rmF)$ 
normalise $N_{Q}(\A)$ sans changer sa mesure et 
il commute avec $A_{\tlQ}^{st,\infty} \sbs A_{Q}^{\infty}$. 
Les facteurs
de \(\Psi_{P,Q,\mathfrak{o}}^{T,T'}(na \delta mk)\) 
 deviennent:
\begin{gather*}
\Gamma_{\tlQ}'(H_{\tlQ}(na \delta mk)-T'_{\tlQ},
T_{\tlQ}-T'_{\tlQ}) = 
\Gamma_{\tlQ}'(H_{\tlQ}(a) + H_{\tlQ}(m) -T'_{\tlQ},
T_{\tlQ}-T'_{\tlQ}), \\
\htau_{\tlP}^{\tlQ}(H_{\tlP}(na\delta mk)- T'_{\tlP})=  
\htau_{\tlP}^{\tlQ}(H_{\tlP}(\delta m) + H_{\tlP}(a)- T'_{\tlP})=  
\htau_{\tlP}^{\tlQ}(H_{\tlP}(\delta m) - T'_{\tlP}).
\end{gather*}
Quant à $k_{\tlP,\mathfrak{o}}(na\eta m k)$,
on  fait le changement de variable
$(a\inv n\inv (\xi + U_{\tlP})na - \xi) \mapsto U_{\tlP}$  
(voir la définition de $k_{\tlP,\ol}$ au début du paragraphe
\ref{par:convergenceG}) et l'on obtient:
\begin{equation*}
k_{\tlP,\mathfrak{o}}(na\delta mk) =e^{2\rho_{\tlQ}(H_{\tlQ}(a))} 
\int_{\nl_{\tlP}(\A)} 
\sum_{\xi \in \ml_{\tlP}(\rmF) \cap \mathfrak{o}}
f((\delta mk)^{-1}(\xi + U_{\tlP} )\delta mk)dU_{\tlP} .
\end{equation*}

Ensuite, comme la mesure de \(N_{Q}(\rmF)\backslash N_{Q}(\A)\) 
vaut \(1\) on a 
en faisant le changement de variable 
$( (\delta m)^{-1}(\xi + U_{\tlQ})\delta m - \xi) \mapsto U_{\tlQ}$ 
\begin{multline}\label{eq:induction}
e^{-2\rho_{\tlQ}(H_{\tlQ}(a))}
\int\limits_{K}
\int\limits_{\mathrlap{[N_{Q}]}}
k_{\tlP,\mathfrak{o}}(na \delta mk) \eta(\det k)dndk \ \ \mathclap{=} 
\int\limits_{K} \ \ \int\limits_{\mathclap{\nl_{\tlP}(\A)}} \
\sum_{\mathrlap{\xi \in \ml_{\tlP}(\rmF) \cap \mathfrak{o}}}
f((\delta mk)^{-1}(\xi + U_{\tlP})\delta mk) \eta(\det k)dU_{\tlP}dk  \\
=
\int_{\nl_{\tlP}^{\tlQ}(\A)}\int_{K}\int_{\nl_{\tlQ}(\A)}\sum_{\xi \in \ml_{\tlP}(\rmF) \cap \mathfrak{o}}
f((\delta mk)^{-1}(\xi + U_{\tlP}^{\tlQ} + U_{\tlQ})\delta mk) \eta(\det k)
dU_{\tlQ}dkdU_{\tlP}^{\tlQ} =\\
= 
e^{2\rho_{\tlQ}(H_{\tlQ}(m))}
\int_{\nl_{\tlP}^{\tlQ}(\A)}\sum_{\xi \in \ml_{\tlP}(\rmF) \cap \mathfrak{o}}
f_{\tlQ}((\delta m)^{-1}(\xi + U_{\tlP}^{\tlQ})\delta m)dU_{\tlP}^{\tlQ}
\end{multline}
où $f_{\tlQ} \in \calS(\ml_{\tlQ}(\A))$ 
est définie par (\ref{eq:fQdefG}) dans le 
paragraphe \ref{par:JTforLevisG}. 

Remarquons qu'on a $2\rho_{\tlQ}(H_{\tlQ}(m)) = 2\rho_{Q}(H_{Q}(m))$.
En utilisant les lemmes \ref{lem:HbarrhosQ} et \ref{lem:pQExplicitG} ainsi que l'égalité 
(\ref{eq:HQtlQ})
 on voit donc que l'intégrale sur $A_{\tlQ}^{st,\infty}$ se réduit alors à
\begin{multline*}
\int\limits_{A_{\tlQ}^{st,\infty}}
\Gamma_{\tlQ}'(H_{\tlQ}(a) + H_{\tlQ}(m) -T'_{\tlQ},
T_{\tlQ}-T'_{\tlQ})e^{(s \det + 2(\rho_{\tlQ} - \rho_{Q}))(H_{\tlQ}(a))}da = 
e^{\upla_{\tlQ,s}(T'_{\tlQ}-  H_{\tlQ}(m))}p_{\tlQ,s}(T_{\tlQ}-T'_{\tlQ}).
\end{multline*}

En utilisant le calcul (\ref{eq:induction}) et 
en regardant la relation (\ref{eq:JMQisThisG}), 
on s'aperçoit qu'avec la notation de l'équation 
(\ref{eq:useGammaJolG}) on a
\begin{equation*}
\int_{Q(\rmF)\backslash  G(\A)}
\sum_{\tlP \subseteq \tlQ }
(-1)^{d_{\tlP}^{\tlQ}}
\sum_{\delta \in (P \cap M_{Q})(\rmF)\backslash M_{Q}(\rmF)}
\Psi^{T,T'}_{\tlP,\tlQ, \ol}(\delta x) \eta_{s}(\det x)dx = 
e^{\upla_{\tlQ,s}(T'_{\tlQ})}p_{\tlQ,s}(T_{\tlQ}-T'_{\tlQ})
I_{\ol}^{M_{\tlQ}, T'}(\eta_{s}, f_{\tlQ}).
\end{equation*}
Ce qu'il fallait démontrer.
\edem
\erop

En utilisant la proposition \ref{prop:mainQualitPropG} démontrée 
ci-dessus et le lemme \ref{lem:pQExplicitG} qui 
décrit les fonctions $p_{\tlQ}$ explicitement on obtient 
le comportement en $T$ des distributions $I_{\ol}^{T}$ et $I^{T}$.

\begin{theo}\label{thm:mainQualitThm}
Soit $f \in \calS(\tlgl(\A))$.
Les fonctions 
$T \mapsto I^{T}_{\mathfrak{o}}(\eta_{s},f)$ et
$T \mapsto I^{T}(\eta_{s},f)$ 
 où
$\mathfrak{o} \in \mathcal{O}$, $s \in \C$ et
$T$ 
parcourt $T_{+} + \all_{\tlzero}^{+}$ 
sont des polynômes-exponentielles. 
De plus, si $s \neq -1,1$ leur
parties purement polynomiales sont constantes 
et données respectivement par  
\begin{gather*}
I_{\ol}(\eta_{s},f) := \sum_{\tlQ \in \relPb}(-1)^{d_{\tlQ}^{\tlG}}j_{\tlQ}^{-1}\hat \theta_{\tlQ}(\upla_{\tlQ,s})^{-1}
e^{\upla_{\tlQ,s}(T'_{\tlQ})}I_{\ol}^{M_{\tlQ},T'}(\eta_{s},f_{\tlQ}), \\
I(\eta_{s},f) := \sum_{\tlQ \in \relPb}(-1)^{d_{\tlQ}^{\tlG}}j_{\tlQ}^{-1}\hat \theta_{\tlQ}(\upla_{\tlQ,s})^{-1}
e^{\upla_{\tlQ,s}(T'_{\tlQ})}I^{M_{\tlQ},T'}(\eta_{s},f_{\tlQ}),
\end{gather*}
pour tout $T' \in T_{+} + \all_{\tlzero}^{+}$. 
En particulier, les distributions $I_{\ol}$ et $I$ 
ne dépendent pas de $T'$.
\end{theo}

\brem\label{rem:JMQisPolExpG}
 Soit $\tlQ$ un sous-groupe parabolique relativement standard 
de $\tlG$. 
Par le même raisonnement que dans la proposition \ref{prop:mainQualitPropG}
 on obtient que pour tout $s \in \C$ les
distributions $I_{\ol}^{M_{\tlQ},T}(\eta_{s}, \cdot)$ et $I^{M_{\tlQ},T}(\eta_{s}, \cdot )$, 
définies dans le paragraphe \ref{par:JTforLevisG}, 
sont des polynômes-exponentielles en $T$ qui ne dépendent pas
de $T_{\tlQ} \in \all_{\tlQ}$. Cependant, 
si $\tlQ \neq \tlG$ le terme purement polynomial n'est pas constant. 
\erem

\brem\label{rem:Iolhomo} 
En vertu de la convergence absolue, pour tout $\tlQ \in \relPb$ et tout $f \in \calS(\ml_{\tlQ}(\A))$ on a que 
les fonctions $\C \ni s \mapsto I^{M_{\tlQ},T}(\eta_{s}, f)$ 
et $\C \ni s \mapsto I_{\ol}^{M_{\tlQ},T}(\eta_{s}, f)$, où $\ol \in \calO$, sont holomorphes.
Le théorème \ref{thm:mainQualitThm} dit alors que pour tout $f \in \calS(\tlgl(\A))$ 
et tout $\ol \in \calO$ les fonctions
$\C \smin \{-1,1\} \ni s \mapsto I^{T}(\eta_{s}, f)$ et 
$\C \smin \{-1,1\} \ni s \mapsto I_{\ol}^{T}(\eta_{s}, f)$ sont holomorphes
et admettent des prolongement méromorphes à $\C$ avec des pôles possibles 
en $-1$ et $1$.
\erem

\subsection{Équivariance}\label{par:compConjugG}

Soient $f \in \calS(\tlgl(\A))$ et $y \in G(\A)$. 
Notons $f^{y} \in \calS(\tlgl(\A))$ la fonction 
définie par $f^{y}(X) = f(\Ad(y)X)$. 

On voit que $I_{\ol}^{T}(\eta_{s},f^{y})$ pour 
$T \in T_{+} + \all_{\tlzero}^{+}$ et $s \in \C$ égale
\begin{equation*}
\int_{G(\rmF)\backslash G(\A)}
\sum_{\tlP \in \relPb}(-1)^{d_{\tlP}^{\tlG}}
\sum_{\delta \in P(\rmF)\backslash G(\rmF)}
k_{\tlP,\ol}(\delta x)\htau_{\tlP}(H_{\tlP}(\delta xy)-T_{\tlP})
\eta_{s}(\det (xy))dx.
\end{equation*}

Pour $x \in G(\A)$ et $P \in \calF(B)$ soit $k_{P}(x)$ un élément 
de $K$ tel que $xk_{P}(x)^{-1} \in P(\A)$. 
Alors, en utilisant l'égalité (\ref{GammaPrimeRecurrenceG}) 
on a:
\begin{equation*}
\htau_{\tlP}(H_{\tlP}(\delta x y)-T_{\tlP}) =
\sum_{\tlQ \supseteq \tlP}(-1)^{d_{\tlQ}^{\tlG}}
\htau_{\tlP}^{\tlQ}(H_{\tlP}(\delta x)-T_{\tlP})
\Gamma_{\tlQ}'(H_{\tlP}(\delta x)-T_{\tlP},
-H_{\tlP}(k_{P}(\delta x)y))
\end{equation*}
d'où on obtient que $I_{\ol}^{T}(\eta_{s},f^{y})$ égale la somme sur 
$\tlQ \in \relPb$ de
\begin{equation*}
 \int\limits_{\mathclap{Q(\rmF)\backslash G(\A)}}
\Gamma_{\tlQ}'(H_{\tlQ}(x)-T_{\tlQ}, -H_{\tlQ}(k_{Q}(x)y))
\sum_{\tlP \subseteq \tlQ}(-1)^{d_{\tlP}^{\tlQ}} 
\sum_{\mathrlap{\delta \in P(\rmF) \cap M_{Q}(\rmF) \backslash M_{Q(\rmF)}}} \
k_{\tlP,\ol}(\delta x)
\htau_{\tlP}^{\tlQ}(H_{\tlP}(\delta x)-T_{\tlP})
\eta_{s}(\det (xy))
dx.
\end{equation*}

Soit
$x = namk$ où $n \in N_{Q}(\rmF)\backslash N_{Q}(\A)$, 
$m \in M_{Q}(\rmF)\backslash (H_{\tlQ}(\A)^{1} \times G_{\tlQ}(\A))$, 
$a \in A_{\tlQ}^{st,\infty}$ et $k \in K$. Donc 
$dx = e^{-2\rho_{Q}(H_{Q}(am))}dndadmdk$
et pour $\delta \in M_{Q}(\rmF)$ on a
\begin{equation*}
\Gamma_{\tlQ}'(H_{\tlQ}(\delta namk)-T_{\tlQ},-H_{\tlQ}(k_{Q}(\delta namk)y))
 = 
\Gamma_{\tlQ}'(H_{\tlQ}(a) + H_{\tlQ}(m)-T_{\tlQ},-H_{\tlQ}(ky)).
\end{equation*}

Ensuite, en faisant les mêmes opérations 
comme dans (\ref{eq:induction}) au début de la preuve
de la proposition \ref{prop:mainQualitPropG} on s'aperçoit  
qu'on a pour $\tlP \subseteq \tlQ$,
$m \in M_{Q}(\rmF)\backslash  (H_{\tlQ}(\A)^{1} \times G_{\tlQ}(\A))$ 
et $\delta \in M_{Q}(\rmF)$
\begin{multline*}
\int\limits_{A_{\tlQ}^{st,\infty}}
\int\limits_{K}
\int\limits_{[N_{Q}]}
e^{(s\det-2\rho_{Q})(H_{Q}(am))}
k_{\tlP,\ol}(\delta nmak)
\Gamma_{\tlQ}'(H_{\tlP}(\delta namk)-T_{\tlQ},-H_{\tlP}(k_{P}(\delta namk)y))\eta(\det k)\\
 dndkda = |\det m|_{\A}^{s_{\tlQ}}  e^{\upla_{\tlQ,s}(T_{\tlQ})}
\int\limits_{\mathclap{\nl_{\tlP}^{\tlQ}(\A)}}
\sum_{\xi \in \ml_{\tlP}(\rmF) \cap \ol}
\int\limits_{K}
\int\limits_{\all_{\tlQ}^{\tlG}}
\int\limits_{\mathrlap{\nl_{\tlQ}(\A)}}
e^{\upla_{\tlQ,s}(H)}
f(k^{-1}((m\delta)^{-1}(\xi + U_{\tlP}^{\tlQ})\delta m + U_{\tlQ})k)\eta (\det k)\\
\Gamma_{\tlQ}'(H,-H_{\tlQ}(ky))
dU_{\tlQ}dHdkdU_{\tlP}^{\tlQ} = 
|\det m|_{\A}^{s_{\tlQ}}
 e^{\upla_{\tlQ,s}(T_{\tlQ})}
\int\limits_{\mathclap{\nl_{\tlP}^{\tlQ}(\A)}}
\sum_{\xi \in \ml_{\tlP}(\rmF) \cap \ol}
f_{\tlQ,s,y}((m\delta)^{-1}(\xi + U_{\tlP}^{\tlQ})m\delta)dU_{\tlP}^{\tlQ}
\end{multline*}
où l'on pose
\begin{equation*}
f_{\tlQ,s,y}(X) = 
\int_{K}
\int_{\nl_{\tlQ}(\A)}
f(k^{-1}(X + U_{\tlQ})k)u'_{\tlQ,s}(k,y) \eta (\det k)
dU_{\tlQ}dk, \quad X \in \ml_{\tlQ}(\A)
\end{equation*}
où 
\begin{equation*}
u'_{\tlQ,s}(k,y) = \int_{\all_{\tlQ}^{\tlG}}
e^{\upla_{\tlQ,s}(H)}\Gamma_{\tlQ}'(H,-H_{\tlQ}(ky))dH, \quad k \in K, \ s \in \C.
\end{equation*}
La fonction $K \ni k \mapsto u'_{\tlQ,s}(k,y)$ étant continue 
on a bien $f_{\tlQ,s,y } \in \calS(\ml_{\tlQ}(\A))$. 
On obtient le théorème suivant.

\begin{theo}\label{invarianceTheoG} Soient $y \in G(\A)$, $s \in \C$ et 
$f \in \calS(\tlgl(\A))$. Les distributions $I_{\ol}^{T}(\eta_{s}, \cdot)$ vérifient
\begin{equation*}
I_{\ol}^{T}(\eta_{s},f^{y}) - \eta_{s}(\det y) I_{\ol}^{T}(\eta_{s},f) = \eta_{s}(\det y)
\sum_{\tlQ \in \relPb \smin \{\tlG\}} 
e^{\upla_{\tlQ,s}(T_{\tlQ})}I_{\ol}^{M_{\tlQ},T}(\eta_{s}, f_{\tlQ,s,y})
\end{equation*}
où les distributions $I_{\ol}^{M_{\tlQ},T}(\eta_{s}, \cdot)$ 
sur $\calS(\ml_{\tlQ}(\A))$ 
sont définies par (\ref{eq:indDistDefG}).
 En particulier, pour $s \neq -1,1$ on a
\[
I_{\ol}(\eta_{s},f^{y}) = \eta_{s}(\det y) I_{\ol}(\eta_{s},f), \quad I(\eta_{s},f^{y}) = \eta_{s}(\det y)I(\eta_{s},f).
\]

\bdem 
La formule pour la différence $I_{\ol}^{T}(\eta_{s},f^{y}) - \eta_{s}(\det y) I_{\ol}^{T}(\eta_{s},f) $ est claire 
après les calculs qu'on a faits. Si $s \neq -1,1$, cette formule-ci démontre aussi l'$\eta_{s}$-invariance de 
la distribution $I_{\ol}(\eta_{s},\cdot)$, car
si $\tlQ \subsetneq \tlG$, d'après la remarque \ref{rem:JMQisPolExpG}, le terme 
$I_{\ol}^{M_{\tlQ},T}(\eta_{s}, f_{\tlQ,y,s})$ est un polynôme-exponentielle en $T$ qui ne dépend pas 
de $T_{\tlQ}$. En outre $\upla_{\tlQ,s}$ est non-trivial
sur $\all_{\tlQ}^{\tlG}$ en vertu du lemme \ref{lem:uplasNonNul}. 
Il en découle que $e^{\upla_{\tlQ,s}(T_{\tlQ})}I_{\ol}^{M_{\tlQ},T}(\eta_{s},f_{\tlQ,y,s})$ n'a pas de terme 
constant dans ce cas et par conséquent les termes constants de $I_{\ol}^{T}(\eta_{s},f^{y})$ 
et de $\eta_{s}(\det y)I_{\ol}^{T}(\eta_{s},f)$ coïncident. 
\edem
\end{theo}

\subsection{Indépendance des choix}\label{par:noChoixMadeG}

Soit $s \in \C \smin \{-1,1\}$ et $\ol \in \calO$.
Dans ce paragraphe on démontrera que la distribution $I_{\ol}(\eta_{s}, \cdot)$ 
ne dépend d'aucun choix, sauf le choix d'une mesure de Haar 
sur $G(\A)$ et les choix des mesures sur les 
sous-$\rmF$-espaces $\calV$ de $\tlgl$, notre choix étant que 
$\calV(\rmF)\bsl \calV(\A)$ soit de volume $1$.

Remarquons d'abord que $I_{\ol}(\eta_{s}, \cdot)$ 
ne dépend pas du choix du sous-groupe de Borel 
$\tlBmin \in \calP(M_{\tlzero})$ choisi au début du paragraphe 
\ref{par:reductionTheory}. En effet, ce choix intervient seulement dans le 
choix du chambre positive. Soit $\tlB \in \calP(M_{\tlzero})$ et $\sigma \in \Omega^{\tlG}$ tel 
que $\sigma\tlBmin = \tlB$. Notons $I_{\tlB, \ol}(\eta_{s}, \cdot)$ et $I_{\tlB, \ol}^{T}(\eta_{s}, \cdot)$ les distributions 
obtenues à partir de $\tlB$. Alors si $T \in \all_{\tlB}^{+}$, on a
$I_{\tlB, \ol}^{T}(\eta_{s}, \cdot) = I_{\ol}^{\sigma^{-1}T}(\eta_{s}, \cdot)$ 
et l'indépendance suit du théorème \ref{thm:mainQualitThm}.

Démontrons que $I_{\ol}(\eta_{s}, \cdot)$ ne dépend pas du choix 
du sous-groupe de Borel $B$ de $G$ contenant $M_{0}$. 
Soient $B' \in \calP(M_{0})$ et $\sigma \in \Omega^{G}$ 
tel que $B' = \sigma^{-1}B$. Notons $\tlB' = \sigma^{-1}\tlBmin$.
Notons $I_{B',\ol}^{T}(\eta_{s}, \cdot)$ et $I_{B',\ol}(\eta_{s}, \cdot)$ les distributions 
construites par rapport au sous-groupes de Borel
$B'$ et $\tlB'$.
C'est une conséquence simple de la relation 
(\ref{eq:weyl2}) qu'on a
$I_{\ol}^{T}(\eta_{s}, \cdot) = I_{B',\ol}^{\sigma^{-1}T}(\eta_{s}, \cdot)$.
Le théorème \ref{thm:mainQualitThm} permet de déduire de nouveau alors que
$I_{\ol}(\eta_{s}, \cdot) = I_{B',\ol}(\eta_{s}, \cdot)$.

Le compact maximal $\tlK$ qu'on a choisi dans le paragraphe \ref{par:glnGlnplus1} dépend du choix 
des vecteurs non-nuls dans les droites $D_{i}$ où $i=1,\ldots,n $. 
Tout autre compact $\tlK^{*}$ de ce type défini par rapport à $M_{\tlzero}$ 
égale donc $\gamma K \gamma^{-1}$ pour un $\gamma \in A_{B}(\rmF)$. 
Si l'on note alors $I_{\tlK^{*}, \ol}$ la distribution défini par rapport au mêmes données 
que $I_{\ol}$ mais pour le compact $\tlK^{*} = \gamma K \gamma^{-1}$ où $\gamma \in A_{B}(\rmF)$ on voit que 
$I_{\rmK^{*}, \ol}(\eta_{s}, f) = I_{\ol}(\eta_{s}, f^{\gamma})$ ce qui égale 
$I_{\ol}(\eta_{s}, f)$ en vertu du théorème \ref{invarianceTheoG} et du fait que 
$\eta_{s}$ est trivial sur $\rmF^{*}$.
La distribution $I_{\ol}(\eta_{s}, \cdot)$ ne dépend pas alors du choix de sous-groupe compact maximal 
"standard" par rapport à $M_{0}$.

Il nous reste à démontrer l'indépendance du choix de sous-groupe de Levi minimal de $G$
(car $M_{\tlzero}$ est déterminé par celui-ci).
Soit $M_{0}' = yM_{0}y^{-1}$ un sous-$\rmF$-groupe de Levi minimal de $G$, 
où $y \in G(\rmF)$. On note alors $I_{M_{0}', \ol}$ la distribution 
définie par rapport au sous-groupe 
parabolique minimal $yBy^{-1}$, le compact maximal $y\tlK y^{-1}$ de $\tlG(\A)$ 
et le sous-groupe de Borel $y\tlBmin y^{-1}$ de $\tlG$. 
On trouve alors $I_{M_{0}',\ol}(\eta_{s},f) = I_{\ol}(\eta_{s},f^{y})$ et le résultat 
découle du théorème \ref{invarianceTheoG} de nouveau.

\subsection{Orbites semi-simples réguli\`eres}\label{par:regOrbsChapG}

Soit $\mathcal{O}_{reg} \subseteq \mathcal{O}$ l'ensemble des 
orbites semi-simples régulières, \cad 
des orbites composées d'éléments 
semi-simples réguliers. 
La preuve de la proposition suivante est quasiment identique  
à la preuve de son homologue dans le paragraphe 3.6 de \cite{leMoi}.
\brop 
Pour tous $f \in \calS(\tlgl(\A))$, $T \in T_{+} + \all_{\tlzero}^{+}$, $s \in \C$, 
$\ol \in \calO_{reg}$
et $X \in\mathfrak{o}$ on a
\begin{displaymath}
I_{\mathfrak{o}}^{T}(\eta_{s},f) = I_{\mathfrak{o}}(\eta_{s},f) = 
\int_{G(\A)}f(x^{-1}Xx)\eta_{s}(\det x) dx
\end{displaymath}
où l'intégrale est absolument convergente.
\erop 

%% file: fourier_gln.tex
\section{Formule des traces infinitésimale}\label{sec:FourierTrans}

Il résulte de l'analyse faite dans le paragraphe \ref{par:lesInvs} qu'on 
a la décomposition de $\tlgl$ en sous-$\rmF$-espaces stables sous l'action de $G$ 
suivante:
\begin{equation}\label{eq:decBeteG}
\tlgl= \tl_{1} \oplus \tl_{2} \oplus \tl_{3}
\end{equation}
où $\tl_{1} = \gl$, $\tl_{2} = V \times V^{*}$ et $\tl_{3} = \Ga$. 
Soit $\tl \sbs \tlgl$ un sous-$\rmF$-espace défini 
comme une somme directe de certains d'entre les $\tl_{i}$, $i = 1,2,3$. 
Il y a donc huit possibilités pour $\tl$. 
Puisque chaque $\tl_{i}$ est $G$-stable et la restriction 
de la forme trace $\bilif$, 
à $\tl_{i}$ est non-dégénérée, il en est de même pour $\tl$. 
Pour $X \in \tlgl(\A)$ soit $X_{\tl}$ la projection 
de $X$ à $\tl(\A)$ selon la décomposition 
(\ref{eq:decBeteG}) ci-dessus.

Fixons $\psi$ 
un caractère non-trivial de $\rmF \bsl \A$. 
Pour $\tl$ comme ci-dessus, notons $\calF_{\tl}$ 
l'opérateur sur $\mathcal{S}(\tlul(\A))$ suivant
\[
\calF_{\tl}(f)(X) = \int_{\tl(\A)}f(X - X_{\tl} + Y_{\tl})
\psi (\langle X_{\tl}, Y_{\tl} \rangle )dY_{\tl}, \quad 
f \in \mathcal{S}(\tlgl(\A)), \ X \in \tlgl(\A)
\]
où $dY_{\tl}$ c'est la mesure de Haar sur $\tl(\A)$ pour 
laquelle le volume de $\tl(\rmF) \bsl \tl(\A)$ vaut $1$. 

\begin{theo}\label{thm:RTFJRIG} Pour tout
$f \in \mathcal{S}(\tlgl(\A))$ et tout $s \in \C \smin \{-1,1\}$ on a
\begin{displaymath}
\sum_{\ol \in \calO} 
I_{\ol}(\eta_{s}, f) = 
\sum_{\ol \in \calO} 
I_{\ol}(\eta_{s},\calF_{\tl}(f)).
\end{displaymath}

\bdem 
Soit $T \in T_{+} + \all_{\tlzero}^{+}$ et $s \in \C$ tel que $-1 < \Rel(s) < 1$.
On se place dans le contexte de la preuve du théorème \ref{thm:MainConvG}. 
En utilisant l'identité (\ref{eq:koAlternating}) et l'analyse qui suit l'équation (\ref{firstSuffit}) on a
que $\sum_{\ol \in \calO} 
I_{\ol}^{T}(\eta_{s}, f) -
\sum_{\ol \in \calO} 
I_{\ol}^{T}(\eta_{s},\calF_{\tl}(f)) = I^{T}(\eta_{s}, f) - I^{T}(\eta_{s}, \calF_{\tl}(f))$ vaut
\begin{multline*}
\int\limits_{[G]}
(k_{\tlG,\tlG}(x,f) - 
k_{\tlG,\tlG}(x, \calF_{\tl}(f)) + 
\sum_{\tlP_{1} \sbn \tlP_{2}}
\sum_{\delta \in P_{1}(\rmF) \backslash G(\rmF)} \\
\chi_{\tlP_{1},\tlP_{2}}^{T}(\delta x)
(k_{\tlP_{1},\tlP_{2}}(\delta x,f) - k_{\tlP_{1},\tlP_{2}}(\delta x,\calF_{\tl}(f))))\eta_{s}(\det x)dx
\end{multline*}
où pour les sous-groupes paraboliques relativement standards 
$\tlP_{1} \sbs \tlP_{2}$ et $\upphi \in \mathcal{S}(\tlgl(\A))$ 
on pose
\[
k_{\tlP_{1},{\tlP_2}}(x, \upphi) = 
\sum_{\tlP_{1} \subseteq \tlP \subseteq \tlP_{2}}
(-1)^{d_{\tlP}^{\tlG}}(\sum_{\ol \in \calO}k_{\upphi, \tlP, \ol}(x)), \ 
x \in P_{1}(\rmF)\backslash G(\A).
\]
On a alors pour tout $x \in G(\A)$
\begin{displaymath}
k_{\tlG,\tlG}(x,f) =
\sum_{\xi \in \tlgl(\rmF)}f(x\inv \xi x) = 
\sum_{\xi \in \tlgl(\rmF)}\calF_{\tl}(f)(x\inv \xi x) = 
k_{\tlG,\tlG}(x, \calF_{\tl}(f)) 
\end{displaymath}
grâce à la formule sommatoire de Poisson. 
On s'aperçoit alors que $\sum_{\ol \in \calO}I_{\ol}^{T}(\eta_{s}, f) -\sum_{\ol \in \calO}I_{\ol}^{T}(\eta_{s},\calF_{\tl}(f))$ 
c'est en fait la somme sur $\tlP_{1} \sbn \tlP_{2}$ 
de 
\[
\int_{P_{1}(\rmF) \bsl G(\A)}
\chi_{\tlP_{1},\tlP_{2}}^{T}(x)
(k_{\tlP_{1},\tlP_{2}}(x,f) - k_{\tlP_{1},\tlP_{2}}(x,\calF_{\tl}(f))))\eta_{s}(\det x)dx.
\]

Fixons $\varepsilon_{0} >0$.
En utilisant le théorème \ref{thm:mainConv2}
pour $f$ et $\calF_{\tl}(f)$ on a donc pour 
tout $N > 0$
\begin{displaymath}
|I^{T}(\eta_{s}, f) - I^{T}(\eta_{s}, \calF_{\tl}(f))|= O(e^{-N\|T\|})
\end{displaymath}
si $T \in T_{+} + \all_{\tlzero}^{+}$ est tel 
que $\forall \al \in \Delta_{\tlBmin}$, $\al(T) > \varepsilon_{0}\|T\|$. 
D'après la proposition \ref{prop:mainQualitPropG} 
la différence
$I^{T}(\eta_{s}, f) - I^{T}(\eta_{s}, \calF_{\tl}(f))$ 
est un polynôme-exponentielle en $T$ qui pour $-1 < \Rel(s) < 1$,
en vertu du lemme \ref{lem:uplasNonNul} \textit{b)}, est soit constant, soit sa norme
tend vers $\infty$ quand le paramètre $T \in \all_{\tlzero}^{+}$ est tel que 
$\al(T)$ tendent vers $\infty$ pour tout $\al \in \Delta_{\tlBmin}$. 
L'égalité ci-dessus 
implique alors $I^{T}(\eta_{s}, f) = I^{T}(\eta_{s}, \calF_{\tl}(f))$. 
En vertu de la remarque \ref{rem:Iolhomo} on obtient donc 
l'égalité des fonctions méromorphes en la variable $s \in \C$
$I^{T}(\eta_{s}, f) = I^{T}(\eta_{s}, \calF_{\tl}(f))$.
En particulier, en vertu du théorème \ref{thm:mainQualitThm}, on a égalité de leur 
termes constants pour $s \neq -1,1$.
 \edem
\end{theo}

%% file: semisimple_orbs.tex
\section{Orbites semi-simples}\label{sec:orbRrssSG}

Soit
$\mathcal{O}_{rs} \subseteq \mathcal{O}$ l'ensemble 
des classes contenant un élément $
\begin{pmatrix}
B & u \\
v & d \\
\end{pmatrix}$ tel que le polynôme caractéristique de 
$B$ est séparable (i.e. $B$ est semi-simple régulier dans $\gl(\rmF)$). 
On appelle de telles 
classes \textit{relativement semi-simples régulières}. 
Si $\mathfrak{o} \in \mathcal{O}_{rs}$, alors tout élément 
$
\begin{pmatrix}
B_{0} & u \\
v & d \\
\end{pmatrix} \in \mathfrak{o}$ a la propriété 
que $B_{0}$ soit semi-simple régulier.

Soit $\ol \in \mathcal{O}_{rs}$.
Le but de cette section est de donner une expression explicite 
pour $I_{\ol}(\eta_{s}, \cdot)$ ce qu'on achève par le théorème
\ref{thm:theThmOrbsG}. Les résultats sont analogues, et parfois même identiques, 
à ceux de la section 5 de \cite{leMoi} et 
par conséquent on renvoie souvent à cet article pour les preuves détaillées. 
En particulier, on omets les résultats de convergence, 
qui sont démontrés dans le paragraphe 5.7, ainsi que les résultats de prolongement 
méromorphe qui se trouvent dans le paragraphe 5.8 de loc. cit.
On utilisera aussi la même notation que dans la section 5 de loc. cit. ce qui rendra l'analogie plus visible.

Voici le plan de la section:
après avoir introduit quelques notations dans le paragraphe 
suivant \ref{par:noteEnsemG}, on décrit 
la décomposition de $\ol$ en $G(\rmF)$-orbites 
dans le paragraphe \ref{par:orbitesDansClasseG}. 
On introduit encore 
un peu plus de notation dans 
\ref{par:defsOrbsG}. Dans le paragraphe \ref{par:LeResultRssG}
on définit une expression $i_{\ol}(x)$ pour laquelle on a
\begin{equation}\label{eq:JolIsjolG}
I_{\ol}(\eta_{s}, f) = \int_{[G]} i_{\ol}(x) \eta_{s}(\det x)dx
\end{equation}
si $-1 < \Rel(s) < 1$.
En supposant cela, on
donne la preuve du théorème \ref{thm:theThmOrbsG} omettant 
les preuves des énoncés techniques. 
Dans la section \ref{par:noyuTronqNouvG} on introduit 
un nouveau noyau tronqué $i_{f,\ol}^{T}(x)$ 
tel que 
$\int_{[G]} i_{f,\ol}^{T}(x)\eta_{s}(\det x)dx = I_{\ol}^{T}(\eta_{s}, f)$ pour tout $s \in \C$.
Ce résultat nous permet de démontrer (\ref{eq:JolIsjolG}) 
 dans le paragraphe suivant 
\ref{par:IntReprDeJG}. 
On finit la preuve dans le 
dernier paragraphe \ref{par:holoResG} 
où on étudie certaines fonctions zêtas - les homologues 
des fonctions zêtas étudies dans \cite{leMoi}.

\subsection{Notations}\label{par:noteEnsemG}

On utilisera les lettres $I, J$, avec de possibles indices, 
pour noter des 
sous-ensembles finis de $\N^{*}$. 
Soit $I \sbs \N^{*}$ fini. On pose 
$-I = \bigcup_{i \in I} \{-i\}$.
On dit que $\calI$ est un $\epsilon$-sous-ensemble 
de $I$, si $\calI \sbs I \cup -I$ et si pour tout $i \in \calI$ 
on a $-i \nin \calI$.
Dans ce cas on 
écrit $\calI \sbse I$. La notation est 
un peu abusive car $\calI$ n'est pas forcément un 
sous-ensemble de $I$. 
On définit aussi 
$\acalI = \{|i| | i \in \calI \} \sbs \N^{*}$ 
et
$\calI^{\sharp} \sbse I$ 
par la propriété $\calI \cup \calI^{\sharp} = 
\acalI \cup -\acalI$. Autrement dit $\calI^{\sharp} = -\calI$ mais 
on écrira $\calI^{\sharp}$ dans ce contexte.
On réserve les lettres $\calI$, $\calJ$ et $\calK$ et seulement ces 
trois lettres avec de possibles indices, 
pour des $\epsilon$-sous-ensembles. 

On utilisera aussi la notation abrégée suivante
soit $I' \sbs \N^{*}$ et $\calI, \calJ \sbse I'$, 
on écrira $\calJ \cup \calI \sbse I'$ pour signifier 
que la réunion ensembliste $\calI \cup \calJ$ est aussi 
un $\epsilon$-sous-ensemble (ce qui n'est pas toujours vrai). 
On utilise le symbole $\sqcup$ pour noter la réunion
disjointe, donc $I \sqcup J = I'$ implique 
$I \cap J = \varnothing$.
On écrira aussi, pour $\calI \sbse I'$ fixés
\[
\sum_{\acalJ = I'} := 
\sum_{\begin{subarray}{c}
\calJ \sbse I' \\
\acalJ = I'
\end{subarray}}
, \ 
\sum_{\calK \sqcup \calJ \sbs \calI} := 
\sum_{\begin{subarray}{c}
\calK, \calJ \sbs \calI \\
\calK \cap \calJ = \varnothing
\end{subarray}}
, \ 
\sum_{\acalK \sqcup \acalJ = I'} := 
\sum_{\begin{subarray}{c}
\calK, \calJ \sbse I' \\
\acalK \sqcup \acalJ = I'
\end{subarray}}.
\]
Finalement, pour $\calJ_{1} \sqcup \calJ_{3} \sqcup \calJ_{4} \sbse I_{0}$ 
et $\calJ_{2} \sbs \calJ_{3}$ on utilisera parfois la notation suivante
\[
\calJ_{3 \smin 2} := \calJ_{3} \smin \calJ_{2}, \ 
\calJ_{13} := \calJ_{1} \cup \calJ_{3}, \ 
\calJ_{134} := \calJ_{1} \cup \calJ_{3} \cup \calJ_{4} \text{ etc}.
\]

\subsection{Orbites dans une classe relativement semi-simple régulière}\label{par:orbitesDansClasseG}

Dans ce paragraphe on décrit une décomposition en orbites 
d'une classe relativement semi-simple régulière. 

En utilisant la décomposition (\ref{eq:xisamatrix}) on va écrire les éléments $X$
de $\tlgl$ sous forme
\[
X = (B, \rmw, d)
\]
où $B \in \gl$, $\rmw \in V \times V^{*}$ et $d \in \Ga$. 
On considère alors $V \times V^{*}$ comme un espace vectoriel 
dont $V$ et $V^{*}$ sont des sous-espaces. 
L'action d'un $g \in G$ s'écrit donc simplement 
\[
\Ad(g)X = (\Ad(g)B, g\rmw, d)
\]
où $G$ agit sur $V \times V^{*}$ par ses actions naturelles sur $V$ et $V^{*}$.

Donnons-nous 
\begin{itemize}
\item Un polynôme séparable $Q \in \rmF[T]$ de degré $n$.
\item Une décomposition de $Q$ en facteurs irréductibles $Q = \prod_{i=1}^{m}Q_{i}$
et notons $I = \{1, \ldots, m\}$.
\item Pour tout $i \in I$ notons $\rmF_{i} = \rmF[T]/ (Q_{i}(T))$ (resp. $\rmF_{i} = \rmF[T]/(-1)^{\deg Q_{i}}Q_{i}(-T)$)
et $F_{I} = \rmF[T]/ (Q(T))$ (resp. $F_{I} = \rmF[T]/ (-1)^{n}(Q(T))$). On identifie $\rmF_{I}$ (resp. $\rmF_{-I}$)
à $\prod_{i \in I}\rmF_{i}$  (resp. à $\prod_{i \in I}\rmF_{-i}$)
à l'aide des homomorphismes des $\rmF$-algèbres $\rmF[T]/ (Q(T)) = \rmF_{I} \rar \rmF_{i} = \rmF[T]/ (Q_{i}(T))$ 
(resp. $\rmF_{-I} \rar \rmF_{-i}$)
induits par $T \mapsto T$.
\item On fixe l'isomorphisme des algèbres étales $\iota_{I}: \rmF_{I} \rar \rmF_{-I}$ 
induit par $T \mapsto -T$. On identifie $\rmF_{-I}$ à $\rmF_{I}^{*}$, en tant que $\rmF$-espace vectoriel 
grâce à la forme bilinéaire non-dégénérée $\rmF_{I} \times \rmF_{-I} \ni (u,v) \mapsto \Tr_{\rmF_{I}/\rmF}(u \iota_{I}(v))$.
\item On note $b_{I} \in \rmF_{I}$ l'image de $T$ dans $\rmF_{I}$. Dans ce cas $b_{I}$ engendre la $\rmF$-algèbre étale 
$\rmF_{I}$.
\end{itemize}

Soit alors $P_{I} = M_{I}N_{I}$ le sous-groupe parabolique standard de $G$ défini comme le stabilisateur 
du drapeau
\[
0 = V_{j_{0}} \sbn V_{j_{1}} \sbn V_{j_{2}} \sbn \cdots \sbn V_{j_{m}} = V
\]
où pour $i \in I$ on a $\deg Q_{i} = j_{i} - j_{i-1}$. Pour $i \in I$ posons $Z_{i} := \bigoplus_{j = j_{i-1} + 1}^{j_{i}}D_{j}$. 
On peut, et on fixe un $B \in \ml_{P_{I}}(\rmF)$ dont $Q_{i}$ est le polynôme caractéristique pour son action 
sur $Z_{i}$.
Alors, $B$ est semi-simple régulier dans $\gl(\rmF)$ 
et tout élément semi-simple régulier de $\gl(\rmF)$ s'obtient par une telle construction.

L'élément $B$ induit une structure de $\rmF[T]$-module sur $V$ par $T.v = Bv$. 
L'algèbre $\rmF_{I}$ est aussi naturellement un $\rmF[T]$-module qui est isomorphe, en tant 
que $\rmF[T]$-module à $V$. On fixe un tel isomorphisme et on identifie $V$ avec $\rmF_{I}$ désormais. Notons 
qu'un tel isomorphisme identifie aussi $\rmF_{i}$ avec $Z_{i}$ pour tout $i \in \rmI$. En plus, l'isomorphisme $V \cong \rmF_{I}$
induit 
l'isomorphisme duale entre $V^{*}$ et $(\rmF_{I})^{*} = \rmF_{-I}$. On identifie aussi ces espaces.
On a donc pour tous 
$(u,v) \in \FFi$, $g \in G(\rmF)$ et $k \in \N$:
\begin{gather*}
v(B^{k}u) = \Tr_{\rmF_{I}/\rmF}(b_{I}^{k} u \iota_{I}(v) ), \\
\Tr_{\rmF_{I}/\rmF}( gu \iota_{I}(gv) ) = \Tr_{\rmF_{I}/\rmF}(u \iota_{I}(v) ).
\end{gather*}

Pour tout $i \in I$ on note $1_{i}$ et $1_{-i}$ les unités de $\rmF_{i}$
et $\rmF_{-i}$ respectivement. Pour tout $\calI \sbse I$ 
on pose $\rmF_{\calI} = \bigoplus_{i \in I} \rmF_{i} \sbs \rmF_{I} \times \rmF_{-I}$ 
et $1_{\calI} = \sum_{i \in I}1_{i} \in \rmF_{\calI}^{*}$.

Soit $T_{I}$ le centralisateur de $B$ dans $G$. C'est un $\rmF$-tore maximal de $G$, et il 
est contenu dans $M_{I}$.
Pour $J \sbs I$ on note $T_{J} \sbs T_{I}$ le plus grand sous-tore 
qui agit trivialement sur $\rmF_{I \smin J}$. Si $J = \{i\}$ on écrit 
simplement $T_{i} = T_{\{i\}}$. Alors $T = \prod_{i \in I}T_{i}$ 
et pour tout $i \in I$ le groupe $T_{i}(\rmF)$ 
agit simplement transitivement sur $\rmF_{i}^{*}$ ainsi que sur $\rmF_{-i}^{*}$.

On fixe une classe $\ol \in \calO_{rs}$ contenant un élément $X$ de type 
$(B, \rmw, d)$. 
Remarquons que $d$ est le même pour tout $X \in \ol$, on le note simplement $d_{\ol}$.

On introduit l'ensemble 
$\rmV_{\ol} = \{(u,v) \in \rmF_{I} \times \rmF_{-I}|  \Tr_{\rmF_{I}/\rmF}(b_{I}^{i}u \iota_{I}(v)) = A_{i}(X)
\ \forall \ 0 \le i \le n-1\}$ où
$X \in \mathfrak{o}$ quelconque. 
Comme $T_{I}(\rmF)$ agit sur $F_{I}\times \rmF_{-I}$, commute à $B$ et préserve l'accouplement naturel 
sur $\rmF_{I} \times \rmF_{-I}$, il agit aussi sur $\rmV_{\ol}$.
On voit que l'ensemble des orbites dans $\rmV_{\ol}$ sous l'action 
de $T_{I}(\rmF)$
est en bijection avec l'ensemble des classes de $G(\rmF)$-conjugaison dans 
$\ol$, la bijection étant induite 
par l'application
$\rmV_{\ol} \ni \rmw \mapsto (B, \rmw, d_{\ol})$. 

Pour $\rmw \in \rmF_{I} \times \rmF_{-I}$ et $\calI \sbse I$ on pose 
$\rmw_{\calI} := 1_{\calI} \rmw$. Si $\calI = \{i\}$ on note simplement 
$\rmw_{i} = \rmw_{\{i\}}$. On va continuer pourtant d'écrire parfois les éléments de 
$\rmF_{I} \times \rmF_{-I}$ comme $(u,v)$.

\blem\label{lem:chaudNagginLemG} 
Il existe un $\al_{I} \in \rmF_{I}$ tel que 
pour tout $(u,v) \in \FFi$
on a
\[
(u,v) \in \rmV_{\ol} \Longleftrightarrow 
u \iota_{I}(v) = \al_{I}.
\]
\bdem 
Pour tous $(u,v), (u',v') \in \rmV_{\ol}$ et $k \in \N$ on a
\[
v(B^{k}u) = 
\Tr_{\rmF_{I}/\rmF}(b_{I}^{k}u \iota_{I}(v))= 
\Tr_{\rmF_{I}/\rmF}(b_{I}^{k}u' \iota_{I}(v')) =v'(B^{k}u')
\]
d'où
\[
\Tr_{\rmF_{I}/ \rmF}
( b^{k}_{I} (u \iota_{I}(v) - u' \iota_{I}(v') )) = 0, \quad 
\forall \ k \in \N.
\]
D'après la proposition (18.3) dans \cite{bInvolutions}, 
la forme
$\Tr_{\rmF_{I}/ \rmF}$ est non-dégénérée, et puisque
les puissances de $b_{I}$ 
engendrent $\rmF_{I}$ sur $\rmF$ on obtient 
\begin{equation}
u \iota_{I}(v) = u' \iota_{I}(v'), \quad
\forall \ (u,v), (u',v') \in \rmV_{\ol}.
\end{equation}
On pose donc $\al_{I} := u \iota_{I}(v)$ où $(u,v) \in \rmV_{\ol}$ 
quelconque. Il reste à démontrer 
que si $(u,v) \in \rmV_{\ol}$ est tel que 
$u \iota_{I}(v)  = \al_{I}$  
alors $(u,v) \in \rmV_{\ol}$. Pour cela il suffit de faire le même calcul dans le sens inverse.
\edem
\elem

Soit $\al_{I}  = (\al_{i})_{i \in I} \in \rmF_{I}$ comme dans le lemme précédant. 
Notons $I_{0} \sbs I$ l'ensemble de $i \in I$ tels 
que $\al_{i} = 0$.

\brop\label{prop:orbitsRRSSG}
Il existe une unique $T_{I}(\rmF)$-orbite 
dans $\rmV_{\ol}$ composée des $(u,v) \in \rmV_{\ol}$ 
tels que $u_{i} = 0$ et $v_{-i} = 0$ pour tout $i \in I_{0}$. On choisit 
$\xi_{\varnothing}$ 
un représentant de cette orbite. Alors,  
les $T_{I}(\rmF)$-orbites dans $\rmV_{\ol}$ 
sont en bijection avec les $\epsilon$-sous-ensembles de $I_{0}$, 
le représentant de l'orbite correspondant à $\calI \sbse I_{0}$ 
étant 
$\xi_{\varnothing} + 1_{\calI}$. Les orbites 
de dimension maximale correspondent aux $\calI \sbse I_{0}$ 
tels que $\acalI = I_{0}$. 

\bdem La preuve est identique à celle de la proposition 5.2 de \cite{leMoi}. 
\edem
\erop


\subsection{Quelques définitions associées aux orbites}\label{par:defsOrbsG}

D'après la proposition \ref{prop:orbitsRRSSG} ci-dessus les  
\begin{equation}\label{eq:XcalIdef}
X_{\calI} := (B, \xi_{\varnothing} +1_{\calI}, d_{\ol})
\end{equation}
où $\calI \sbse I_{0}$, sont des représentants des orbites 
pour l'action de $G(\rmF)$ à $\ol$. On les considère fixés désormais.

Pour un sous-$\rmF$-groupe $H$ de $G$, $X \in \tlgl(\rmF)$ 
et une $\rmF$-algèbre $R$, notons $H(R,X)$ le groupe des $R$-points 
du stabilisateur de $X$ dans $H$.
Alors, pour tout $\calI \sbse I_{0}$ et toute $\rmF$-algèbre $R$ on a:
\begin{equation}\label{eq:stabXcalIG}
G(R,X_{\calI}) = T_{I_{0} \smin \acalI}(R).
\end{equation}

Pour tout $i \in I$ on choisit une mesure de Haar sur $T_{i}(\A)$ et pour tout 
$I' \sbs I$ on met la mesure produit sur $T_{I'}(\A) = \prod_{i \in I'}T_{i}(\A)$. 
On note $H_{I}$ l'application de Harish-Chandra,  introduite 
dans le paragraphe \ref{par:prelimstraceG}, 
par rapport au sous-groupe parabolique standard $P_{I}$.
Pour $I' \sbs I$ on note $\all_{I'}$ l'image et $T_{I'}(\A)^{1}$ le noyau 
de la restriction de l'application $H_{I}$ à $T_{I'}(\A)$ 
(la restriction étant un homomorphisme car $T_{I'}(\A) \sbs M_{P_{I}}(\A)$).
En particulier $\all_{I} = \all_{P_{I}}$ et le quotient 
$T_{I'}(\rmF) \bsl T_{I'}(\A)^{1}$ est compact. 
Puisque $H_{I}$ restreint à $T_{I'}(\A) \cap A_{P_{I}}^{\infty}$ est un isomorphisme 
on obtient une unique mesure de Haar sur $T_{I'}(\A)^{1}$ compatible avec celles sur 
$T_{I'}(\A)$ et $\all_{I'}$.

Si $I'' \sbs I' \sbs I$, les espaces $\all_{I''}$ et $\all_{I' \smin I''}$ sont en somme directe. 
On voit alors $\all_{I'}^{*} := \Hom_{\R}(\all_{I'}, \R)$ comme un sous-espace 
de $\all_{0}^{*}$, et donc de $\all_{\tlzero}^{*}$ aussi.

En particulier, puisque $A_{P_{I}} \sbs T_{I} \sbs M_{P_{I}}$, on voit que 
$\all_{I}^{*}$ s'identifie à $\Hom_{\rmF}(T_{I}, \Gm) \otimes_{\Z} \R$. 
Pour tout $i \in I \cup - I$ 
on note alors $\upla_{i} \in \all_{I}^{*}$ le caractère par lequel $T_{I}$ agit sur $\rmF_{i}$. 
On a alors $\upla_{i} = -\upla_{-i}$. 
Soit $\calI \sbse I_{0}$. On a donc 
que $\{\upla_{i}\}_{i \in \calI}$ est une base de 
$\all_{\acalI}^{*}$. On pose aussi
\[
\upla_{\calI}  = \sum_{i \in \calI} \upla_{i} \in \all_{\acalI}^{*}.
\]
Notons que $\upla_{\calI^{\sharp}} = -\upla_{\calI}$.

On définit encore:
\begin{itemize}
\item pour $I' \sbs I$,  l'espace $\all_{I',\C}^{*} := \all_{I'}^{*} \otimes_{\R} \C$;
\item pour $\calI \sbse I_{0}$, 
$\ind_{\calI}$ la fonction caractéristique de 
$H \in \all_{I_{0}}$ tels que 
\[
\upla_{i}(H) \le 0, \  \forall \ i \in \calI \cap I_{0} \quad 
\text{et} \quad \upla_{i}(H) < 0, \  \forall \ i \in \calI \cap -I_{0}.
\] 
\end{itemize}

\subsection{Le résultat principal}\label{par:LeResultRssG}

Dans ce paragraphe on énonce et on démontre le théorème \ref{thm:theThmOrbsG}. 
Cependant, certains résultats seront seulement énoncés avec les renvois  
vers leurs démonstrations dans les paragraphes suivants. 

On note $M_{\tlI}$ le sous-groupe de Levi de $\tlG$ stabilisant 
$\rmF_{i}$ pour tout $i \in I$ ainsi que $D_{0}$ et l'on note 
$M_{\tlI_{0}}$ le plus grand sous-groupe de Levi de $\tlG$
qui stabilise les espaces $\rmF_{i}$ 
pour tout $i \in I_{0}$. On a donc $M_{\tlI} \sbs M_{\tlI_{0}}$.

Soit $\tlQ \in \calF(M_{\tlI})$, remarquons alors qu'il existe un unique
$\calI \sbse I$ tel que l'espace 
$\calV_{\tlQ}$ défini dans le paragraphe \ref{par:glnGlnplus1} 
égale $\rmF_{\calI}$. On pose dans ce cas $\calI_{\tlQ} = \calI$. 
D'ailleurs, il est clair que si l'orbite d'un $X_{\calI}$
(voir (\ref{eq:XcalIdef})) où $\calI \sbse I_{0}$ 
intersecte non-trivialement 
$\ml_{\tlP}(\rmF)$ où $\tlP \in \relPb$
alors celui-ci est conjugué 
à un élément de $\calF(M_{\tlI})$.
On peut alors reformuler 
le lemme \ref{shortLemmeOrbG} ainsi:
\blem\label{orbsAsSetsG}
 Soit $\tlQ \in \calF(M_{\tlI})$. Alors $X_{\calI} \in \ml_{\tlQ}(\rmF)$ 
si et seulement si $(\acalI \cup (I \smin I_{0})) \cap |\calI_{\tlQ}| = \varnothing$.
\elem

On a en particulier donc que si $\ml_{\tlQ}(\rmF) \cap \ol \neq \varnothing$, 
alors $(I \smin I_{0})  \cap |\calI_{\tlQ}| = \varnothing$ autrement dit 
$\tlQ \in \calF(M_{\tlI_{0}})$.

Pour un sous-groupe parabolique relativement standard $\tlP$ de $\tlG$ on pose
\[
\calF(M_{\tlI_{0}}, \tlP)^{rel} = \{\tlQ \in \calF(M_{\tlI_{0}})|\, 
\exists \gamma \in G(\rmF), \ \gamma \tlQ \gamma^{-1} = \tlP \}.
\]
Soit $\tlQ \in \calF(M_{\tlI_{0}}, \tlP)^{rel}$ et $\gamma \in G(\rmF)$ 
tel que $\gamma \tlQ \gamma^{-1} = \tlP$. Grâce à la décomposition de Bruhat, 
on choisi alors un $s_{\tlQ} \in \Omega^{G}/\Omega^{Q}$ tel 
que $\gamma \in B(\rmF)w_{s_{\tlQ}}Q(\rmF)$. On a alors $w_{s_{\tlQ}}\tlQ w_{s_{\tlQ}}^{-1} = \tlP$. 

On a alors:
\blem\label{lem:orbitCapPG} 
 Soient $\tlP \in \relPb$ et 
$\calI \sbse I_{0}$. 
Alors, l'intersection de la $G(\rmF)$-orbite de $X_{\calI}$ 
avec $\ml_{\tlP}(\rmF)$ égale
\[
\coprod_{
\begin{subarray}{c}
\tlQ \in \calF(M_{\tlI_{0}}, \tlP)^{rel} \\
\acalI \cap |\calI_{\tlQ}| = \varnothing
\end{subarray}
}
\coprod_{\delta \in M_{P}(\rmF, \Ad(w_{s_{\tlQ}})X_{\calI}) \bsl M_{P}(\rmF)}
\{\Ad(\delta^{-1}w_{s_{\tlQ}})X_{\calI}\}.
\]
\bdem
Voir la preuve du lemme 5.4 de \cite{leMoi}.
\edem
\elem

Pour une fonction $\upphi$ sur $\tlgl(\A)$ et 
$x \in G(\A)$ on définit
\[
\upphi_{x}(X) := \upphi(\Ad(x^{-1})X).
\]

Fixons $f \in \calS(\tlgl(\A))$.
Soit $\tlP$ un sous-groupe parabolique semi-standard de $\tlG$.
D'après le lemme \ref{shortLemmeOrbG} \textit{b)} 
on a un isomorphisme $N_{P}$-équivariant
\begin{equation}\label{eq:ntlPisnpcalVP}
\nl_{\tlP} \cong \nl_{P} \oplus \calV_{\tlP}
\end{equation}
On pose 
\begin{equation}\label{eq:fhatPG}
f^{\tlP}(X) = \int_{\calV_{\tlP}(\A)}f(X + Y)dY, \quad 
X \in \tlgl(\A).
\end{equation}
Supposons en plus que $\tlP$ est relativement standard et posons
\begin{equation}\label{eq:IPG}
i_{\tlP,\ol}(x) = i_{f,\tlP,\ol}(x) =
\sum_{\xi \in \ml_{\tlP}(\rmF) \cap \ol}
\sum_{\gamma \in N_{P}(\rmF)}f_{\gamma x}^{\tlP}(\xi), \quad 
x \in P(\rmF)\bsl G(\A),
\end{equation}
où $f_{x}^{\tlP} = (f_{x})^{\tlP}$.
En vertu du lemme \ref{lem:orbitCapPG} on a alors
\[
i_{\tlP,\ol}(x) = \sum_{\tlQ \in \calF(M_{\tlI_{0}}, \tlP)^{rel}}
\sum_{
\begin{subarray}{c}
\calI \sbs_{\epsilon} I_{0} \\
\acalI \cap |\calI_{\tlQ}| = \varnothing
\end{subarray}}
\sum_{\gamma \in P(\rmF, \Ad w_{s_{\tlQ}} X_{\calI}) \bsl P(\rmF)}
f_{\gamma x}^{\tlP}(\Ad (w_{s_{\tlQ}})X_{\calI}).
\]
On voit donc que
\begin{equation}\label{IPrewritten1G} 
\sum_{\delta \in P(\rmF) \bsl G(\rmF)}
i_{\tlP,\ol}(\delta x)
= 
\sum_{\tlQ \in \calF(M_{\tlI_{0}}, \tlP)^{rel}}
\sum_{
\begin{subarray}{c}
\calI \sbse I_{0} \\
\acalI \cap |\calI_{\tlQ}| = \varnothing
\end{subarray}} 
\sum_{\delta \in G(\rmF,X_{\calI}) \bsl G(\rmF)}
f_{\delta x}^{\tlQ}(X_{\calI}).
\end{equation}
On expliquera bientôt pourquoi la somme ci-dessus converge.

Pour tout $\calJ \sbse I_{0}$ posons $\A_{\calJ} = \A \otimes_{\rmF} \rmF_{\calJ}$. 
Notons aussi $\A_{I'} = \A \otimes_{\rmF} \rmF_{I'}$ 
et $\A_{-I'} = \A \otimes_{\rmF} \rmF_{-I'}$ pour tout $I' \sbs I$. 
Fixons un caractère additif continu non-trivial $\psi$ sur $\rmF \bsl \A $.
Pour une fonction $\upphi \in \calS(\tlgl(\A))$ et 
$\calJ \sbse I_{0}$ on définit la "transformée 
de Fourier" de $\upphi$ par rapport à $\calJ$
\begin{equation}\label{eq:trFourcalJ}
\hat \upphi^{\calJ}(X) = 
\int_{\A_{\calJ}}
\upphi(B', u_{\calJ} + \rmw' -\rmw'_{\calJ^{\sharp}}, d')
\psi (\langle u_{\calJ}, \rmw_{\calJ^{\sharp}}' \rangle )
du_{\calJ}
\end{equation}
où $X = (B', \rmw', d') \in \tlgl(\A)$, $du_{\calJ}$ c'est la mesure 
de Haar sur $\A_{\calJ}$
pour laquelle $\rmF_{\calJ} \bsl \A_{\calJ}$ est de volume $1$ et 
$\bilif$ c'est l'accouplement défini par la forme trace.
En fait, si l'on écrit 
$\rmw' = (u',v') \in \A_{I} \times \A_{-I}$ et 
$u_{\calJ} = (u, v)$ où $u \in \A_{\calJ \cap I}$ et
$v \in \A_{\calJ \cap -I}$ on a 
\[
\langle u_{\calJ}, \rmw_{\calJ^{\sharp}}' \rangle =
\langle u_{\calJ}, \rmw'  \rangle =
\Tr_{\rmF_{I} / \rmF}(u'\iota_{I}(v) + u\iota_{I}(v')).
\]
La restriction de la fonction $\hat \upphi^{\calJ}$ à 
$\gl(\A) \times \A_{I \smin (\calJ \cap I)} \times \A_{-I \smin (\calJ \cap -I)} \times \A$ est alors de type Bruhat-Schwartz. 
Dans tous ce qui suit les arguments des fonctions de type $\hat \upphi^{\calJ}$ seront toujours 
dans $\gl(\A) \times \A_{I \smin (\calJ \cap I)} \times \A_{-I \smin (\calJ \cap -I)} \times \A$ de façon qu'on pourrait les traiter 
comme les fonctions de type Bruhat-Schwartz.

Soient alors $\tlQ \in \calF(M_{\tlI_{0}}, \tlP)^{rel}$ et $\calI \sbse I_{0}$ tel que $\acalI \cap |\calI_{\tlQ}| = \varnothing$. 
On a $f_{x}^{\tlQ}(X_{\calI}) = \widehat{(f_{x})}^{\calI_{\tlQ}}(X_{\calI}) = \hf_{x}^{\calI_{\tlQ}}(X_{\calI})$. 
En utilisant (\ref{eq:stabXcalIG}), on voit donc qu'on peut réécrire (\ref{IPrewritten1G}) comme
\begin{equation}\label{eq:IpolrewrittenG}
\sum_{\tlQ \in \calF(M_{\tlI_{0}}, \tlP)^{rel}}
\sum_{
\begin{subarray}{c}
\calI \sbse I_{0} \\
\acalI \cap |\calI_{\tlQ}| = \varnothing
\end{subarray}} 
\sum_{\delta \in T_{I_{0} \smin \acalI}(\rmF) \bsl G(\rmF)}
\hf_{\delta x}^{\calI_{\tlQ}}(X_{\calI}).
\end{equation}
En particulier, la somme (\ref{IPrewritten1G}) est convergente et par conséquent 
bien définie.

On pose:
\[
i_{\ol}(x) = i_{f, \ol}(x) =
\sum_{\tlP \in \relPb}(-1)^{d_{\tlP}^{\tlG}}
\sum_{\delta \in P(\rmF) \bsl G(\rmF)}
i_{\tlP,\ol}(\delta x), \quad x \in G(\rmF) \bsl G(\A).
\]

\brop[cf. \ref{prop:newjolG}]\label{prop:newjol0G} On a 
pour tout $s \in \C$ tel que $-1 < \Rel(s) < 1$:
\[
\int_{G(\rmF) \bsl G(\A)}|i_{f,\ol}(x)\eta_{s}(\det x)|dx < \infty
\quad \text{et} \quad \int_{G(\rmF) \bsl G(\A)}i_{f,\ol}(x)\eta_{s}(\det x)dx = I_{\ol}(\eta_{s}, f).
\]
\erop
\nident
En utilisant le résultat ci-dessus,
on a $\int_{[G]}i_{\ol}(x)\eta_{s}(x)dx = I_{\ol}(\eta_{s}, f)$ pour $-1 < \Rel(s) < 1$,  où, grâce à
la formule (\ref{eq:IpolrewrittenG}), on a
\[
i_{\ol}(x) = \sum_{\tlP \in \relPb}(-1)^{d_{\tlP}^{\tlG}}
\sum_{\tlQ \in \calF(M_{\tlI_{0}}, \tlP)^{rel}}
\sum_{
\begin{subarray}{c}
\calI \sbse I_{0} \\
\acalI \cap |\calI_{\tlQ}| = \varnothing
\end{subarray}} 
\sum_{\delta \in T_{I_{0} \smin \acalI}(\rmF) \bsl G(\rmF)}
\hf_{\delta x}^{\calI_{\tlQ}}(X_{\calI}).
\]
En inversant l'ordre de sommation on a aussi
\[
i_{\ol}(x) =
\sum_{\calI \sqcup \calJ \sbse I_{0}}
\ \ \ \
\mu_{\calJ}
\sum_{\mathclap{\delta \in T_{I_{0} \smin \acalI}(\rmF) \bsl G(\rmF)}}
\ \ 
\hf^{\calJ}_{\delta x}(X_{\calI})
\]
où 
\[
\mu_{\calJ} = 
\sum_{\tlP \in \relPb}(-1)^{d_{\tlP}^{\tlG}}
\sum_{
\begin{subarray}{c}
\tlQ \in \calF(M_{\tlI_{0}}, \tlP)^{rel}\\
\calI_{\tlQ} = \calJ
\end{subarray}
}1.
\]
\blem Soit $\calJ \sbse I_{0}$. Alors $\mu_{\calJ} = (-1)^{\# \calJ}$.
\bdem
Remarquons d'abord qu'on a $\mu_{\calJ} = 
\sum_{
\begin{subarray}{c}
\tlQ \in \calF(M_{\tlI_{0}})\\
\calI_{\tlQ} = \calJ
\end{subarray}} (-1)^{d_{\tlQ}^{\tlG}}$. 
Soit $J_{1} = \calJ \cap I_{0}$ et $J_{2} = -(\calJ \cap -I_{0})$.
Soit $\tlQ \in \calF(M_{\tlI_{0}})$ et notons $k = d_{\tlQ}^{\tlG}$. 
Alors, la condition $\calI_{\tlQ} = \calJ$ veut dire
qu'il existe des sous-ensembles $\varnothing = I_{0} \sbn I_{1} \sbn \cdots \sbn I_{i} = J_{1}$ et 
$\varnothing = I_{0}' \sbn I_{1}' \sbn \cdots \sbn I_{k-i}' = J_{2}$ tels que 
$\tlQ$ est le stabilisateur du drapeau
\begin{multline*}
0 = \rmF_{I_{0}} \sbs \rmF_{I_{1}} \sbn \cdots \sbn \rmF_{I_{i}} = \rmF_{J_{1}}
\sbn \rmF_{I \smin J_{2}} \oplus D_{0} = 
\rmF_{(I \smin J_{2}) \cup I_{0}'} \oplus D_{0} \sbn \\
\rmF_{(I \smin J_{2}) \cup I_{1}'} \oplus D_{0} \sbn \cdots 
\sbn 
\rmF_{(I \smin J_{2}) \cup I_{k-i}'} \oplus D_{0} = \rmF_{I} \oplus D_{0}.
\end{multline*}

Si l'on pose alors $j_{1} = \# J_{1}$,  $j_{2} = \# J_{2}$ et
\[
a_{i}^{m} = \#\{(I_{0}, I_{1}, \ldots, I_{i})| 
\varnothing = I_{0} \sbn I_{1} \sbn \cdots \sbn 
I_{i} = \{1, 2, \ldots, m\}\}, \quad i,m  \in \N.
\]
On voit que le nombre des $\tlQ \in \calF(M_{\tlI_{0}})$ tels que 
$\calI_{\tlQ} = \calJ$ et  $k = d_{\tlQ}^{\tlG}$ égale
\[
c_{k}^{j_{1},j_{2}} := \sum_{i=0}^{k}a_{i}^{j_{1}}a_{k-i}^{j_{2}}.
\]
On veut alors montrer que l'expression
\[
\sum_{k=0}^{j_{1} + j_{2}}(-1)^{k}c_{k}^{j_{1},j_{2}} = 
(\sum_{i=0}^{j_{1}}(-1)^{i}a_{i}^{j_{1}})(\sum_{j=0}^{j_{2}}(-1)^{j}a_{j}^{j_{2}}) 
\]
égale $(-1)^{j_{1} + j_{2}}$.
Or, comme il est expliqué dans le lemme 5.6 de \cite{leMoi}, 
on a pour tout $m \in \N$ que $\sum_{i=0}^{m}(-1)^{i}a_{i}^{m} = (-1)^{m}$ donc, le résultat 
suit de l'identité ci-dessus.
\edem
\elem

On vient d'obtenir alors
\begin{equation}\label{eq:jIknewAllAlongG}
i_{\ol}(x) = 
\sum_{\calI \sqcup \calJ \sbse I_{0}}
(-1)^{\# \calJ}
\sum_{\mathclap{\delta \in T_{I_{0} \smin \acalI}(\rmF) \bsl G(\rmF)}}
\ \ 
\hf^{\calJ}_{\delta x}(X_{\calI}).
\end{equation}

\blem[cf. corollaire \ref{cor:UpsPros}]\label{lem:UpsPros0G}
 Soient $\calJ, \calJ_{1}, \calJ_{2} \sbse I_{0}$ 
tels que $\calJ_{1} \sqcup \calJ_{2} \sbs \calJ$.
L'intégrale suivante
\begin{equation*}
\int\limits_{\mathclap{
T_{I_{0} \smin |\calJ_{12}|}(\rmF)
\bsl G(\A)}}
\ind_{(\calJ \smin \calJ_{2}) \cup \calJ_{2}^{\sharp}}(H_{I}(x))
\hf_{x}^{\calJ \smin \calJ_{1}}(X_{\calJ_{1} \cup \calJ_{2}^{\sharp}})\eta_{s}(\det x)dx, \quad 
s \in \C
\end{equation*}
converge absolument et uniformément sur tous les compacts de $-1 < \Rel(s) < 1$ et 
admet un prolongement méromorphe à $\C$, noté $\brLa_{\calJ_{1},\calJ_{2}}^{\calJ}(\eta, f)(s)$, 
holomorphe sur $s \in \C \smin \{-1,1\}$.
\elem

Le rapport entre les $\brLa_{\calJ_{1},\calJ_{2}}^{\calJ}(\eta, f)(s)$ et $I_{\ol}(\eta_{s}, f)$ est donné 
par le lemme suivant.

\blem\label{lem:distIsUpsilonG} Pour tout  $s \in \C \smin \{-1,1\}$ on a:
\[
I_{\ol}(\eta_{s}, f) = \sum_{\acalJ = I_{0}}
\sum_{\calJ_{1} \sqcup \calJ_{2} \sbs \calJ}
(-1)^{\# (\calJ \smin \calJ_{12})}
\brLa_{\calJ_{1},\calJ_{2}}^{\calJ}(\eta, f)(s).
\]
\bdem 
En raisonnant comme dans la preuve du lemme 5.9 de \cite{leMoi}, on obtient 
pour tout $x \in G(\A)$:
\begin{equation*}
\sum_{\acalJ = I_{0}}
\sum_{\calJ_{1} \sqcup \calJ_{2} \sbs \calJ}
(-1)^{\# (\calJ \smin \calJ_{12})}
\sum_{\mathclap{\delta \in T_{|\calJ_{12}|}(\rmF)}}
\hf_{\delta x}^{\calJ \smin \calJ_{1}}
(X_{\calJ_{1} \cup \calJ_{2}^{\sharp}}) 
\ind_{(\calJ \smin \calJ_{2}) \cup \calJ_{2}^{\sharp}}(H_{I}(x))
=
\sum_{\mathclap{\calI \sqcup \calJ \sbse I_{0}}}
(-1)^{\# \calJ}
\sum_{\mathclap{\delta \in T_{|\calI|}(\rmF)}}
\hf_{\delta x}^{\calJ}
(X_{\calI}).
\end{equation*}
On multiplie cette égalité par $\eta_{s}(\det x)$, où $s \in \C$ est tel que $-1 < \Rel(s) < 1$.
En regardant le côté droit de cette égalité et
en utilisant la formule (\ref{eq:jIknewAllAlongG}), 
on voit que, en vertu de la proposition \ref{prop:newjol0G}, 
l'intégrale de cette expression sur 
$T_{I_{0}}(\rmF) \bsl G(\A)$ égale $I_{\ol}(\eta_{s}, f)$.
Or, en vertu du corollaire
\ref{lem:UpsPros0G}
l'intégrale 
du côté gauche sur le même quotient 
donne l'égalité cherchée pour $s \in \C$ tels que $-1 < \Rel(s) < 1$.
Ceci suffit pour conclure car $s \mapsto I_{\ol}(\eta_{s}, f)$ est une fonction 
méromorphe en vertu de la remarque \ref{rem:Iolhomo}.
\edem
\elem

On introduit maintenant les fonctions zêta.
\brop[cf. \ref{prop:zetaDefPropsG} et lemme \ref{lem:detIsNotSing}]
Soit $\calJ \sbse I_{0}$. Alors l'intégrale
\begin{equation*}
\int\limits_{\mathclap{
G(\A, X_{\calJ})\bsl G(\A)}}
f(x^{-1}X_{\calJ}x)
e^{\la(H_{I}(x))}
\eta(\det x)
dx, \quad \la \in \all_{\acalJ,\C}^{*}
\end{equation*}
converge absolument 
et uniformément sur tous les compacts d'un 
ouvert non-vide de $\all_{\acalJ,\C}^{*}$ 
et admet un prolongement méromorphe 
à $\all_{\acalJ,\C}^{*}$, noté $\zeta_{\calJ}(\eta, f)$. En plus, 
la droite $\C \det \sbs \all_{\acalJ,\C}^{*}$ est non-singulière.
\erop

\bcor\label{cor:ZetaPros0}
Pour $\calJ \sbse I_{0}$
on définit la fonction 
méromorphe sur $\C$ suivante:
\[
\zeta_{\calJ}(\eta, f)(s) := \zeta_{\calJ}(\eta, f)(s \det), 
\quad s \in \C
\]
où l'on voit $\det$ comme un élément de $\all_{\acalJ}^{*}$ par restriction.
\ecor

\brem 
Il découle de la proposition \ref{prop:zetaDefPropsG}, que
pour certains $\calJ \sbse I_{0}$, pour tout $s \in \C$
l'élément $s \det \in \all_{\acalJ, \C}^{*}$ n'appartient pas 
au domaine de convergence de l'intégrale définissant 
$\zeta_{\calJ}(\eta, f)$.
\erem

Soit $I_{\eta} \sbs I_{0}$ l'ensemble de $i \in I_{0}$ tels que la restriction de $\eta$ au groupe 
de normes de $\A_{i}^{*}$ dans $\A^{*}$ est non-triviale.
On est prêt à démontrer le résultat principal de cette section.

\begin{theo}\label{thm:theThmOrbsG} 
Soit $f \in \calS(\tlgl(\A))$. 
\begin{enumerate}[i)]
\item Pour tout $\calJ_{1} \sbse I_{\eta}$ tel que $|\calJ_{1}| = I_{\eta}$ la fonction 
de la variable $s \in \C$ suivante
\[
\sum_{|\calJ'| = I_{0} \smin I_{\eta}} \zeta_{\calJ_{1} \cup \calJ'}(\eta, f)
\]
est holomorphe sur  $\C \smin \{-1,1\}$.
\item
On a pour tout $s \in \C \smin \{-1, 1\}$:
\[
I_{\ol}(\eta_{s}, f) = \sum_{|\calJ_{1}| = I_{\eta}} \dsl (
\sum_{|\calJ'| = I_{0} \smin I_{\eta}} \zeta_{\calJ_{\eta} \cup \calJ'}(\eta, f) 
)(s) \rb.
\]
\end{enumerate}
\bdem 
Le point \textit{i)} est démontré dans le corollaire \ref{cor:zetasCsum} ci-dessous.
En outre, en vertu du lemme 5.27 de \cite{leMoi}, on a l'égalité suivante des fonctions méromorphes sur $\C$: 
\[ 
\sum_{\acalJ = I_{0}}
\sum_{\calJ_{1} \sqcup \calJ_{2} \sbs \calJ}
(-1)^{\#(\calJ \smin \calJ_{12})}
\brLa_{\calJ_{1}, \calJ_{2}}^{\calJ}(\eta, f) = 
\sum_{|\calJ| = I_{0}} \zeta_{\calJ}(\eta, f).
\]
En utilisant alors l'égalité démontrée dans le 
lemme \ref{lem:distIsUpsilonG} ainsi que le point \textit{i)} de ce théorème
on peut conclure. 
\edem
\end{theo}

%% file: koSansT_gln.tex
\subsection{Deuxième formule pour le noyau tronqué}\label{par:noyuTronqNouvG}

On considère $f \in \calS(\tlgl(\A))$ fixée. Pour $\tlP \in \relPb$, 
dans le paragraphe \ref{par:LeResultRssG}, équation (\ref{eq:fhatPG}), 
on a introduit la fonction $f^{\tlP}$ ainsi que $f_{x}(X) := f(\Ad(x^{-1})X)$. 
Par $f_{x}^{\tlP}$ on entend toujours $(f_{x})^{\tlP}$. 
Soient $\tlP, \tlQ \in \relPb$ tels que $\tlP \sbs \tlQ$.  
On pose alors
\begin{equation}\label{eq:calVPQ}
\calV_{\tlP}^{\tlQ} := \calV_{\tlP} \cap \calZ_{\tlQ}.
\end{equation}
Alors, pour $\tlQ$ fixé, les espaces $\calV_{\tlP}^{\tlQ}$, où $\relPb \ni \tlP \sbs \tlQ$, 
jouent le rôle des espaces $\calV_{\tlP}$ 
dans le contexte de l'inclusion $G_{\tlQ} \hrar \tlG_{\tlQ}$. 
En plus,  on a les relations suivantes:
\begin{equation}\label{eq:notForMe}
\nl_{Q} \oplus \calV_{\tlP} = \nl_{\tlQ} \oplus \calV_{\tlP}^{\tlQ}, \quad
\nl_{\tlP}^{\tlQ} \oplus \nl_{Q} = \nl_{P} \oplus \calV_{\tlP}^{\tlQ}.
\end{equation}

Les lemmes suivants sont des analogues 
des corollaires 5.13 et 5.14 de \cite{leMoi}.
\blem\label{lem:lemcor1G}
 Soient $\tlP \in \relPb$,
 $\xi \in  \ml_{\tlP}(\rmF) \cap \mathfrak{o}$ et $x \in G(\A)$, alors
\begin{displaymath}
\sum_{\gamma \in N_{P}(\rmF)} f_{\gamma x}^{\tlP}(\xi)= 
\sum_{\zeta \in \nl_{P}(\rmF)}
f_{x}^{\tlP}(\xi +\zeta).
\end{displaymath}
\elem

\blem\label{lem:lemcor2G}
 Soient $\tlP \in \relPb$,
 $\xi \in  \ml_{\tlP}(\rmF) \cap \mathfrak{o}$ et $x \in G(\A)$, alors
\begin{displaymath}
\int_{N_{P}(\A)} 
f_{nx}^{\tlP}(\xi)dn = \int_{\nl_{\tlP}(\A)}
f_{x}(\xi+U)dU.
\end{displaymath}
\elem

Pour $T \in \all_{\tlzero}^{+}$ posons
\begin{displaymath}
i_{\mathfrak{o}}^{T}(x) = i_{f,\mathfrak{o}}^{T}(x) = 
\sum_{\tlP \in \relPb}(-1)^{d_{\tlP}^{\tlG}}
\sum_{\delta \in P(\rmF)\backslash G(\rmF)}
\htau_{\tlP}(H_{\tlP}(\delta x)-T_{\tlP})i_{\tlP,\mathfrak{o}}(\delta x), \ 
x \in G(\rmF)\backslash G(\A)
\end{displaymath}
où $i_{\tlP, \ol}$ est définie dans la section \ref{par:LeResultRssG} 
par (\ref{eq:IPG}).
La fonction $i_{\mathfrak{o}}^{T}$ 
est une variante de $k_{\ol}^{T}$ définie 
au début du paragraphe \ref{par:convergenceG}. 

\btheo\label{thm:convThmNotMainG} 
Pour tous $T \in T_{+} + \all_{\tlzero}^{+}$, 
$f \in \mathcal{S}(\tlgl(\A))$, $\sigma \in \R$ et 
$\mathfrak{o} \in \mathcal{O}_{rs}$, on a:
\[
\int_{G(\rmF)\backslash G(\A)}|i_{f,\mathfrak{o}}^{T}(x)| |\det x|_{\A}^{\sigma} dx < \infty.
\]
\bdem
 En procédant comme au début 
de la preuve du théorème \ref{thm:MainConvG} on montre 
que l'intégrale $\int_{[G]}|i_{f,\mathfrak{o}}^{T}(x)||\det x|_{\A}^{\sigma}dx$
est majorée par la somme sur les sous-groupes 
paraboliques relativement standards 
$\tlP_{1} \subseteq \tlS \subseteq \tlP_{2}$ 
de 
\[
\int_{P_{1}(\rmF) \bsl G(\A)}
\sum_{
\xi \in (\ml_{\tlone}^{\tlS})'(\rmF) \cap \mathfrak{o}
}
\chi^{T}_{\tlone,\tltwo}(x)
|\sum_{\tlS \subseteq \tlP \subseteq \tlP_{2}}
(-1)^{d_{\tlP}^{\tlG}} 
\sum_{\zeta \in \nl_{\tlS}^{\tlP}(\rmF)}
\sum_{\gamma \in N_{P}(\rmF)}f_{\gamma x}^{\tlP}(\xi+\zeta)||\det x|_{\A}^{\sigma}dx.
\]

En utilisant le lemme \ref{lem:lemcor1G} ci-dessus (pour 
la fonction $X \mapsto f(X + V)$ o\`u 
$V \in \nl_{\tlS}^{\tlP}(\rmF)$)
et ensuite la formule sommatoire de 
Poisson on s'aperçoit que la somme entre la valeur absolue 
dans l'intégrale ci-dessus égale
\begin{equation}\label{eq:newKernelConv}
|\sum_{\tlS \subseteq \tlP \subseteq \tlP_{2}}
(-1)^{d_{\tlP}^{\tlG}} 
\sum_{\zeta_{1} \in \bar \nl_{\tlS}^{\tlP}(\rmF)}
\sum_{\zeta_{2} \in \bar \nl_{P}(\rmF)}
\widetilde{\phi}_{\tlS}(x,\xi,\zeta_{1}+\zeta_{2}) |
\end{equation}
où:
\begin{displaymath}
\widetilde{\phi}_{\tlS}(x,X,Y) = 
\int_{\nl_{\tlS}(\A)}f(\Ad(x^{-1})(X+U))\Psi(\langle U,Y\rangle)dU, 
\ x \in G(\A), X \in \ml_{\tlP_{1}}^{\tlS}(\A), 
Y \in \bar \nl_{\tlS}(\A)
\end{displaymath}
et $\bar \nl_{P}$ est le radical nilpotent de l'algèbre 
de Lie du sous-groupe 
parabolique opposé à $P$.

On pose 
$\fl_{\tlS}^{\tlP} = \brnl_{\tlS}^{\tlP} \oplus \brnl_{P}$. Grâce aux décompositions (\ref{eq:notForMe})
on voit qu'on a $\fl_{\tlS}^{\tlR} \sbs \fl_{\tlS}^{\tlP}$ pour tout 
$\tlS \sbs \tlR \sbs \tlP$.
On pose donc
\[
(\fl_{\tlS}^{\tlP})' = \fl_{\tlS}^{\tlP} \smin \bigcup_{\tlS \sbs \tlR \sbn \tlP}\fl_{\tlS}^{\tlR}.
\]
On a donc $\fl_{\tlS}^{\tlP} = \coprod_{\tlS \sbs \tlR \sbs \tlP}(\fl_{\tlS}^{\tlR})' $ et
par un raisonnement habituel basé sur l'identité (\ref{basicidentityG}) 
on s'aperçoit que (\ref{eq:newKernelConv}) égale
\[
|
\sum_{\zeta \in (\fl_{\tlS}^{\tltwo})'(\rmF)}
\widetilde{\phi}_{\tlS}(x,\xi,\zeta) |.
\]

Posons alors
\begin{displaymath}
\Psi_{\tlS}(x,X,Y) = 
\sum_{\zeta_{2} \in \bar \nl_{2}(\rmF)}
\widetilde{\phi}_{\tlS}(x,X,Y+\zeta_{2}), \quad 
x \in G(\A), X \in \ml_{\tlone}^{\tlS}(\A), 
Y \in \bar \nl_{\tlS}^{\tltwo}(\A).
\end{displaymath}
Alors $\Psi_{S}(x,X,Y) \in 
\mathcal{S}((\ml_{\tlone}^{\tlS} \oplus \bar \nl_{\tlS}^{\tltwo})(\A))$
pour un $x$ fixé et on a:
\[
\sum_{\xi \in (\ml_{\tlone}^{\tlS})'(\rmF) \cap \mathfrak{o}}
|
\sum_{\zeta \in (\fl_{\tlS}^{\tltwo})'(\rmF)}
\widetilde{\phi}_{\tlS}(x,\xi,\zeta) 
| \le 
\sum_{\xi \in (\ml_{\tlone}^{\tlS})'(\rmF) \cap \mathfrak{o}}
\sum_{\zeta_{1} \in (\bar \nl_{\tlS}^{\tltwo})'(\rmF)}
|\Psi_{S}(x,\xi,\zeta_{1})|
\]
car l'image de $(\fl_{\tlS}^{\tltwo})'$ par la projection naturelle sur $\bar \nl_{\tlS}^{\tltwo}$ 
est contenue dans $(\bar \nl_{\tlS}^{\tltwo})'$. 

On se ramène alors à borner, pour $\tlP_{1} \sbs \tlS \sbs \tlP_{2}$ fixés:
\[
\int_{P_{1}(\rmF)\backslash G(\A)}
\chi_{\tlone,\tltwo}^{T}(x)
\sum_{\xi \in (\ml_{\tlone}^{\tlS})'(\rmF) \cap \mathfrak{o}}
\sum_{\zeta_{1} \in (\bar \nl_{\tlS}^{\tltwo})'(\rmF)}
|\Phi_{\tlS}(x,\xi,\zeta_{1})||\det x|_{\A}^{\sigma}dx.
\]
Cette intégrale est identique à (\ref{mainThmConv3G}) du théorème 
\ref{thm:MainConvG}. Cela conclut la preuve.
\edem
\etheo

\brop\label{prop:smalljIntG}
Pour $T \in T_{+} + \all_{\tlzero}^{+}$, $s \in \C$,
$f \in \mathcal{S}(\tlgl(\A))$ et $\mathfrak{o} \in \mathcal{O}_{rs}$, on a:
\begin{displaymath}
I^{T}_{\mathfrak{o}}(\eta_{s}, f) = \int_{G(\rmF)\backslash G(\A)}
i_{f, \mathfrak{o}}^{T}(x)\eta_{s}(\det x)dx.
\end{displaymath}
\bdem 
Dans la preuve on utilisera la notation 
du paragraphe
\ref{par:convergenceG} introduite au début de la preuve du 
théorème \ref{thm:MainConvG}. Donc, en raisonnant comme au début 
de la preuve de ce théorème-là, on voit que
\[
\int_{[G]}i_{f, \mathfrak{o}}^{T}(x) \eta_{s}(\det x)dx = 
\sum_{\tlP_{1} \in \relPb}
\sum_{\tlP_{2} \sps \tlP_{1}}
\int_{P_{1}(\rmF) \bsl G(\A)}
\chi_{\tlone,\tltwo}^{T}(x)
\sum_{\tlP_{1} \sbs \tlP \sbs \tlP_{2}}
(-1)^{d_{\tlP}^{\tlG}}
i_{\tlP,\ol}(x)
\eta_{s}(\det x)
dx.
\]
Fixons les sous-groupes paraboliques relativement standards $\tlP_{1} \sbs \tlP_{2}$ et 
décomposons l'intégrale sur 
$P_{1}(\rmF)\bsl G(\A)$ 
en une double intégrale 
sur $x \in M_{1}(\rmF)N_{1}(\A)\bsl G(\A)$ 
et $n_{1} \in N_{1}(\rmF) \bsl N_{1}(\A)$. 
Ensuite on fait passer cette dernière intégrale à l'intérieure 
de la somme sur $\tlP$. On peut le faire car pour tout 
$\tlP_{1} \sbs \tlP \sbs \tlP_{2}$ la fonction $N_{1}(\A) \ni n_{1} \mapsto i_{\ol, \tlP}(n_{1}x)$ 
est $N_{1}(\rmF)$-invariante et continue. 
Comme le volume de $N_{P}(\rmF) \bsl N_{P}(\A)$ vaut $1$ 
et $N_{P} \sbs N_{1}$ on a 
\begin{equation*}
\int\limits_{[N_{1}]}i_{\tlP,\ol}(n_{1}x)dn_{1} \mathrlap{=}
\int\limits_{[N_{1}]}
\sum_{\xi \in \ml_{\tlP}(\rmF) \cap \ol}
\int\limits_{[N_{P}]}
\sum_{\mathrlap{\gamma \in N_{P}(\rmF)}}
f_{\gamma n n_{1} x}^{\tlP}( \xi)dndn_{1} \mathrlap{=} 
\int\limits_{[N_{1}]}
\sum_{\xi \in \ml_{\tlP}(\rmF) \cap \ol}
\int\limits_{\mathrlap{N_{P}(\A)}}
f_{n n_{1} x}^{\tlP}( \xi)dndn_{1}.
\end{equation*}
La dernière expression égale $\int_{[N_{1}]}k_{\tlP,\ol}(n_{1}x)dn_{1}$ en vertu du lemme
\ref{lem:lemcor2G}. On intervertit de nouveau 
la somme qui porte sur $\tlP$ avec l'intégrale 
sur $N_{1}(\rmF) \bsl N_{1}(\A)$ et l'on 
recombine cette dernière avec l'intégrale sur 
$M_{1}(\rmF)N_{1}(\A)\bsl G(\A)$. On retrouve donc
\[
\int_{G(\rmF) \bsl G(\A)}i_{f, \mathfrak{o}}^{T}(x)\eta_{s}(\det x)dx = 
\sum_{\tlP_{1} \in \relPb}
\sum_{\tlP_{2} \sps \tlP_{1}}
\int_{P_{1}(\rmF) \bsl G(\A)}
\chi_{\tlone,\tltwo}^{T}(x)
k_{\tlone,\tltwo,\ol}(x)
\eta_{s}(\det x)
dx.
\]
Chaque intégrale dans la somme ci-dessus converge d'après le théorème 
\ref{thm:mainConv2} ce qui justifie l'intégration et de surcroît, 
la somme elle-même
égale $I_{\ol}^{T}(\eta_{s}, f)$ 
grâce à l'identité (\ref{eq:koAlternating}) et 
la définition 
de $I_{\ol}^{T}(\eta_{s},f)$ donnée 
au début de la section \ref{quantSect},
 ce qu'il fallait démontrer. 
 \edem 
 \erop

 Plaçons-nous maintenant dans le cadre du paragraphe \ref{par:JTforLevisG}. 
 On veut généraliser le théorème \ref{thm:convThmNotMainG} et la proposition 
 \ref{prop:smalljIntG} au cas de l inclusion $H \times G' \hrar H \times \tlG'$. 
 
Notons $\calO^{H \times \tlG'}_{rs}$ l'ensemble de classes 
contenant un élément $X_{1} + X_{2} \in \hl(\rmF) \times \tlgl'(\rmF)$ 
tel que le polynôme caractéristique de $X_{1}$ est séparable 
et tel que $X_{2} \in \tlgl'(\rmF)$ appartient à une classe relativement semi-simple régulière
 dans le contexte d'inclusion $G' \hrar \tlG'$. 
 
Pour $f \in \calS((\hl \times \tlgl')(\A))$, $\ol \in \calO^{H \times \tlG'}_{rs}$ 
et $\tlP \in \calF_{H \times G'}(M_{\tlzero}, B)$ 
soit
 \[
i_{f,\tlP,\ol}(x) = 
\sum_{\xi \in \ml_{\tlP}(\rmF) \cap \ol}
\sum_{\gamma \in N_{P}(\rmF)}\int_{\calV_{\tlP}'(\A)}
f(\Ad((\gamma x)^{-1})(\xi + Y_{\tlP}'))dY_{\tlP}'
 \]
 où $x \in P(\rmF) \bsl (H(\A)^{1} \times \tlG'(\A))$
 et $\calV_{\tlP}'$ c'est le plus grand sous-espace de $V' \times (V')^{*}$ 
 stabilisé par $\tlP$.
 
 Pour $T \in \all_{\tlBmin}^{+}$ on pose aussi
 \[
i_{f,\ol}^{T}(x) = 
\sum_{\tlP \in \calF_{H \times G'}(M_{\tlzero}, B)}(-1)^{d_{\tlP}^{H \times \tlG'}}
\sum_{\mathclap{\delta \in P(\rmF) \bsl (H \times G')(\rmF)}}
\
\htau_{\tlP}^{H \times \tlG'}(H_{\tlP}(\delta x) - T_{\tlP})
i_{f, \tlP,\ol}(\delta x), \quad 
x \in (H \times G')(\rmF) \bsl (H(\A)^{1} \times G'(\A)).
 \] 
 La preuve du théorème \ref{thm:convThmNotMainG} s'étend sans problème 
 dans ce cas donnant pour tout $\sigma \in  \R$
 \[
\int_{ (H \times G')(\rmF) \bsl (H(\A)^{1} \times G'(\A))}
|i_{f,\ol}^{T}(x)| |\det x|_{\A}^{\sigma}dx < \infty 
\]
pour $T$ suffisamment régulier. De même la preuve 
de la proposition \ref{prop:smalljIntG} 
s'étend aussi bien et l'on obtient, 
avec la notation du paragraphe \ref{par:JTforLevisG}, 
pour $\ol\in \calO^{H \times \tlG'}_{rs}$ et $s \in \C$
\[
\int_{ (H \times G')(\rmF) \bsl (H(\A)^{1} \times G'(\A))}
i_{f,\ol}^{T}(x)\eta_{s}(\det x)dx = 
I_{\ol}^{H \times \tlG',T}(\eta_{s}, f).
\]

Soient maintenant $\tlQ \in \relPb$, 
$\ol \in \calO_{rs}$ et 
$\ol_{\tlQ,1}, \ldots, \ol_{\tlQ,m} \in \calO^{M_{\tlQ}}$
comme dans l'équation (\ref{eq:mltlQcapolG}). Alors $\ol_{\tlQ,i} \in \calO^{M_{\tlQ}}_{rs}$
pour $i = 1, \ldots, m$. 
On a dans ce cas l'analogue suivant de 
l'égalité (\ref{eq:JMQisThisG}): pour
tout sous-groupe de Borel 
relativement standard $\tlB \sbs \tlQ$ et
tous $f \in \calS(\tlgl(\A))$ et 
$s \in \C$
\begin{multline}\label{eq:JMQisThisRsG}
I^{M_{\tlQ},T}_{\mathfrak{o}}(\eta_{s}, f_{\tlQ}) = 
\int\limits_{\mathclap{
M_{Q}(\rmF) \backslash H_{Q}(\A)^{1} \times G_{\tlQ}(\A)}} \qquad \qquad
\sum_{i=1}^{m}i_{f_{\tlQ},\ol_{\tlQ,i}}^{T_{\tlB}}(m)\eta_{s_{\tlQ}}(\det m)dm  = 
\int\limits_{\mathclap{M_{Q}(\rmF) \backslash H_{Q}(\A)^{1} \times G_{\tlQ}(\A)}}
e^{-\upla_{\tlQ,s}(H_{\tlQ}(m))}
\sum_{\tlP \sbs \tlQ}(-1)^{d_{\tlP}^{\tlQ}} \\
\sum_{\mathclap{\delta \in (P \cap M_{Q})(\rmF) \bsl M_{Q}(\rmF)}} \ \
\htau_{P}^{Q}(H_{P}(\delta m)-T)
\dsl
\sum_{
\begin{subarray}{c}
\xi \in \ml_{\tlP}(\rmF) \cap \ol \\
\gamma \in N_{P}^{Q}(\rmF)
\end{subarray}}
\int\limits_{\mathrlap{\calV_{\tlP}^{\tlQ}(\A)}}
f_{\tlQ}(\Ad((\gamma \delta m)^{-1})(\xi + Y))dY
\rb \eta_{s}(\det m)dm.
\end{multline}


\subsection{Expression intégrale de \texorpdfstring{$J_{\ol}$}{Jo}}\label{par:IntReprDeJG}

Dans ce paragraphe, on démontre 
la proposition suivante, énoncée dans le paragraphe \ref{par:LeResultRssG} 
par la proposition \ref{prop:newjol0G}.

\brop\label{prop:newjolG} Pour tout $s \in \C$  tel que $-1 < \Rel(s) < 1$ on a
\[
\int_{G(\rmF) \bsl G(\A)}|i_{f,\ol}(x)\eta_{s}(\det x)|dx < \infty
\quad \text{et} \quad \int_{G(\rmF) \bsl G(\A)}i_{f,\ol}(x)\eta_{s}(\det x)dx = I_{\ol}(\eta_{s},f).
\]

\bdem 
Pour tout $\tlQ \in \relPb$ 
soit $\brtau_{\tlQ}$ la fonction caractéristique de 
$H \in \all_{\tlQ}$ tels que $\al(H) \le 0$ 
pour tout $\al \in \Delta_{\tlQ}$. 
Il résulte du lemme combinatoire de Langlands 
(Proposition 1.7.2 de \cite{labWal}) que pour tout 
$\tlP \in \relPb$ et tout $H \in \all_{\tlP}$ on a
\[
\sum_{\tlQ \sps \tlP}\htau_{\tlP}^{\tlQ}(H)\brtau_{\tlQ}(H) = 1.
\]
En utilisant cette identité, on a 
pour tout $T \in T_{+} + \all_{\tlzero}^{+}$
\begin{multline*}
i_{f,\ol}(x) =
\sum_{\tlP \in \relPb}(-1)^{d_{\tlP}^{\tlG}}
\sum_{\gamma \in P(\rmF) \bsl G(\rmF)}
i_{\tlP,\ol}(\gamma x) = 
\sum_{\tlP \in \relPb}(-1)^{d_{\tlP}^{\tlQ}} \\
\sum_{\mathclap{\gamma \in P(\rmF) \bsl G(\rmF)}} \ \
i_{\tlP,\ol}(\gamma x)
(\sum_{\tlQ \sps \tlP}\htau_{\tlP}^{\tlQ}(H_{\tlP}(\gamma x) - T_{\tlP})
\brtau_{\tlQ}(H_{\tlQ}(\gamma x) - T_{\tlQ})) = \ 
\sum_{\mathclap{\tlQ \in \relPb}} \ (-1)^{d_{\tlQ}^{\tlG}} \\
\sum_{\mathclap{\gamma \in Q(\rmF) \bsl G(\rmF)}} \ \
\brtau_{Q}(H_{Q}(\gamma x) - T_{Q}) \ 
\sum_{\mathclap{\relPb \ni \tlP \sbs \tlQ}} \ \quad
\quad (-1)^{d_{\tlP}^{\tlQ}} \
\sum_{\mathclap{\delta \in (M_{Q} \cap P)(\rmF) \bsl M_{Q}(\rmF)}} \quad
i_{\tlP,\ol}(\eta \gamma x)\htau_{\tlP}^{\tlQ}(H_{\tlP}(\delta \gamma x) - T_{\tlP})
\end{multline*}
où les sommes sont absolument convergentes.
Il suffit de montrer que pour tout $\tlQ \in \relPb$ 
et tout $s \in \C$ tel que $-1 < \Rel(s) < 1$
l'intégrale
\begin{equation}\label{eq:suffitPourQG}
\int_{Q(\rmF) \bsl G(\A)}
\brtau_{\tlQ}(H_{\tlQ}(x) - T_{\tlQ})
\sum_{\relPb \ni \tlP \sbs \tlQ}
(-1)^{d_{\tlQ}^{\tlP}} \ \ 
\sum_{\mathclap{\delta \in (M_{Q} \cap P)(\rmF) \bsl M_{Q}(\rmF)}} \
i_{\tlP,\ol}(\delta x)\htau_{\tlP}^{\tlQ}(H_{\tlP}(\delta x) - T_{P})
\eta_{s}(\det x)
dx
\end{equation}
converge absolument. 
L'analyse va être analogue à celle de la preuve du théorème 
\ref{thm:mainQualitThm}.

Posons
$x = namk$ o\`u $n \in N_{Q}(\rmF)\backslash N_{Q}(\A)$, 
$m \in M_{Q}(\rmF)\backslash H_{\tlQ}(\A)^{1} \times G_{\tlQ}(\A)$, 
$a \in A_{\tlQ}^{st, \infty}$ et $k \in K$. Donc 
$dx = e^{-2\rho_{Q}(H_{Q}(am))}dndadmdk$. 

Fixons $\tlP \in \relPb$ tel que $\tlP \sbs \tlQ$.
Pour $n$, $m$, $k$ et $a$ comme ci-dessus et $\delta \in M_{Q}(\rmF)$ on a
$\brtau_{\tlQ}(H_{\tlQ}(\delta namk) - T_{\tlQ}) = \brtau_{\tlQ}(H_{\tlQ}(ma) - T_{\tlQ})$
et, en faisant les changements de variable 
$\delta a^{-1}na \delta^{-1} \mapsto n$ et $ a^{-1}Y_{\tlP} \mapsto Y_{\tlP}$, 
\begin{multline}\label{eq:passAuLeviG}
\int\limits_{\mathrlap{[N_{Q}]}}
i_{\tlP,\ol}(\delta n amk)\htau_{\tlP}^{\tlQ}(H_{\tlP}(\delta namk) - T_{\tlP})
dn = 
\int\limits_{\mathclap{K}}\int\limits_{\mathrlap{[N_{Q}]}}
i_{\tlP,\ol}(\delta n amk)\htau_{\tlP}^{\tlQ}(H_{\tlP}(\delta m) - T_{\tlP})
dn = \\
\htau_{\tlP}^{\tlQ}(H_{\tlP}(\delta m) - T_{\tlP}) e^{2\rho_{\tlQ}(H_{Q}(a))}
\sum_{\xi \in \ml_{\tlP}(\rmF) \cap \ol}
\sum_{\gamma \in N_{P}^{Q}(\rmF)}
\int\limits_{N_{Q}(\A)}
\int\limits_{\mathrlap{\calV_{\tlP}(\A)}}
f_{n \gamma \delta mk}(\xi + Y)dYdn
\end{multline}
où l'on utilise le fait que la somme de poids pour $A_{H_{\tlQ}} \sbs A_{\tlQ}$ 
agissant sur $\nl_{Q} \oplus \calV_{\tlP} = \nl_{\tlQ} \oplus \calV_{\tlP}^{\tlQ}$ 
(voir paragraphe \ref{par:noyuTronqNouvG}, equations (\ref{eq:calVPQ}) et (\ref{eq:notForMe}))
égale $2\rho_{\tlQ}$ car $A_{H_{\tlQ}}$ agit trivialement sur $\calZ_{\tlQ} \sps \calV_{\tlP}^{\tlQ}$.

En utilisant le lemme \ref{lem:lemcor2G},
on a pour tout $\xi \in \ml_{\tlP}(\rmF) \cap \ol$ fixé 
\begin{equation*}
\int_{N_{Q}(\A)}
\int_{\calV_{\tlP}(\A)}
f_{n \gamma \delta mk}(\xi + Y_{\tlP})dY_{\tlP}dn = 
\int_{\calV_{\tlP}^{\tlQ}(\A)}
\int_{\nl_{\tlQ}(\A)}
f_{\gamma \delta mk}(\xi +  U_{\tlQ}+Y_{\tlP}^{\tlQ})dU_{\tlQ} dY_{\tlP}^{\tlQ}.
\end{equation*}
On voit alors que l'intégrale sur $k \in K$ de (\ref{eq:passAuLeviG}) fois $\eta(\det k)$ devient:
\[
e^{2\rho_{\tlQ}(H_{Q}(a))} \htau_{\tlP}^{\tlQ}(H_{\tlP}(\delta m) - T_{\tlP})
\sum_{\xi \in \ml_{\tlP}(\rmF) \cap \ol}
\sum_{\gamma \in N_{P}^{Q}(\rmF)}
\int_{\calV_{\tlP}^{\tlQ}(\A)}
f_{\tlQ}(\Ad((\gamma \delta m)^{-1})(\xi + Y))dY
\]
où $f_{\tlQ} \in \calS(\ml_{\tlQ}(\A))$ est définie 
par (\ref{eq:fQdefG}) dans le 
paragraphe \ref{par:JTforLevisG}.
D'autre part, l'intégrale sur $A_{\tlQ}^{st,\infty}$ se réduit, 
grâce à l'identité (\ref{eq:HQtlQ}) et le lemme \ref{lem:HbarrhosQ}, à
\[
\int_{\all_{\tlQ}^{st}}
e^{(2\rho_{\tlQ} - 2\rho_{Q} + s \det)(H)}\brtau_{\tlQ}(H + H_{\tlQ}(m)-T_{\tlQ})dH = 
j_{\tlQ}^{-1}
e^{-\upla_{\tlQ,s}(H_{\tlQ}(m))}
\int_{\all_{\tlQ}^{\tlG}}
e^{\upla_{\tlQ,s}(H)}
\brtau_{\tlQ}(H -T_{\tlQ})dH. 
\]
L'intégrale ci-dessus converge pour $s \in \C$ tels que $-1 < \Rel(s) < 1$ 
en vertu du lemme 
\ref{lem:uplasNonNul} \textit{ii)} et donne 
précisément 
$\hat \theta_{\tlQ}(\upla_{\tlQ,s})^{-1} e^{\upla_{\tlQ,s}(T_{\tlQ})}$
où $\hat \theta_{\tlQ} = \hat \theta_{\tlQ}^{\tlG} $ est définie 
par (\ref{eq:thetaHatDefG}) dans le paragraphe \ref{par:fonsPolExpG}.

En utilisant alors l'égalité (\ref{eq:JMQisThisRsG}), 
on s'aperçoit que l'intégrale (\ref{eq:suffitPourQG})
égale
$
J_{\ol}^{M_{\tlQ},T}(\eta_{s}, f_{\tlQ})
$
fois $j_{\tlQ}^{-1}\hat \theta_{\tlQ}(\upla_{\tlQ,s})^{-1} e^{\upla_{\tlQ,s}(T_{\tlQ})}$. 

On a donc la convergence de l'intégrale dans le théorème
ainsi que 
pour tout $T \in \all_{\tlzero}^{+}$ 
suffisamment régulier et tout $s \in \C$ tel que $-1 < \Rel(s) < 1$
\[
\int_{[G]}i_{\ol}(x) \eta_{s}(\det x)dx = 
\sum_{\tlQ \in \relPb}(-1)^{d_{\tlQ}^{\tlG}}
j_{\tlQ}^{-1}\hat \theta_{\tlQ}(\upla_{\tlQ,s})^{-1}e^{\upla_{\tlQ,s}(T_{\tlQ})}I_{\ol}^{M_{\tlQ},T}(\eta_{s}, f_{\tlQ}).
\]
D'après le théorème \ref{thm:mainQualitThm} la somme ci-dessus égale 
$I_{\ol}(\eta_{s}, f)$, ce qu'il fallait démontrer. 
\edem
\erop

\subsection{Résultats d'holomorphie}\label{par:holoResG}

Soient $\{e_{i}^{\vee}\}_{i \in I_{0} \cup -I_{0}} \sbs \all_{I_{0}}$
les vecteurs tels que
\[
\upla_{j}(e_{i}^{\vee}) = 
\begin{cases} 
1 \text{ si } j=i,\\
-1 \text{ si } j= -i,\\
0 \text{ sinon, }
\end{cases} \quad i,j \in I_{0} \cup -I_{0}.
\]
On a donc $e_{i}^{\vee} = -e_{-i}^{\vee}$ pour tout $i \in I_{0}$. En plus, 
il est clair que pour tout $\calI \sbse I_{0}$ l'ensemble $\{e_{i}^{\vee}\}_{i \in \calI}$ 
est une base de $\all_{\acalI}$. 

Soit $\calI \sbse I_{0}$.
On introduit
le cône ouvert $\all_{\calI}^{*}$ défini comme
\[
\all_{\calI}^{*} = 
\displaystyle \left \{
\sum_{i \in \calI} a_{i}\upla_{i}| a_{i} > 0
\right \} 
\sbs \all_{\acalI}^{*} \sbs \all_{\IoB}^{*}.
\]
Donc, pour qu'un $\la \in \all_{\acalI,\C}^{*}$ vérifie
$\Rel(\la) \in \all_{\calI}^{*}$ il faut et suffit que
$\Rel (\la(e_{i}^{\vee})) >0$ pour tout $ i \in \calI$.

Soit $I' \sbs I$. Pour $\la \in \all_{I,\C}^{*}$ on note 
$\la_{I'}$ sa restriction à 
$\all_{I',\C}^{*}$. Remarquons aussi que le caractère $\det$ est naturellement 
un élément de $\all_{I', \C}^{*}$. Il vérifie $\det = \sum_{i \in I'} \upla_{i}$.

On rappelle que nous avons défini $I_{\eta}$ comme l'ensemble de $i \in I_{0}$ tels que la restriction de $\eta$ au groupe 
de normes de $\A_{i}^{*}$ dans $\A^{*}$ est non-trivial.
Pour $I' \sbs I_{0}$ 
notons 
\[
c_{I'} = \int_{T_{I'}(\rmF) \bsl T_{I'}(\A)^{1}}\eta(\det t)dt.
\]
On voit donc que si $I' \cap I_{\eta} \neq \varnothing$ alors $c_{I'} = 0$.
Notons aussi
$v_{I'}$ le volume dans $\all_{I'}$ du 
parallélotope déterminé par les vecteurs 
$\{e_{i}^{\vee}\}_{i \in \calI}$ où 
$\calI \sbse I'$ est tel que $\acalI = I'$. 
Cela ne dépend pas du choix de $\calI$. Notons aussi 
pour tout $i \in I'$:
\begin{equation}\label{eq:hyperplSing}
\calD_{I' ,i, \C}^{*} := \{\la\in \all_{I', \C}^{*} | \la(e_{i}^{\vee}) = 0\} = \{\la\in \all_{I', \C}^{*} | \la(e_{-i}^{\vee}) = 0\}.
\end{equation} 

\blem\label{lem:upetaG} 
Soit $\calI \sbse I_{0}$.  Alors, 
pour tout $x \in G(\A)$, l'intégrale
\[
\int\limits_{\mathclap{T_{\acalI}(\rmF) \bsl T_{\acalI}(\A)}} \ \ 
\ind_{\calI}(H_{I}(tx))e^{\la(H_{I}(tx))}
\eta(\det t)
dt, \quad 
\la \in \all_{\acalI,\C}^{*}
\]
converge absolument et uniformément sur tous les 
compacts de $\Rel(\la) \in \all_{\calI}^{*}$. Elle admet 
un prolongement méromorphe, noté
$\upeta_{\calI}$, égale à
\[
\upeta_{\calI}(\la) = (-1)^{\# \calI}
c_{\acalI}v_{\acalI}\prod_{i \in \calI}\la(e_{i}^{\vee})^{-1}.
\]
En particulier, $\upeta_{\calI}$ ne dépend pas de $x \in G(\A)$ 
et vaut $0$ si $\acalI \cap I_{\eta} \neq \varnothing$. 
Si $\acalI \cap I_{\eta} = \varnothing$ l'ensemble
$\{\calD_{\acalI, i, \C}^{*}\}_{i \in \acalI}$ est l'ensemble de tous les
 hyperplans singuliers de $\upeta_{\calI}$.

\bdem 
Calcul direct. 
\edem
\elem

\blem\label{lem:LamProsG} Soient $\calJ, \calJ_{1}, 
\calJ_{2},\calJ_{3} \sbse I_{0}$ 
tels que 
$\calJ_{1} \sqcup \calJ_{3} \sbs \calJ$ et
$\calJ_{2} \sbs \calJ_{3}$.
\begin{enumerate}[i)]
\item L'intégrale suivante
\[
\Upsilon_{\calJ_{1},\calJ_{2},\calJ_{3}}(\eta, f)(\la) = 
\int\limits_{\mathclap{
T_{I_{0} 
\smin |\calJ_{12}|}(\A) \bsl G(\A)}}  \ 
\ind_{\calJ_{1} \cup \calJ_{2}^{\sharp}}(H_{I}(x))
e^{(\la + \upla_{\calJ_{3 \smin 2}^{\sharp}})(H_{I}(x))}
\hf_{x}^{\calJ_{3}}(X_{\calJ_{1} \cup \calJ_{2}^{\sharp}})\eta(\det x)
dx,
\]
où $\la \in \all_{|\calJ_{12}|,\C}^{*}$, 
converge absolument et uniformément sur tous les compacts de $\Rel(\la) \in \all_{|\calJ_{12}|}^{*}$ 
et définit une fonction holomorphe sur $\all_{|\calJ_{12}|,\C}^{*}$.
\item L'intégrale suivante
\begin{equation*}
\int\limits_{\mathclap{
T_{|\calJ \smin \calJ_{12}|}(\rmF)
T_{I_{0} \smin \acalJ}(\A)\bsl G(\A)}}
\ind_{(\calJ \smin \calJ_{2}) \cup \calJ_{2}^{\sharp}}(H_{I}(x))
e^{\la(H_{I}(x))}
\hf_{x}^{\calJ_{3}}(X_{\calJ_{1} \cup \calJ_{2}^{\sharp}})\eta(\det x)dx, \quad 
\la \in \all_{\acalJ, \C}^{*}
\end{equation*}
converge absolument et uniformément sur tous les compacts de
\begin{equation*}
\Rel(\la) \in 
\all_{|\calJ_{12}|}^{*} \times 
(\upla_{\calJ_{3 \smin 2}^{\sharp}} + 
\all_{\calJ \smin \calJ_{12}}^{*})
\end{equation*}
et elle admet un prolongement méromorphe
à $\all_{\acalJ,\C}^{*}$, noté
$\brLa_{\calJ_{1},\calJ_{2},\calJ_{3}}^{\calJ}(\eta, f)$
qui vérifie
\begin{equation*}
\brLa_{\calJ_{1},\calJ_{2},\calJ_{3}}^{\calJ}(\eta, f)(\la) = 
\upeta_{\calJ \smin \calJ_{12}}
(\la_{|\calJ \smin \calJ_{12}|} + \upla_{\calJ_{3 \smin 2}})
\Upsilon_{\calJ_{1},\calJ_{2},\calJ_{3}}(\eta, f)(\la_{|\calJ_{12}|}).
\end{equation*}
Si $|\calJ \smin \calJ_{12}| \cap I_{\eta} \neq \varnothing$ 
on a $\brLa_{\calJ_{1},\calJ_{2},\calJ_{3}}^{\calJ}(\eta, f) \equiv 0$ 
et dans le cas contraire, l'ensemble:
\[
\{\calD_{\acalJ, i, \C}^{*}\}_{i \in |\calJ \smin \calJ_{13}|} 
\cup \{\calD_{\acalJ, |i|, \C}^{*} + \upla_{-i}\}_{i \in \calJ_{3 \smin 2}}, 
\]
où les hyperplans $\calD_{\acalJ, i, \C}^{*}$ 
sont définies par (\ref{eq:hyperplSing}), 
est l'ensemble de tous les hyperplans singuliers de $\brLa_{\calJ_{1},\calJ_{2},\calJ_{3}}^{\calJ}(\eta, f)$.
\end{enumerate}

\bdem 
Le premier point c'est le lemme 5.22 de \cite{leMoi}.
Les assertions sur la convergence, le prolongement méromorphe
ainsi que l'équation fonctionelle 
dans le deuxième point 
découlent des lemmes 5.24 et 5.25 de loc. cit. Finalement, la propriété d'annulation 
de $\brLa_{\calJ_{1},\calJ_{2},\calJ_{3}}^{\calJ}(\eta, f)$ 
et l'assertion sur ses éventuels hyperplans singuliers
découlent de l'équation
fonctionnelle qu'elle vérifie et de la propriété analogue 
de la fonction $\upeta_{\calJ \smin \calJ_{12}}$ décrite dans le 
lemme \ref{lem:upetaG} ci-dessus.
\edem
\elem

\brop\label{prop:zetaDefPropsG} 
Soit $\calJ \sbse I_{0}$.
Alors, l'intégrale 
\[
\int\limits_{\mathclap{
G(\A, X_{\calJ})\bsl G(\A)}}
f(x^{-1}X_{\calJ}x)
e^{\la(H_{I}(x))}\eta(\det x)dx, \quad \la \in \all_{\acalJ,\C}^{*}
\]
converge absolument 
et uniformément sur tous les compacts de
$\Rel(\la) \in \upla_{\calJ^{\sharp}} + \all_{\calJ^{\sharp}}^{*}$
et admet un prolongement méromorphe 
à $\all_{\acalJ,\C}^{*}$, noté $\zeta_{\calJ}(\eta, f)$,
qui vérifie
\[
\zeta_{\calJ}(\eta, f) = 
\sum_{\calJ_{1} \sqcup \calJ_{2} \sqcup \calJ_{3} \sbs \calJ}
(-1)^{\# \calJ_{3}}
\brLa_{\calJ_{1},\calJ_{2},\calJ_{23}}^{\calJ}(\eta,f).
\]
En plus, l'ensemble
\[
\{\calD_{\acalJ, i, \C}^{*}\}_{i \in \acalJ} 
\cup \{\calD_{\acalJ, |i|, \C}^{*} + \upla_{-i}\}_{i \in \calJ}, 
\]
où les hyperplans $\calD_{\acalJ, i, \C}^{*}$ 
sont définies par (\ref{eq:hyperplSing}), contient tous les 
hyperplans singuliers de $\zeta_{\calJ}(\eta, f)$.
\bdem 
Convergence, prolongement méromorphe et l'équation fonctionnelle 
découlent de la proposition 5.26 accouplée 
avec le lemme 5.25 de \cite{leMoi}. L'assertion sur 
les hyperplans singuliers est une conséquence de l'équation fonctionnelle 
vérifiée par $\zeta_{\calJ}(\eta, f)$ et du lemme \ref{lem:LamProsG} ci-dessus.
 \edem
 \erop
 
 Dans le contexte du lemme \ref{lem:LamProsG} ci-dessus, 
 posons $\brLa_{\calJ_{1},\calJ_{2}}^{\calJ}(\eta, f) = \brLa_{\calJ_{1},\calJ_{2},\calJ \smin \calJ_{1}}^{\calJ}(\eta, f)$.
 
 \blem\label{lem:sumZetasUpsilonG}
 Soit $I' \sbs I$ et $\calJ_{0} \sbse I' \cap I_{\eta}$ 
 tel que $|\calJ_{0}| =  I' \cap I_{\eta}$.
On a alors l'égalité de fonctions méromorphes
sur $\all_{I',\C}^{*}$ suivante
 \[
\sum_{|\calI_{0}| = I' \smin |I' \cap I_{\eta}|}\zeta_{\calJ_{0} \cup \calI_{0}}(f)(\eta, \la) =  
\sum_{\calJ_{0} \sbs \calJ_{1} \sqcup \calJ_{2} \sbse I'}
\sum_{|\calJ_{3}| = I' \smin |\calJ_{12}|}
(-1)^{\#\calJ_{3}}
\brLa_{\calJ_{1}, \calJ_{2}}^{\calJ_{123}}(f)(\eta, \la).
 \] 
 \bdem 
 Soient $\calJ_{1}, \calJ_{2}, \calJ_{3}, \calI_{0} \sbse I' \smin |I' \cap I_{\eta}|$ tels que $|\calI_{0}| =  I' \smin |I' \cap I_{\eta}|$
 et $\calJ_{1} \sqcup \calJ_{2} \sqcup \calJ_{3} \sbs \calI_{0}$.  
 En vertu du lemme \ref{lem:LamProsG} \textit{ii)} on a que si 
 $\calJ_{0} \nsbs \calJ_{1} \sqcup \calJ_{2}$ alors
 $\brLa_{\calJ_{1},\calJ_{2},\calJ_{23}}^{\calJ_{0} \cup \calI_{0}}(\eta,f) \equiv 0$.
 En vertu de la proposition \ref{prop:zetaDefPropsG} on a donc:
 \begin{multline*}
 \sum_{\mathclap{|\calI_{0}| = I' \smin |I' \cap I_{\eta}|}}
 \zeta_{\calJ_{0} \cup \calI_{0}}(f)(\eta, \la) = 
 \sum_{\mathclap{|\calI_{0}| = I' \smin |I' \cap I_{\eta}|}}
 \quad \quad
  \sum_{\calJ_{0} \sbs \calJ_{1} \sqcup \calJ_{2} \sbs \calJ_{0} \cup \calI_{0}}
  \sum_{\calJ_{3} \sbs (\calJ_{0} \cup \calI_{0}) \smin \calJ_{12}}
(-1)^{\# \calJ_{3}}
\brLa_{\calJ_{1},\calJ_{2},\calJ_{23}}^{\calJ_{0} \cup \calI_{0}}(\eta,f) = \\ 
  \sum_{\calJ_{0} \sbs \calJ_{1} \sqcup \calJ_{2} \sbse I'}
  \sum_{\calJ_{3} \sbse I' \smin |\calJ_{12}|}
  (-1)^{\# \calJ_{3}}
  \sum_{\calJ_{4} \sbse I' \smin |\calJ_{123}|}
\brLa_{\calJ_{1},\calJ_{2},\calJ_{23}}^{\calJ_{1234}}(\eta,f).
 \end{multline*}
Il résulte de l'égalité 
(5.27)
du lemme 5.27 de \cite{leMoi} que si $I' \smin |\calJ_{123}| \neq \varnothing$ 
 on a
 \[
   \sum_{\calJ_{4} \sbse I' \smin |\calJ_{123}|}
\brLa_{\calJ_{1},\calJ_{2},\calJ_{23}}^{\calJ_{1234}}(\eta,f) \equiv 0.
 \]
 D'où le résultat voulu.
\edem
 \elem
 
\blem\label{lem:detIsNotSing}
 Pour tout $\calI \sbse I_{0}$ non-vide et tout $i \in \calI$ 
la droite $\C \det$ n'est pas contenu 
dans $\calD_{\acalJ, |i| ,\C}^{*} \cup (\calD_{\acalJ, |i| ,\C}^{*} + \upla_{-i})$.
\bdem
Il suffit de vérifier qu'il existe un $s \in \C$ tel que 
$s\det (e_{|i|}^{\vee}) \neq 0$ et $(s\det + \upla_{i}) (e_{|i|}^{\vee}) \neq 0$. Mais $s\det (e_{|i|}^{\vee}) = s$
et $(s\det + \upla_{i}) (e_{|i|}^{\vee}) = s + \frac{i}{|i|}$, d'où le résultat.
\edem
\elem
 
 Soient $\calJ, \calJ_{1}, \calJ_{2} \sbse I_{0}$ 
tels que $\calJ_{1} \sqcup \calJ_{2} \sbs \calJ$. Posons:
 \begin{equation}\label{eq:uplaDet}
\brLa_{\calJ_{1},\calJ_{2}}^{\calJ}(\eta, f)(s) = 
\brLa_{\calJ_{1},\calJ_{2}}^{\calJ}(\eta, f)(s \det ), \quad s \in \C.
\end{equation}
En vertu des lemmes \ref{lem:LamProsG} et \ref{lem:detIsNotSing}, la 
fonction $\brLa_{\calJ_{1},\calJ_{2}}^{\calJ}(\eta, f)(s)$ de la variable $s \in \C$ est méromorphe 
sur $\C$. On a alors le corollaire du lemme \ref{lem:LamProsG} suivant:
 \bcor\label{cor:UpsPros}
 Soient $\calJ, \calJ_{1}, \calJ_{2} \sbse I_{0}$ 
tels que $\calJ_{1} \sqcup \calJ_{2} \sbs \calJ$.
\begin{enumerate}[i)]
\item
Pour $-1 < \Rel(s) < 1$, la fonction $\brLa_{\calJ_{1},\calJ_{2}}^{\calJ}(\eta, f)(s)$ est 
donnée par l'intégrale absolument convergente suivante:
\begin{equation*}
\int\limits_{\mathclap{
T_{I_{0} \smin |\calJ_{12}|}(\rmF)
\bsl G(\A)}}
\ind_{(\calJ \smin \calJ_{2}) \cup \calJ_{2}^{\sharp}}(H_{I}(x))
\hf_{x}^{\calJ \smin \calJ_{1}}(X_{\calJ_{1} \cup \calJ_{2}^{\sharp}})\eta_{s}(\det x)dx.
\end{equation*}
\item La fonction $\brLa_{\calJ_{1},\calJ_{2}}^{\calJ}(\eta, f)(s)$ est holomorphe sur $s \in \C \smin \{-1,1\}$. 
\item Si $|\calJ \smin \calJ_{12}| \cap I_{\eta} \neq \varnothing$ 
on a $\brLa_{\calJ_{1},\calJ_{2}}^{\calJ}(\eta, f)(s) \equiv 0$.
\end{enumerate}
\bdem 
Pour démontrer le point \textit{i)}, en vertu du lemme \ref{lem:LamProsG} il suffit de vérifier 
\begin{equation}\label{eq:detIsConv}
\Rel(s \det) \in 
\all_{|\calJ_{12}|}^{*} \times 
(\upla_{(\calJ \smin \calJ_{12})^{\sharp}} + 
\all_{\calJ \smin \calJ_{12}}^{*}) \quad -1 < \Rel(s) < 1.
\end{equation}
Autrement dit, on doit avoir $\Rel(s) \det - \upla_{(\calJ \smin \calJ_{12})^{\sharp}} \in \all_{\calJ \smin \calJ_{12}}^{*}$.
Or, on a que la restriction de $\det$ à $\all_{|\calJ \smin \calJ_{12}|}$ égale
$\sum_{i \in \calJ \smin \calJ_{12}} \upla_{|i|} = \sum_{i \in \calJ \smin \calJ_{12}} \frac{i}{|i|}\upla_{i}$ 
et $\upla_{(\calJ \smin \calJ_{12})^{\sharp}} = \sum_{i \in \calJ \smin \calJ_{12}} -\upla_{i}$. 
Par définition de $\all_{\calJ \smin \calJ_{12}}^{*}$ on a donc que (\ref{eq:detIsConv}) est vérifié si et seulement 
si pour tout $i \in \calJ \smin \calJ_{12}$ on a $\Rel(s)\frac{i}{|i|} > -1$. Comme 
$\frac{i}{|i|} \in \{-1,1\}$, l'assertion suit.

Pour le point \textit{ii)}, en utilisant l'équation fonctionnelle donnée dans 
le lemme \ref{lem:LamProsG}, on a
\begin{equation}\label{eq:EQfoncAux}
\brLa_{\calJ_{1},\calJ_{2}}^{\calJ}(\eta, f)(s) = 
\upeta_{\calJ \smin \calJ_{12}}
(s \det+ \upla_{\calJ \smin \calJ_{12}})
\Upsilon_{\calJ_{1},\calJ_{2},\calJ_{3}}(\eta, f)(s \det)
\end{equation}
et la description explicite de $\upeta_{\calJ \smin \calJ_{12}}$ donnée dans le lemme \ref{lem:upetaG} nous donne
\[
\upeta_{\calJ \smin \calJ_{12}}
(s \det+ \upla_{\calJ \smin \calJ_{12}}) = 
 (-1)^{\# (\calJ \smin \calJ_{12})}
c_{|\calJ \smin \calJ_{12}|}v_{|\calJ \smin \calJ_{12}|}\prod_{j \in \calJ \smin \calJ_{12}}
\dsl 1 + \dfrac{i}{|i|}s\rb^{-1}.
\]
Puisque $1 + \frac{i}{|i|}s = 0$ implique $s \in \{-1,1\}$ on a l'holomorphie de $\brLa_{\calJ_{1},\calJ_{2}}^{\calJ}(\eta, f)(s)$ 
sur $\C \smin \{-1,1\}$, vu l'holomorphie 
de $\C \ni s \mapsto \Upsilon_{\calJ_{1},\calJ_{2},\calJ_{3}}(\eta, f)(s \det )$.

Le point \textit{iii)} découle de l'équation fonctionnelle (\ref{eq:EQfoncAux}) ci-dessus.
\edem
\ecor

Pour $\calJ \sbse I$ on définit la fonction méromorphe $\zeta_{\calJ}(\eta, f)$ sur $\C$ 
par
\[
\zeta_{\calJ}(\eta, f)(s) = \zeta_{\calJ}(\eta, f)(s \det), \quad s \in \C.
\]
En vertu de la proposition \ref{prop:zetaDefPropsG} et du lemme \ref{lem:detIsNotSing} 
la définition est licite.

On a alors le corollaire suivant du lemme \ref{lem:sumZetasUpsilonG} et du corollaire 
\ref{cor:UpsPros}.

 \bcor\label{cor:zetasCsum}
 Soit $I' \sbs I$ et $\calJ_{0} \sbse I' \cap I_{\eta}$ 
 tel que $|\calJ_{0}| =  I' \cap I_{\eta}$.
On a alors l'égalité de fonctions méromorphes
sur $\C$ suivante
 \[
\sum_{|\calI_{0}| = I' \smin |I' \cap I_{\eta}|}\zeta_{\calJ_{0} \cup \calI_{0}}(\eta, f) =  
\sum_{\calJ_{0} \sbs \calJ_{1} \sqcup \calJ_{2} \sbse I'}
\sum_{|\calJ_{3}| = I' \smin |\calJ_{12}|}
(-1)^{\#\calJ_{3}}
\brLa_{\calJ_{1}, \calJ_{2}}^{\calJ_{123}}(\eta, f).
 \] 
 En particulier, en vertu du point \textit{ii)} du corollaire \ref{cor:UpsPros}, 
 la somme ci-dessus est holomorphe sur $\C \smin \{-1,1\}$.
 \ecor